\documentclass[usletter,12pt]{article}
\usepackage{amsmath}
\numberwithin{equation}{section}
 \usepackage{color}
\usepackage{graphicx}
\usepackage{appendix}
\usepackage{enumerate}
\usepackage{natbib}
\usepackage{url} 
\usepackage{amssymb}
\usepackage{amsthm}
\usepackage{subcaption}

\usepackage{cancel}
\usepackage{multirow}
\usepackage[colorlinks=true,linkcolor=black,anchorcolor=black,citecolor=black,filecolor=black,menucolor=black,runcolor=black,urlcolor=black]{hyperref}
\usepackage{placeins}

\theoremstyle{definition}
\newtheorem{theorem}{Theorem}
\newtheorem{corollary}{Corollary}
\newtheorem{proposition}{Proposition}

\newtheorem{lemma}{Lemma}

\newtheorem{remark}{Remark}[section]

\usepackage[margin=1in]{geometry}
\usepackage{setspace}
\DeclareCaptionFormat{custom}
{
    \textbf{#1#2}\textit{\small #3}
}
\usepackage{titlesec}
\titlespacing*{\section}{3pt}{0pt}{5pt}
\titlespacing*{\subsection}{3pt}{0pt}{5pt}
\titlespacing*{\subsubsection}{3pt}{0pt}{5pt}

\usepackage{needspace}
\usepackage{titling}
\newcommand{\labthis}{\stepcounter{equation}\tag{\theequation}}

\usepackage{comment}
\usepackage{authblk}
\usepackage{xcolor}
\definecolor{dpink}{rgb}{1,0,1}

\definecolor{hypcol}{rgb}{0.5,0,1}
\hypersetup{colorlinks=true,linkcolor=hypcol,anchorcolor=hypcol,citecolor=hypcol,filecolor=hypcol,urlcolor=hypcol,bookmarksnumbered=true}

\title{Reconciliation of Bayes and empirical Bayes interval estimation with application to small area estimation }

\author[1]{Aditi Sen}
\author[2]{Masayo Y. Hirose}
\author[1]{Partha Lahiri}

\affil[1]{Department of Mathematics, University of Maryland, College Park, USA.}
\affil[2]{Institute of Mathematics for Industry, Kyushu University, Japan.} 

\date{}

\begin{document}
\maketitle
\vspace{-2.5em}

\begin{abstract}

Multi-level normal hierarchical models, also interpreted as mixed effects models, play an important role in developing statistical theory in multi-parameter estimation for a wide range of applications, including small area estimation. In this article, we propose a novel reconciliation framework of the empirical Bayes (EB) and hierarchical Bayes approaches for interval estimation of random effects under a two-level normal model. Our framework shows that a second-order efficient empirical Bayes confidence interval, with  EB coverage error of order $O(m^{-3/2})$, $m$ being the number of areas in the area-level  model, can also be viewed as a credible interval whose expected posterior coverage is close to the nominal level, provided a carefully chosen prior, referred to as a ‘matching prior’, is placed on the hyperparameters. While existing literature has examined matching priors that reconcile frequentist and Bayesian inference in various settings, this paper is the first to study matching priors with the goal of interval estimation of random effects in a two-level model. We obtain a matching prior on the variance component that achieves a proper posterior under mild regularity conditions. The theoretical results in the paper are corroborated through a Monte Carlo simulation study and real data analysis.

\end{abstract}

\noindent%

\section{Introduction}

Multi-level modeling provides a flexible framework that incorporates uncertainties from different sources while combining information using multiple databases. The problem of simultaneously estimating several independent normal means received intense attention in the decades following the introduction of the James–Stein estimator in \citet{james1961estimation}. Since then, multi-level models have been extensively developed from both theoretical and applied perspectives with applications spanning diverse domains -- from  compound decision problems to meta-analysis. Relevant references include \citet{MorrisJASA,  casella1983empirical, EversonMorris2000, carlin1997bayes, Hwang2008, gelman2007data, efron2012large, dersimonian1986meta, raudenbush2002hierarchical}, among others. A particularly important application of multi-level modeling is in Small Area Estimation (SAE), where estimates are obtained at granular levels to facilitate important policy making decisions -- see \citet{Rao2015, Ghosh2020, Franco, otto1995sampling, DATTA2009251}, among others, for related applications. 

In this article, we assume that the data follows a two-level normal model - the $1^\text{st}$ level considers input data at the area level and brings in sampling variability, while the $2^\text{nd}$ level borrows strength from other areas connecting area-level means to covariates. This model is commonly referred to as the Fay–Herriot model in SAE literature -- \citet{fay1979estimates} proposed it in the context of estimating per-capita income for small places, using housing data and tax records as supplementary information. This two-level normal model is presented below. For $i=1,\cdots,m,$
\begin{eqnarray}
\mbox{Level 1:}\;&&y_i|\theta_i \stackrel{\text{ind}}\sim N(\theta_i,D_i), 
\label{model_1} \\ 
\mbox{Level 2:}\;&&\theta_i|\beta,A \stackrel{\text{ind}}\sim N(x_i'\beta,A), \label{model_2}
\end{eqnarray}
where $m$ is the total number of small areas, $y_i$ is the direct survey estimate (usually weighted averages of the sample observations in $i^\text{th}$ area) of the true small area mean $\theta_i$, $D_i$ is the sampling variance of $y_i$, $A$ is the prior variance, $x_i$ and $\beta$ are $p\times 1$ vectors of known auxiliary variables from area-level supplementary data and their respective regression coefficients. In this model we assume $\beta \in \mathbb{R}^p, $ the $p$-dimensional Euclidean space, and $p < m$. Following the Bayesian framework of \citet{fay1979estimates}, \eqref{model_1} is referred to as the sampling distribution where variances ($D_i$'s) are assumed to be known (although in practice, they are estimated using within area observations) and \eqref{model_2} is referred to as the prior distribution, where hyperparameters $\beta$ and $A$ are assumed to be unknown. 

While the present work focuses on an SAE context, this two-level normal model has broad applicability beyond this domain. \citet{efron1973stein} studied a special case of this model without covariates, with the aim of providing an EB justification of the James–Stein estimator using the prior $\theta_i \stackrel{\text{iid}}\sim N(0,A),\; i=1,\cdots,m$. 
Variants of this model including covariates have been used in the estimation of poverty rates and proportions at the lowest level of literacy for states and counties -- see \citet{kalton2000small, mohadjer2012hierarchical}, among others. In recent times, this model finds applications in high-dimensional data -- \citet{hwang2013empirical} studied empirical Bayes confidence interval (EBCI) for selected parameters in this context. In fact, in a broader sense, the model given by (\ref{model_1}) and (\ref{model_2}) can be viewed as a linear mixed-effects model, including some distinct area-level variation. In this representation of the model, we consider $i^\text{th}$ area-level mean $\theta_i$ as $x_i'\beta + v_i$, where $v_i$ is the random area-level effect and $e_i$ is the sampling error associated with $y_i$. Thus, we have 
$$y_i = x_i'\beta + v_i + e_i,$$ where $\{v_i,\;i=1,\cdots,m\}$ and  $\{e_i,\;i=1,\cdots,m\}$ are independent with $v_i \stackrel{\text{iid}}\sim N(0,A)$ and $e_i \stackrel{\text{ind}}\sim N(0,D_i), \; i = 1, \cdots, m$.

For the aforementioned model, the empirical best linear unbiased prediction (EBLUP) in the classical linear mixed model context is identical to the empirical Bayes (EB) prediction approach of $\theta_i$ in the two-level normal model. EBLUP of $\theta_i$ is obtained by replacing the unknown 
$A$ with its estimator in the best linear unbiased predictor (BLUP) that minimizes the mean squared error (MSE) among all linear unbiased predictors. Considering $A$ to be known but $\beta$ unknown, the formula of BLUP of $\theta_i$ for this model is $$\hat{\theta}^{\text{BLUP}}_i = (1-B_i) y_i + B_i x_i'\Bar{\beta},$$ where $B_i \equiv B_i(A) = D_i/(A+D_i)$ is the shrinkage factor, $\Bar{\beta} \equiv \Bar{\beta}(A) = 
(X'V^{-1}X)^{-1}X'V^{-1}y$ is the weighted least square estimator of $\beta$ with $V \equiv \mbox{diag}(A+D_1,\dots,A+D_m)$, an $m \times m$ diagonal matrix and $X$, an $m\times p$ known matrix of auxiliary data and $y = (y_1,\cdots,y_m)'$, an $m\times1$ vector of direct estimates. In this case, the predictor error variance is
\begin{align}\label{bp_dist}
    \mbox{Var}(\theta_i - \hat\theta_i^{\text{BLUP}}) = g_{1i} + g_{2i}, \; \; i =1,\cdots,m,
\end{align}
where $ g_{1i} \equiv g_{1i}(A) = A D_i /(A+D_i)$ and $g_{2i} \equiv g_{2i}(A) = B_i^2 \; r_i$ with $r_i = x_i' (X'V^{-1}X)^{-1} x_i$. An EBLUP of $\theta_i$ is then obtained by using consistent estimators of $A$, such as Maximum Likelihood (ML), Residual Maximum Likelihood (REML), Adjusted REML (AML), etc. We refer the readers to \citet{PR1990, ghosh1994small, jiang1997derivation, jiang1998asymptotic, datta2000unified, das2004mse, jiang2006mixed, Jiang2007}, among others for advancements of linear mixed effects models and small area applications.

As an alternative to the aforementioned EBLUP or the EB procedure, the hierarchical Bayes (HB) procedure uses the posterior mean and posterior standard deviation for estimating parameter of interest ($\theta$) and associated the standard error. We refer the readers to \cite{DattaGhosh1991, Datta01121999, datta2011bayesian} for HB methods in SAE applications. The HB analysis of the Fay-Herriot model considers a $3^\text{rd}$ level that assumes prior distributions on the parameters from the $2^\text{nd}$ level (i.e., hyper-parameters). In this article, we focus on a dual justification -- both EB and HB, for interval estimation of random effects in the two-level normal model. We emphasize that prior interval estimation studies in SAE have primarily adopted the parametric bootstrap (see \citet{Hall2006, CLL2008}, among others) and EBCI approach (see \citet{datta2002, YoshimoriLahiri2014, Hirose2017}, among others). In this article, we explore second-order correct EBCIs from a new perspective (HB), with the goal of achieving expected posterior coverage close to the nominal level. Our results indicate that a valid HB justification for an EBCI can be achieved only when an appropriate matching prior is placed on the hyperparameter. Matching priors bridging Bayesian and frequentist paradigms have been studied extensively in literature -- see the book article \citet{datta2004book}. \citet{DMGS2000} considered a sequence of independent and identically distributed random variables and obtained probability matching priors for prediction region of a future observation. While matching priors in two-level normal models have been studied previously by \citet{DGM2000, DattaRaoSmith2005, ganesh2008new, HiroseLahiri2021}, the emphasis in those works differs from ours. For example, \citet{DattaRaoSmith2005} aimed to reconcile measures of variability, i.e., expected posterior variance with the MSE. In the next section, we review existing confidence intervals (CI) and matching priors in this model and describe our contribution in detail.

\section{Existing literature}

We dedicate this section to a review of existing CIs and matching priors for the two-level normal model, followed by our contributions.

\subsection{Interval estimation of small area means}
\label{sec:int}

In this section, we discuss in detail the construction of widely used Cox-type EBCIs. To explain its origin, we begin by discussing the direct CI. For a confidence level $\alpha \in (0,1)$, the direct CI of $\theta_i$ is given by $$I_i^\text{D}: y_i \pm z \sqrt{D_i},$$ 
where $z \equiv z_{\alpha/2}$ is the $100(1-\alpha/2)\%$ point of the standard normal distribution $N(0,1)$. Since this interval can have an unacceptably large length for large values of $D_i$'s, \citet{Cox_1975} introduced an EBCI of $\theta_i$, for a special case of the two-level normal model (for $p=1$, $x_i\beta=\beta$ and $D_i=D, \; i = 1,\cdots,m$)  
given by $$I_i^{\text{Cox}} (\bar y,\hat{A}_\text{ANOVA}): \hat{\theta}^\text{EB}_i(\bar y,\hat{A}_{\text{ANOVA}})  \pm z \sigma(\hat{A}_{\text{ANOVA}}),$$ 
where 
$\sigma \equiv \sigma(A) = \sqrt{AD/(A+D)}$, \; $\hat\theta_i^\text{EB}(\bar y, \hat A_{\text{ANOVA}}) =(1- \hat B)y_i + \hat B \bar y$ is an EB estimator of $\theta_i$  and $\hat{A}_{\text{ANOVA}} = \text{max} \lbrace (m-1)^{-1} \sum_{i=1}^{m} (y_i - \bar y)^2, 0 \rbrace.$ However, the EB coverage error of $I_i^{\text{Cox}}$ being of the order $O(m^{-1})$ could be large in an SAE problem, as the number of areas ($m$) to be pooled is not large. Consequently, \citet{datta2002, YoshimoriLahiri2014} introduced new EBCIs of $\theta_i$ reducing the EB coverage error to $o(m^{-1})$. Especially, the latter one obtained an EBCI that is always of smaller length than $I_i^\text{D}$ with the help of a carefully devised AML estimator of $A$ given by $\hat{A}_{i;\text{gls}}$ or $\hat{A}_{i;\text{ols}}$, depending on the generalized least square (GLS) or ordinary least square (OLS) estimator for $\beta$. For brevity, hereon we use $\hat{A}_{i;\text{gls}}$ and denote this interval as
$$I_{i}^\text{YL} \equiv I_{i}^{\text{Cox}}(\hat \beta, \hat{A}_{i;\text{gls}}): \hat{\theta}^\text{EB}_i(\hat \beta, \hat{A}_{i;\text{gls}}) \pm z \sigma_i (\hat{A}_{i;\text{gls}}),$$ 
where $\sigma_i \equiv \sigma_i(A) = \sqrt{g_{1i}},$ and $\hat\theta_i^\text{EB}(\hat \beta, \hat A_{i;\text{gls}}) =(1- \hat B_i)y_i + \hat B_i x_i' \hat \beta$ is an EB estimator of $\theta_i$ with $\hat{\beta} \equiv 
\Bar {\beta}(\hat A_{i;\text{gls}}) = 
(X'\hat V^{-1}X)^{-1}X'\hat V^{-1}y$. 
Following this, with a further objective of allowing high leverage and reducing the number of iterations required to estimate $A$ (for large $m$), \citet{Hirose2017} proposed an EBCI (see Corollary 1(b)), which we denote as  
$$I_i^\text{N}: \hat{\theta}_i^\text{EB}(\Tilde \beta, \Tilde{A}_i) \pm z \delta_i (\Tilde{A}_i),$$ 
where $\hat{\theta}_i^\text{EB}(\Tilde \beta, \Tilde{A}_i)$ is an EB estimator of $\theta_i$ using $\Tilde{A}_i$ -- the proposed AML estimator of $A$, as one of new EBCIs. Here, $\delta_i \equiv \delta_i(A) = \sqrt{g_{1i} + g_{2i}}$ serves as the standard error factor, arising from $\delta_i^2$ being the MSE of BLUP of $\theta_i$ (cf. \eqref{bp_dist}), in contrast to $\sigma_i$ used in $I_i^{\text{YL}}$. Consequently, $I_i^\text{N}$ is not a Cox-type interval. Since, both $I_{i}^{\text{YL}}$ and $I_{i}^\text{N}$ achieve second-order correct EB coverage, under the model  given by \eqref{model_1} and \eqref{model_2} we have 
\begin{align}
    P(\theta_i \in I_{i}|\beta,A) = 1 - \alpha + O(m^{-3/2}), \quad i=1,\cdots m, \label{EB cov}
\end{align} 
where $I_i \in \{ I_{i}^{\text{YL}}, I_{i}^\text{N} \}$ and the probabilities are with respect to the joint distribution of $(\theta,y)$, conditional on the hyperparameters $(\beta, A)$. This article introduces a Bayesian framework that constructs matching priors for EBCIs, providing a dual justification for their use, with specific implementation for the aforementioned EBCIs. While the comparison of these intervals is not our main focus, Proposition \ref{cor5} in Section \ref{sec:asym} yields an interesting finding that upto $O_p(m^{-1})$, $I_{i}^\text{N}$ and $I_{i}^{\text{YL}}$ are equivalent in terms of length, and the difference in lengths is of the smaller order $O_p(m^{-3/2})$.

\subsection{Matching priors in SAE}
\label{sec:mat}
In this section, we review existing literature involving matching priors in the SAE context. \cite{DGM2000} developed a quantile matching prior in the balanced common mean case of the two-level normal model given by
\begin{eqnarray}
\mbox{Level 1:}\;&&y_i|\theta_i \stackrel{\text{ind}}\sim N(\theta_i,D),\quad D \text{ assumed known}, \label{DGM1} \\ 
\mbox{Level 2:}\;&&\theta_i|\mu,A \stackrel{\text{iid}}\sim N(\mu,A), \quad i =1,\cdots,m, \label{DGM2}\\
\mbox{Level 3:}\;&&\pi(\mu,A) \propto \pi(A), \label{DGM3}
\end{eqnarray}
where $D$ is the common sampling variance in Level 1 and $\mu$ is the common mean in Level 2.
Assuming a uniform prior on $\mu$ in Level 3, the authors derived a prior $\pi(A)$ via posterior and frequentist quantile matching of $\theta_i$. In other words, the authors considered a one-sided credible interval in a special case and their objective was to achieve nominal frequentist coverage (under Levels 1 and 2) of posterior quantile (calculated under Levels 1, 2 and 3 above). This prior is given by
\begin{equation}\label{dgmprior}
   \pi_{\text{DGM}}(A) \propto \frac{A}{A + D}, 
\end{equation}   
which provides a proper posterior for $m>3$.

Later \citet{DattaRaoSmith2005} considered the Bayesian implementation of the general two-level model, given by \eqref{model_1} and \eqref{model_2}, that assumes a flat prior on $\beta$ and a prior $\pi(A)$ on $A$ in the $3^{\text{rd}}$ level as follows: 
\begin{align} \label{model_3}
\mbox{Level 3:}\; (\beta,A) \sim p(\beta,A) \propto \pi(A), \quad \beta \in \mathbb{R}^p, \; A >0.
\end{align}
The authors provided a frequentist justification of the posterior variance of $\theta_i$, denoted as $V_i^\text{HB} \equiv \mbox{Var}(\theta_i|y)$, by deriving a prior on $A$ with the matching criteria that the expected value of $V_i^\text{HB}$ is a second-order unbiased estimator for the MSE of EBLUP of $\theta_i$ using REML estimator $\hat{A}$, for large $m$. This matching criteria is mathematically given by
$$E(V_i^\text{HB}) = \mbox{MSE}[\hat{\theta}_i^\text{EB}(\hat{A})] + o(m^{-1}).$$ The prior obtained in \citet{DattaRaoSmith2005} is area-specific -- a function of $D_i$'s, which, in turn, depend on the area-specific sample sizes. This prior is given by
\begin{equation}\label{drsprior}
   \pi_{i;\text{DRS}}(A) \propto (A+D_i)^2 \sum_{u=1}^m (A+D_u)^{-2}, 
\end{equation}   
and provides a proper posterior for $m>p+2$.
\citet{HiroseLahiri2021} also obtained the same prior as in \eqref{drsprior} with different matching criteria involving the shrinkage parameter $B_i$ and EBLUP. The authors referred to this as the `multi-goal' prior and used an AML estimator of $A$ having several advantages over REML, proposed in \citet{HiroseLahiri2018}. For further technical details on the AML estimator and matching criteria, 
we direct the interested readers to \citet{HiroseLahiri2018, HiroseLahiri2021}.

It is important to note here that since \cite{DGM2000} considered a balanced case of the two-level model, their prior $\pi_{\text{DGM}}$ does not depend on area $i$. Whereas in general two-level models, area-specific matching priors, varying with index $i$, are widely accepted (see \cite{DattaRaoSmith2005, HiroseLahiri2021}, among others). Originally defined probability matching priors by \cite{welch1963formulae} (we refer to as WP) depend on the Fisher information matrix, which in turn depends on the model, and critically, on the parameter of interest. If the interest parameter ($\theta_i$ in this case) is changed for the same model, different WP probability matching priors are obtained. Thus, the matching motivation renders the inference to be unlike regular Bayesian or HB inference, where for a given problem, inference results are usually generated for one prior.

\subsection{Our contribution}

In this section, we address the reconciliation of EB and HB approaches for interval estimation of random effects using matching priors. Despite continued interest in matching priors, to the best of our knowledge, none of the existing literature has explored matching priors for interval estimation of random effects in the two-level normal model. The closest to our work is that of \cite{DGM2000}, discussed in the previous section. However, they consider a special case of the two-level model with a different matching criterion, motivation and proof techniques. For comparison purposes, their quantile matching criterion can be interpreted as a one-sided credible interval.
In this article, we demonstrate that an EBCI $I_i$ admits a dual justification: in addition to achieving the second-order correct EB coverage, its { expected} posterior coverage under the HB framework also approximates the nominal level when an appropriate prior is placed on the hyperparameter. Specifically, given the data $y$ and the model described by \eqref{model_1}, \eqref{model_2} and \eqref{model_3}, our goal is to find a prior $\pi(A)$, such that the following holds for a second-order correct EBCI $I_i$ of $\theta_i$: 
\begin{align*} 
    E[P^{\pi}(\theta_i \in I_{i} |y)] &\stackrel{(a)}= P(\theta_i \in I_{i}|\beta,A) + o(m^{-1}) \\&\stackrel{(b)}= 1 - \alpha + o(m^{-1}), \quad i=1,\cdots, m, \labthis \label{posexp}
\end{align*}
where the expectation on the left hand side (LHS) is with respect to the joint distribution of $(\theta,y)$, i.e., considering Levels 1 and 2 in the model given by \eqref{model_1} and \eqref{model_2}. The LHS in \eqref{posexp} can be interpreted as the expected posterior coverage of $\theta_i$, where $P^{\pi}(\theta_i \in I_{i} |y)$ is defined as the posterior coverage of $I_i$, considering Levels 1, 2 and 3. This criteria is in the same spirit as that of the expected posterior variance of $\theta_i$ considered in \cite{DattaRaoSmith2005}. 
The probability after equation $(a)$ in \eqref{posexp} is with respect to the distribution of $\theta_i$ and $y$, given $\beta$ and $A$ i.e., considering levels 1 and 2 in the model given by \eqref{model_1} and \eqref{model_2} and is termed as the EB coverage. The steps from $(a)$ to $(b)$ in \eqref{posexp} can be explained leveraging \eqref{EB cov}, owing to the second-order correct property of EBCIs such as $I_i^\text{YL}$ and $I_i^\text{N}$. Overall, our matching criteria is given by \eqref{posexp}, that matches expected posterior coverage with EB coverage of an EBCI in a two-level normal model.

We refer to the prior that satisfies \eqref{posexp} as the `matching prior', which is derived in Section \ref{sec:prior} using a higher order asymptotic expansion (upto $O(m^{-1})$) of the left hand side of (\ref{posexp}). We demonstrate the derivation of $\pi(A)$ for the EBCIs $I_{i}^{\text{YL}}$ and $I_{i}^\text{N}$, introduced in Section \ref{sec:int}. Our Bayesian framework is general and can be extended to provide dual justification for any EBCI following this two-step procedure:
\begin{itemize}
    \item[\hypertarget{Step-1}{Step 1}.] Expand the posterior coverage of chosen EBCI and collect the terms of higher order (upto $O_p(m^{-1})$ for large m), including that affecting the choice of prior. 
    \item[\hypertarget{Step-2}{Step 2}.] { Calculate expectation of the posterior coverage and equate the $O(m^{-1})$ terms with zero to reduce the expected posterior coverage error (to $o(m^{-1})$) and obtain a differential equation whose solution gives the expression of the matching prior.} 
\end{itemize}

The rest of this article is organized as follows. In Section \ref{sec:not_reg}, we introduce necessary notations and regularity conditions. In Section \ref{sec:asym}, we provide our result on posterior coverage expansion in Theorem \ref{Th1} and the bias and length comparison result in Proposition \ref{cor5}. In Section \ref{sec:prior}, we derive the matching prior in Theorem \ref{Th2} and provide its propriety condition in Theorem \ref{Th3}. In Section \ref{sec:sp_case} we discuss two important special cases -- two-level model without covariates and the balanced case. In Section \ref{simstudy}, we provide a Monte Carlo simulation study followed by two real data analysis in Section \ref{realdata}. The data and code availability is provided in Section \ref{sec:code}. We end with concluding remarks in Section \ref{sec:dis}, deferring the proofs of main results and related detailed calculations to Appendices \ref{appA} and \ref{appB}, respectively. Some additional simulation results are provided in Appendix \ref{appC}.

\section{Notations and regularity conditions}
\label{sec:not_reg}

We provide below a list of notations that are frequently used in this article:\\
$y = (y_1,\cdots,y_m)'$, an $m\times1$  vector of direct estimates;\\
$X = (x_1, \cdots, x_m)'$, an $m\times p$ known matrix of rank $p$;\\
$\hat{A}$ -- REML estimator of $A$ and
$\Tilde{A}_{i}$ -- AML estimator of $A$ from \citet{Hirose2017} with area-specific adjustment factor $h_{i} \equiv h_{i}(A)$, $l_{i} \equiv l_{i} (A) = \log (h_i)$ and 
$l_{i}^{(1)} \equiv l_{i}^{(1)}(A) = \dfrac{\partial}{\partial A} \log (h_{i}) $;\\ 
$\hat{A}_{i;\text{gls}}$ -- AML estimator of $A$ from \citet{YoshimoriLahiri2014} with area-specific adjustment factor $h_{i;\text{YL}}$, $l_{i;\text{YL}} = \log (h_{i;\text{YL}})$ and 
$l_{i;\text{YL}}^{(1)} = \dfrac{\partial}{\partial A} \log (h_{i;\text{YL}})$; \\
$V \equiv V(A) = \mbox{diag}(A+D_1,\dots,A+D_m)$ and
$\tilde{V} \equiv V(\Tilde{A}_i)$ and $\hat{V} \equiv V(\hat{A})$;\\
$r_i = x_i'(X'V^{-1}X)^{-1}x_i$ and $q_i = x_i'(X'X)^{-1}x_i$, leverage of area $i$;\\
$B_i \equiv B_i(A) = D_i/(A+D_i)$ and $\Tilde{B}_i \equiv B_i(\Tilde{A}_i)$ and
$\hat{B}_i \equiv B_i(\hat{A})$, shrinkage factor for area $i$;\\ 
$g_{1i} \equiv AD_i/(A+D_i), \; g_{2i} \equiv D_i^2 (A+D_i)^{-2} r_i$;\\
$\sigma_i \equiv \sigma_i(A) = \sqrt{g_{1i}} = \sqrt{AD_i/(A+D_i)}$ and $\Tilde{\sigma}_i \equiv \sigma_i(\Tilde{A}_i)$ and
$\hat{\sigma}_i \equiv \sigma_i(\hat{A})$;\\
$\delta_i \equiv \delta_i(A) = \sqrt{g_{1i} + g_{2i}} = \sqrt{AD_i/(A+D_i) + B_i^2\times r_i}$ and $\Tilde{\delta}_i \equiv \delta_i(\Tilde{A}_i)$ and
$\hat{\delta}_i \equiv \delta_i(\hat{A})$;\\
$\Bar{\beta} \equiv \Bar{\beta}(A) = (X'V^{-1}X)^{-1}X'V^{-1}y$, GLS estimator and $\hat\beta_\text{ols} = (X'X)^{-1}X'y$, OLS estimator of $\beta$, and $\Tilde{\beta} \equiv \bar{\beta}(\Tilde{A}_i)$ and $\hat{\beta} \equiv \bar{\beta}(\hat{A})$;\\
$\hat{\theta}^{\text{BLUP}}_i(A) = (1-B_i)y_i + B_i x_i'\Bar{\beta}$ and $\Tilde{\theta}_i \equiv \hat{\theta}^{\text{EB}}_i(\Tilde{A}_i)$ and $\hat{\theta}_i \equiv \hat{\theta}_i^{\text{EB}}(\hat{A}).$

We provide below the regularity conditions required to state our theoretical results.
\begin{enumerate}
    \item[R1:] The logarithm of the adjustment term, i.e., $\log ({h_{i}})$ and $\log ({h_{i;\text{YL}}})$, is free of $y$ and is four times continuously differentiable with respect to $A$. Moreover, $\dfrac{\partial^k}{\partial A^k} \log ( h_{i})$
    and $\dfrac{\partial^k }{\partial A^k} \log (h_{i;\text{YL}})$ 
    are of order $O(1)$, for large $m$, with $k=0,1,2,3$.
    \item[R2:] $|\Tilde{A}_i| < C m^\lambda$ where $C$ is a generic positive constant and $\lambda$ is a small positive constant.
    \item[R3:] Rank($X$) = $p$.
    \item[R4:] The elements of $X$ are uniformly bounded, implying that $\mbox{sup}_{j \ge 1} q_j = O(m^{-1})$.
    \item[R5:] $0 < \mbox{inf}_{j \ge 1} D_j \le \mbox{sup}_{j \ge 1} D_j < \infty, A \in (0,\infty)$.
    \item[R6:] $0 < \hat A < \infty$ and $E[\hat A^{-k}] < \infty$ for $ k = 1, \ldots, 6$.  
    \item[R7:] 
    $\sup_{A \in (0, \infty)} |\rho_k(A)| <\infty$ for $k =1,2$, where $\rho_k(A) \equiv \dfrac{\partial^k}{\partial A^k} \log(\pi(A))$.
\end{enumerate}

\section{Asymptotic expansion of posterior coverage and its expectation}
\label{sec:asym}
This section onward, we refer to the model given by \eqref{model_1} and \eqref{model_2} as M$1$ and the model including \eqref{model_3} in addition to \eqref{model_1} and \eqref{model_2} as M$2$. M$1$ is utilized in the EB approach, whereas M$2$ is leveraged in the HB approach. Throughout this article, we use the symbols $\Phi(\cdot)$ and $\phi(\cdot)$ for the cumulative density function (cdf) and the probability density function (pdf) of $N(0, 1)$, respectively. In Theorem \ref{Th1}, we provide a higher-order asymptotic expansion of the posterior coverage of $I_{i}^\text{N}$. The theorem holds for any area $1 \le i \le m$, for large $m$. 

\begin{theorem}\label{Th1}
Suppose the regularity conditions R$1$-R$5$ hold. Then the posterior coverage of $I_{i}^\text{N}$ under M$2$ is given by,
\begin{align}
    P^{\pi}(\theta_i \in I_{i}^\text{N} |y) = 1 - \alpha + \phi(z)\times \hat c_i + o_p(m^{-1}), \label{th1-2}
\end{align}
where  
$\hat c_i$ is a function of $(y, z,D_i)$ given by
\begin{equation}
\begin{split}
    & \hat c_i = \frac{2 z \hat{B}_i}{\hat{A} \; \mbox{tr}[\hat{V}^{-2}]} \Bigg [ 
        \hat l^{(1)}_{i} - 
        \frac{1}{2}\Bigg \lbrace \Big (\frac{\hat{B}_i}{2} - 2 \Big ) \frac{1}{\hat{A}} + 4 \frac{\mbox{tr}[\hat{V}^{-3}]}{\mbox{tr}[\hat{V}^{-2}]} \Bigg \rbrace 
        - \hat{\rho}_1 - \Bigg ( \frac{y_i - x_i'\hat{\beta}}{\hat{A} + D_i} \Bigg )^2 - \frac{z^2 \hat{B_i}}{4 \hat{A}} \Bigg], \\
    &\hat l^{(1)}_{i} \equiv l^{(1)}_{i} |_{A=\hat{A}} = \frac{2}{\hat{A} + D_i} + \frac{(1+z^2)D_i}{ 4 \hat{A} (\hat{A} + D_i)}\quad \text{and}\quad \hat{\rho}_1 = \frac{\partial}{\partial A} \log \pi(A) \Big |_{A=\hat{A}}.
\end{split}
\label{eq:ci}
\end{equation}
\end{theorem}
The proof of Theorem \ref{Th1} is deferred to the Appendix \ref{sec:app_th}. The posterior coverage of $I_i^{\text{YL}}$ can be obtained similarly as in Theorem \ref{Th1} and we comment on this in the remark below.

\begin{remark} \label{rem:post}
\textbf{[Posterior coverage of $I_i^{\text{YL}}$]}
We note that the posterior coverage expansion of $I_i^{\text{YL}}$ yields a similar expression as in \eqref{th1-2} for which the terms $\hat c_i$ and $\hat{l}_i^{(1)}$ (cf. \eqref{eq:ci}) are different. For the posterior coverage expansion of $I_i^{\text{YL}}$, these are given by 
\begin{equation}
\begin{split}
    & \hat c_{i;\text{YL}} = \hat c_i - \frac{ z \; \hat{B}_i \; \hat r_i }{\hat{A}},\\
    & \hat l^{(1)}_{i;\text{YL}} \equiv l^{(1)}_{i;\text{YL}}|_{A=\hat A} = \hat l^{(1)}_{i} + \mbox{tr}[\hat V^{-2}]\; \frac{\hat r_i}{2},
\end{split}
\label{adj_lYL}
\end{equation}
where $\hat r_i \equiv r_i|_{A=\hat A}$. Interestingly, the additional term involving a leverage type factor ($\hat r_i$) in $\hat c_{i;\text{YL}}$ and $l^{(1)}_{i;\text{YL}}$ from \eqref{adj_lYL} eventually get canceled with each other in the prior derivation. The expressions of the adjustment factor $h_{i;\text{YL}}$ of $I_i^{\text{YL}}$ and $h_i$ of $I_i^\text{N}$ as functions of $A$ are given by 
\begin{equation}
    \begin{split}
    & h_i \equiv h_i(A) = A^{(1+z^2)/4}(A+D_i)^{(7-z^2)/4} \; \text{ [Corollary $1(b)$ of \citet{Hirose2017}]},\\
    & h_{i;\text{YL}} = h_i (A) \times \exp \Big [\int \frac{r_i}{2} \mbox{tr} [V^{-2}] \; dA \Big ] \text{ [Theorem 2 of \citet{YoshimoriLahiri2014}].}
    \end{split}\label{adj_term}
\end{equation} 
In conclusion, $h_{i;\text{YL}}$ includes the leverage term whereas $h_i$ does not, making $h_i$ more straightforward than $h_{i;\text{YL}}$. Moreover, in the balanced case, the adjustment factor $h_i$ and the respective AML estimator $\tilde{A}_i$ reduce to non–area-specific forms, which is not the case for $h_{i;\text{YL}}$ and $\hat{A}_{i;\text{gls}}$, due to the leverage term. We later observe in Proposition \ref{cor5} that the term $l^{(1)}_{i}$ plays a crucial role in the bias and length comparison of EBCIs.
\end{remark}

Having obtained the posterior coverage expansion of $I_i^\text{N}$ in Theorem \ref{Th1}, which is the Step \hyperlink{Step-1}{1} in our two-step procedure, we next proceed to Step \hyperlink{Step-2}{2}, wherein we calculate expected posterior coverage. To this end, we calculate $E[\hat c_i]$ upto the order $O(m^{-1})$ and note that the expectation of the remainder terms of order $o_p(m^{-1})$ in \eqref{th1-2}, can be shown to be of order $o(m^{-1})$. Notice that $\hat c_i$ in \eqref{eq:ci} has a prefactor with $\mbox{tr} [V^{-2}]$ in the denominator which aids in achieving the order $O_p(m^{-1})$ of $\hat c_i$. Considering a term by term expectation approach and using Taylor series expansion, we obtain the expression for $E[\hat c_i]$, which is provided in the following proposition.
\begin{proposition}
\label{propEci}
Under regularity conditions R$1$–R$5$ 
and for $\hat c_i$ as in Theorem \ref{Th1}, we have upto $O(m^{-1})$
\begin{align}
    E[\hat c_i] = \frac{2 z B_i}{{A} \; \mbox{tr}[V^{-2}]} \Bigg [ 
        {l}^{(1)}_{i} - 
        \frac{1}{2}\Bigg \lbrace \Big (\frac{{B}_i}{2} - 2 \Big ) \frac{1}{{A}} + 4 \frac{\mbox{tr}[{V}^{-3}]}{\mbox{tr}[{V}^{-2}]} \Bigg \rbrace 
        - \rho_1 - \frac{1}{A + D_i} - \frac{z^2 {B_i}}{4 {A}} \Bigg], \label{Eci}
\end{align}
where the expectation is with respect to M$1$ and $${l}^{(1)}_{i} = \frac{2}{{A} + D_i} + \frac{(1+z^2)D_i}{ 4 {A} ({A} + D_i)}\quad \text{and}\quad {\rho}_1 = \frac{\partial}{\partial A} \log \pi(A).$$
\end{proposition}
The proof of Proposition \ref{propEci} is deferred to the Appendix \ref{sec:app_cr}. Thus, we have the expression for the LHS of our matching condition given in \eqref{posexp}, i.e., $E[P^{\pi}(\theta_i \in I_{i} |y)]$. Before moving onto the remaining part of Step \hyperlink{Step-2}{2}, i.e., deriving the matching prior, we discuss EB coverage, bias and length of EBCIs in the next section.

\subsection{EB coverage, bias and length of EBCIs}
Notwithstanding the fact that posterior coverage is the main focus of this article, we discuss EB coverage briefly, following new results on bias and length comparison of EBCIs. 
We provide below a remark on the EB coverages of $I_{i}^\text{N}$ and $I_{i}^{\text{YL}}$. The EB coverage expansion is critical for the determination of AML estimators in the interval estimation context, as illustrated in \citet{YoshimoriLahiri2014, Hirose2017}. To facilitate the next discussion on AML estimators, we discuss the EB coverages before, for completeness. 

\begin{remark} \label{rem EB}
\textbf{(i) [EB coverage of $I_i^\text{N}$]} 
For completeness, we explicitly provide the EB coverage expansion of $I_i^\text{N}$ (cf. (A.1) in the proof of Theorem 1 of \citet{Hirose2017} for a general case).
Under regularity conditions stated in \citet{Hirose2017} and model M$1$, we have
\begin{align} \label{th1-1}
P(\theta_i \in I_{i}^\text{N}|\beta,A) = 1 - \alpha + z \; \phi(z) \; \frac{a_i + b_i}{m} + O(m^{-3/2}),
\end{align}
where $a_i$ and $b_i$ are functions of $A, \; D_i$ given by
$$  a_i = -\frac{m}{\mbox{tr}[V^{-2}]}\Big [ \frac{4B_i}{A(A+D_i)} + \frac{(1+z^2)B_i^2}{2A^2}\Big] \quad \text{and} \quad b_i = \frac{2m}{\mbox{tr}[V^{-2}]}\frac{B_i}{A}\times l_i^{(1)}.$$

\textbf{(ii)[EB coverage of $I_i^{\text{YL}}$]} 
Under regularity conditions and model M$1$, Theorem 1 of \citet{YoshimoriLahiri2014} states that
\begin{align} \label{YL exp}
P(\theta_i \in I_{i}^{\text{YL}}|\beta,A) = 1 - \alpha + z \; \phi(z) \; \frac{a_{i;\text{YL}} + b_{i;\text{YL}}}{m} + O(m^{-3/2}),
\end{align}
where $a_{i;\text{YL}}$ and $b_{i;\text{YL}}$ are given by $$a_{i;\text{YL}} = a_i - \frac{m B_i}{A} \; x_i' \;\text{Var}(\Bar{\beta})\; x_i \quad \text{and} \quad  b_{i;\text{YL}} = \frac{2m}{\mbox{tr}[V^{-2}]}\frac{B_i}{A}\times l_{i;\text{YL}}^{(1)}.$$

Comparing \eqref{th1-1} and \eqref{YL exp}, we see that the latter has an additional term in $a_{i;\text{YL}}$ of order $O(m^{-1})$, which, in turn, appears in the adjustment factor discussed in Remark \ref{rem:post} (cf. \eqref{adj_term}). However, both $I_{i}^\text{N}$ and $I_{i}^{\text{YL}}$ choose adjustment terms for their respective AML estimators so that the $O(m^{-1})$ terms vanish in the EB coverage expansions. Hence, both these EBCIs are second-order correct in terms of EB coverage. 
\end{remark}

Following the previous discussion on coverages, we delve into a comparison of the AML estimators $\tilde{A}_i$ and $\hat{A}_{i;\text{gls}}$, used to construct $I_{i}^\text{N}$ and $I_{i}^{\text{YL}}$, respectively. We provide below a new result on the bias comparison of $\tilde{A}_i$ and $\hat{A}_{i;\text{gls}}$. A comparison on the length of the intervals they produce is also provided. 

\begin{proposition}
\label{cor5}
Under regularity conditions R$1$–R$5$ and model M$1$
    \begin{enumerate}
    \item [(i)] $ \text{Bias} [\tilde{A}_i]  < \text{Bias}[\hat{A}_{i;\text{gls}}] $, upto the order of $O(m^{-1})$. 
    \item [(ii)] Length [$I_i^\text{N} (\tilde{A}_i)]$ - Length [$I_i^{\text{YL}} (\hat{A}_{i;\text{gls}})] = O_p(m^{-3/2}).$
    \end{enumerate} 
\end{proposition}
The proof of the above proposition has been deferred to Appendix \ref{sec:app_cr}. However, we comment on its implication in the following remark.

\begin{remark}\textbf{[Bias and length comparison]}
It follows from Proposition \ref{cor5} that $\tilde{A}_i$ exhibits improved performance compared to $\hat{A}_{i;\mathrm{gls}}$ in terms of bias. In contrast, $I_i^\text{N}$ may apparently seem wider than $I_i^\text{YL}$ in terms of length, due to the additional term $g_{2}$ being included in $I_i^\text{N}$. Surprisingly, our calculations yield that the lengths of the two intervals can be regarded as asymptotically equivalent upto the order of $O_p(m^{-1})$. This is due to the fact that, even if
$\delta_i - \sigma_i = O(m^{-1})$ (cf. \eqref{sig_del}), the difference in the lengths of $I_i^\text{N}$ and $I_i^{\text{YL}}$  evaluates $\tilde{\delta}_i - \hat{\sigma}_i \equiv \delta_i (\tilde{A}_i) - \sigma_i (\hat A_{i;\text{gls}})$. In this expression, the terms of order $O_p(m^{-1})$ get canceled (cf. \eqref{length3}) and the length difference arises from the remainder terms in the corresponding higher-order expansions. We also observe that the terms of the order $O_p(m^{-3/2})$ in \eqref{length3} do not have the same sign, i.e., the resulting higher-order length difference expression cannot be shown to be either or negative. In other words, neither of the EBCIs cannot be analytically shown to be better in terms of length upto the order $O_p(m^{-3/2})$. 

Before ending this section, we provide a concluding remark on the condition for positivity of $\tilde{A}_i$ and $\hat{A}_{i;\text{gls}}$, that is when unique and strictly positive estimates exist. 

\begin{remark}\textbf{[Positivity of AML estimators]} \label{rem: pos}
A comparison of the conditions for strict positivity of the AML estimators $\tilde{A}_i$ (Theorem 3 of \citet{Hirose2017}) and $\hat{A}_{i;\text{gls}}$ (Corollary to Theorem 4 of \citet{YoshimoriLahiri2014}), reveals that the condition for $\tilde{A}_i$ is less restrictive, as it does not involve the additional leverage term $q_i \in (0,1)$. Specifically, $\hat{A}_{i;\text{gls}}$ requires $m > (p+4)/(1-q_i)$ for positivity in the balanced case, whereas $\tilde{A}_i$ requires $m > p+4$, which is a milder condition than that imposed for $\hat{A}_{i;\text{gls}}$.
\end{remark}
\end{remark}

\section{Desired prior and its propriety}
\label{sec:prior}
In this section, we elucidate the unique matching prior for which EBCIs $I^\text{N}_i$ and $I^{\text{YL}}_i$ possess a dual justification. Moreover, we show that the posterior obtained from this matching prior is proper, under mild conditions. 

\subsection{Matching prior}

The matching prior of $I_i^\text{N}$ is obtained by setting $E[\hat c_i] = 0$ upto the order $O(m^{-1})$ in Proposition \ref{propEci} (cf. \eqref{Eci}), which is derived by computing the expectation after the posterior expansion of $I_i^\text{N}$ in Theorem \ref{Th1}. Surprisingly, we observe that both $I_i^\text{N}$ and $I_i^\text{YL}$ have the same matching prior, achieving propriety under the same condition. The expression of the unique matching prior for both $I_i^\text{N}$ and $I_i^{\text{YL}}$ is provided next in Theorem \ref{Th2}, which holds for any area $1 \le i \le m$, for large $m$. 

\begin{theorem} \label{Th2}
Suppose that the regularity conditions R$1$-R$5$ hold. Then under the model M$2$, we have the unique matching prior of $I^\text{N}_i$ and $I^{\text{YL}}_i$ as 
    \begin{align} \label{finalprior}
     \pi_i(A) & \propto \mbox{tr} [V^{-2}] (A+D_i) \; A.
\end{align}
\end{theorem}
The proof of Theorem \ref{Th2} is deferred to Appendix \ref{sec:app_th}. The matching prior is improper and is bounded, approaching $0$ as $A$ tends to $0$ and approaching a positive constant $m$ as $A$ tends to infinity. We next examine how this prior compares with those in prior work.

\begin{remark} \textbf{[Comparison with $\pi_{i;\text{DRS}}$]}
Under model M2, we see that the proposed prior (cf. \eqref{finalprior}) is different from the prior $\pi_{i;\mathrm{DRS}}(A)$ in \eqref{drsprior}.
Specifically, we can express our prior as
$$\pi_i(A) = \pi_{i;{\text{DRS}}}(A) \times \frac{A}{A + D_i}.$$
As noted earlier in Section \ref{sec:mat}, $\pi_{i;\mathrm{DRS}}(A)$ depends both collectively and individually on $D_i$'s. The proposed prior in \eqref{finalprior}, although introduces an additional factor $A/(A+D_i)$, has similar area-specific properties as $\pi_{i;\mathrm{DRS}}(A)$. 
In conclusion, the proposed matching prior in interval estimation procedure is area-specific for a general two-level model, following motivation of earlier literature, and differs slightly from the matching prior of \cite{DattaRaoSmith2005, HiroseLahiri2021}.
\end{remark}  

\subsection{Propriety of matching prior}

In this section, we provide our result on the propriety of the matching prior in \eqref{finalprior}. For prior $\pi_i(A),$ we define its propriety as the condition when the posterior is proper, given by
\begin{equation} \label{post}
\begin{split}
    \int_{0}^{\infty}  \pi_i(A|y) \; dA  < \infty,
\end{split}
\end{equation}
where $\pi_i(A|y)$ is the posterior of $A$ defined as 
$$\pi_i(A|y) \propto \pi_i(A) \times L_\text{RE}(A),$$ 
with $L_\text{RE}(A)$ being the REML likelihood of $A$ (cf. \eqref{prop1}). We provide the result on propriety in Theorem \ref{Th3}. To this end, we show that for a specific range of the number of small areas, the matching prior from Theorem \ref{Th2}, yields a proper posterior under mild  conditions.
\begin{theorem} \label{Th3}
    Suppose that that the regularity conditions R$3$-R$5$ hold. Then under the model M$2$, the prior $\pi_i(A)$ in \eqref{finalprior} satisfies the propriety condition in \eqref{post} provided $m > p + 2$.
\end{theorem} 

The proof of Theorem \ref{Th3} is deferred to Appendix \ref{sec:app_th}. Since the matching prior in \eqref{finalprior} is same for both $I_i^\text{N}$ and $I_i^{\text{YL}}$, the propriety condition also remains the same. Moreover, our propriety condition of $m>p+2$ matches exactly with that of \cite{DattaRaoSmith2005}. Interestingly, this propriety condition is different than the positivity condition of $\tilde{A}_i$, discussed in Remark \ref{rem: pos}.

\section{Special cases}
\label{sec:sp_case}
In this section we discuss two important special cases of the model M$1$ -- the two-level model with no covariates given by M$3$, and the balanced case given by M$4$. We provide below in Table \ref{tab:model_summary}, a short summary of the four models M1 to M4, discussed in this article.

\begin{table}[!htb]
\centering
\caption[Model characteristics of M1--M4]{Model characteristics of M1--M4 across Levels 1--3, design features, and existence of EBCIs $I_i^\text{YL}$ and $I_i^\text{N}$.}
\label{tab:model_summary}
\begin{tabular}{|c|c|c|c|c|c|c|c|}
\hline
Model & Level 1 & Level 2 & Level 3 & Balanced & Covariate & $I_i^\text{YL}$ & $I_i^\text{N}$ \\
\hline
M1 & $\checkmark$ & $\checkmark$ &  &  & $\checkmark$ & $\checkmark$ & $\checkmark$ \\
M2 & $\checkmark$ & $\checkmark$ & $\checkmark$ &  & $\checkmark$ & $\checkmark$ & $\checkmark$ \\
M3 & $\checkmark$ & $\checkmark$ & $\checkmark$ &  &  & $\checkmark$ &  \\
M4 & $\checkmark$ & $\checkmark$ & $\checkmark$ & $\checkmark$ & $\checkmark$ & $\checkmark$ & $\checkmark$ \\
\hline
\end{tabular}
\end{table}

\subsection{Two-level model without covariates}
\label{sec:nocov}

We consider the two-level model with $x_i'\beta = 0$, also termed as the common mean model with mean zero. This special case is of particular importance, since \citet{efron1973stein} provided an EB justification of the James–Stein estimator using a balanced case of this model ($D_i = 1,\; i = 1, \cdots, m$). However, we consider the unbalanced case, where the $1^\text{st}$ level of the model M$1$ in \eqref{model_1} remains the same, whereas the $2^\text{nd}$ and $3^\text{rd}$ levels given by \eqref{model_2} and \eqref{model_3}, respectively, simplify to 
\begin{eqnarray}
\mbox{Level 2:}\;&&\theta_i|A \stackrel{\text{iid}}\sim N(0,A); \quad  i=1,\cdots,m, \label{model_22} \\
\mbox{Level 3:}\;&& A \sim \pi(A), \; A >0 \label{model_33}. 
\end{eqnarray}
For brevity, we refer to the model given by \eqref{model_1}, \eqref{model_22} and \eqref{model_33} as M$3$ and the estimators and prior with an added suffix of `sp' to denote the special case. Since the parameter $\beta$ is no longer present, we use the Best Predictor (BP) for $\theta_i$ given by 
$\hat \theta_i^\text{BP} = (1-B_i) y_i,$ considering $A$ to be known. EB estimator based on the respective AML estimator $\tilde{A}_{i;\text{sp}}$ 
simplifies to $\Tilde{\theta}_{i;\text{sp}} = (1 - \Tilde{B}_i)y_i,$ where $\tilde{B}_{i;\text{sp}} \equiv B_i(\tilde{A}_{i;\text{sp}})$. Despite M3 being a special case of M2, $\Tilde{\theta}_{i;\text{sp}}$ is not a special case of $\Tilde{\theta}_{i}$, as we do not get $\Tilde{\theta}_{i;\text{sp}}$ by simply substituting the condition $x_i'\beta = 0$ in $\Tilde{\theta}_{i}$. In this special case, the term $g_{2}$, which is a function of covariates, and consequently the interval $I^\textup{N}$, which is a function of $g_{2}$, are no longer applicable. We only consider one interval in this case, the version of $I_i^\text{YL}$ given by $$I_{i;\text{sp}}^{\text{YL}}: \Tilde{\theta}_{i;\text{sp}} \pm z \Tilde{\sigma}_{i;\text{sp}},$$ where $\Tilde{\sigma}_{i;\text{sp}} \equiv \sigma({\Tilde{A}_{i;\text{sp}}})$. 
Following the proposed two-step framework, similar to the general case, we obtain a higher-order asymptotic expansion of the posterior coverage of $I_{i;\text{sp}}^{\text{YL}}$, which leads to the same matching prior as in \eqref{finalprior}, but with a simplified matching criterion $m>2$, as there are no covariates.



\subsection{Balanced case of the two-level model}
\label{sec:bal}

The balanced case of the two-level model have been extensively discussed in literature -- \citet{YoshimoriLahiri2014, Hirose2017} obtained interesting findings that their AML estimators have unique solutions for $A > 0$ in the balanced case. In this case, in M$1$ we assume that $D_i = D, \; i=1, \cdots, m$, and denote this model as M4. Consequently, GLS and OLS estimators of $\beta$ are identical. In the result given in Corollary \ref{cor1}, we observe that the matching prior given in Theorem \ref{Th2} also reduces to a simpler form, denoted as $\pi_{\text{bal}}(A)$ using the suffix `bal' to signify balanced. 
    \begin{corollary} \label{cor1}
        Under regularity conditions R$1$-R$5$ and the model M$4$, 
        the matching prior for $I_{i}^\text{N}$ and $I_{i}^{\text{YL}}$ is given by
\begin{align} \label{bal}
     \pi_{\text{bal}}(A)
    & \propto \frac{A}{A+D}, 
\end{align}
achieving propriety under the condition $m > p + 2$. 
\end{corollary}
    
The proof of the above corollary is provided in Appendix \ref{sec:app_cr}. Following this, we provide some remarks about the implication of this prior in comparison to existing prior in literature for the balanced case. 

\begin{remark} \textbf{[Comparison with $\pi_{\text{DGM}}$]}
To compare proposed prior with that of $\pi_{\text{DGM}}$, we need to consider the same model as \cite{DGM2000}, given by \eqref{DGM1}-\eqref{DGM3}. This balanced common mean model is equivalent to considering $x_i'\beta = \mu, \; i = 1,\cdots, m$, in M4. In this case, our prior simplifies to \eqref{bal}, which is of the same form as that of $\pi_{\text{DGM}}$, provided in \eqref{dgmprior}. Thus, our prior for two-sided interval estimation of $\theta_i$, matches with the quantile matching prior of \cite{DGM2000}, considering one-sided credible interval in the balanced common mean model. This provides us with the scope to investigate the form of the prior, if a two-sided interval (like the EBCI discussed) is considered instead, in the setup of \cite{DGM2000}. 

Note that the authors derive the quantile matching prior from the coverage expansion of a one-sided credible interval. As per model given by \eqref{DGM1}-\eqref{DGM3}, this is expressed as 
\begin{align}
    P_{A} [ \theta_i > h_i(\pi,\hat A,\alpha;Y)] = \alpha + \frac{z_{\alpha} \phi (z_{\alpha})}{2 m \pi(A)} \frac{\partial}{\partial A} \Bigg [ \frac{\pi(A) B'(A)}{J(A)\{1 - B(A)\}} \Bigg ] + o(m^{-1}), \label{DGM4}
\end{align}
where $P_{A}(.)$ denotes the probability measure based on the joint distribution of $(Y_i, \theta_i),\; i=1,\cdots,m$, as specified by \eqref{DGM1} and \eqref{DGM2}, 
$h_i(\pi,\hat A,\alpha;Y)$ is such that $P^{\pi} \{ \theta_i > h_i(\pi,\hat A,\alpha;Y) | Y \} = \alpha + o_p(m^{-1})$,
$J(A) = (A + D)^{-2}/2, \; B(A) = D (A + D)^{-1}$ and $\hat A$ is the REML estimator of $A$.
Considering a two-sided interval, following \eqref{DGM4} we have
\begin{align*}
    &P_{A} [h_i(\pi,\hat A, 1-\alpha/2 ;Y) < \theta_i < h_i(\pi,\hat A, \alpha/2; Y)]\\ 
    &=  P_{A} [\theta_i > h_i(\pi, \hat A, 1-\alpha/2; Y)] -  P_{A} [\theta_i > h_i(\pi,\hat A, \alpha/2; Y)] 
    \\
    &= 1 - \alpha + \frac{1}{2m \pi(A)} \Big[ z_{1-\alpha/2} \phi (z_{1-\alpha/2}) - z_{\alpha/2} \phi (z_{\alpha/2}) \Big]
    \frac{\partial}{\partial A} \Bigg [ \frac{\pi(A) B'(A)}{J(A)\{1 - B(A)\}} \Bigg ] + o(m^{-1})\\
    &= 1 - \alpha -  \frac{z_{\alpha/2} \phi (z_{\alpha/2})}{m \pi(A)} \frac{\partial}{\partial A} \Bigg [ \frac{\pi(A) B'(A)}{J(A)\{1 - B(A)\}} \Bigg ] + o(m^{-1}), \labthis \label{DGM5}
\end{align*}
using the fact that $\phi (z_{1-\alpha/2}) = \phi (z_{\alpha/2})$ and $z_{1-\alpha/2} = - z_{\alpha/2}$. Hence, RHS of \eqref{DGM5} will be equal to $1-\alpha$ to the order $o(m^{-1})$ if and only if $\pi (A) \propto A (D + A)^{-1},$ which is same as the prior in the one-sided case of \cite{DGM2000}, given in \eqref{dgmprior}. 
\end{remark}

\begin{remark} \textbf{[Comparison with $\pi_{i;\text{DRS}}$]}
The matching prior of \citet{DattaRaoSmith2005, HiroseLahiri2021}, reduces to the Stein's superharmonic prior in the balanced case, i.e., becomes non informative and non area-specific. Whereas, the proposed prior in \eqref{bal} is non area-specific in the balanced case, but is not a flat prior like the case of $\pi_{i;\text{DRS}}(A)$.  
Recall that to achieve the expected posterior coverage to be approximately $1-\alpha$, we set $E[\hat c_i] = 0$ upto $O(m^{-1})$. For the uniform prior, however $E[\hat c_i] \ne 0$ upto $O(m^{-1})$, and we evaluate this in the balanced case. Note that with $\rho_1 = 0$, in the balanced case M4, we have from \eqref{Eci} upto $O(m^{-1})$ 
\begin{align*}
    E[\hat c_i] &= \frac{2z}{A}\frac{D}{(A+D)} \frac{(A+D)^2}{m} \Bigg[ \frac{2}{A+D} + \frac{(1+z^2) D}{4A(A+D)} - \frac{D}{4A(A+D)} + \frac{1}{A} -\frac{2m}{(A+D)^3}\frac{(A+D)^2}{m} \\ &\qquad \qquad \qquad \qquad \qquad \qquad  - \frac{1}{A+D}  -\frac{z^2D}{4A(A+D)}  \Bigg] \\
    &= \frac{2z}{m} \frac{D^2}{A} (A+D) \Bigg [ \frac{2}{(A+D)} - \frac{2}{(A+D)} - \frac{1}{(A+D)} + \frac{1}{A}  \Bigg]\\
    &= \frac{2z}{m} \frac{D^2}{A^2}.\labthis \label{Eciflat}
\end{align*}
Therefore, from \eqref{Eciflat} we see that for flat prior, $E[\hat c_i] > 0$ under M4. For the proposed matching prior, we do not have this $O(m^{-1})$ term. This in turn implies that the flat prior will always give over-coverage, upto $O(m^{-1})$. Moreover, when the higher order terms in the coverage $1 - \alpha + \phi(z) E[\hat c_i] < 1$, using the expression of $E[\hat c_i]$ upto $O(m^{-1})$ from \eqref{Eciflat}, for $\alpha = 0.05$ we arrive at the condition for over-coverage to be 
$\dfrac{D^2}{mA^2} < 0.218$.
\end{remark}  

\section{Simulation Study}
\label{simstudy}
In this section, we present simulation studies reporting expected posterior coverage (EPC) and EB coverage (EBC) of the two intervals of interest -- $I_i^\text{N}$ and $I_i^{\text{YL}}$. Our simulation results are presented to support the theoretical findings that both intervals achieve the nominal EPC and EBC and not for the purpose of a comparative study of existing methods. In this Monte Carlo (MC) simulation study, we provide results from $100$ independent datasets and $M = 200$ MC simulations for each dataset, considering two values of the nominal level: $\alpha = \{0.05, 0.1\}$. We take $m = 15$ and $30$, since in SAE problems the total number of areas is considered to be small and higher values of $m$ generally yield similar conclusions from reasonable methods. 

\subsection{Simulation setting}

We portray mainly two types of simulation settings -- the two-level model, followed by the balanced case discussed in Section \ref{sec:bal}. The details are as follows: 
\paragraph{Unbalanced case - S1 :} In this case, for model M$1$, we consider the number of covariates to be $p = 3$, where we take $x_1$ as a unit vector owing to the intercept and generate  auxiliary variables $x_2$ and $x_3$ independently from $N(0,1)$. The regression coefficient values considered are $\beta = (2,-5,8)'$. To maintain parity with the existing literature \citep{datta2000unified, DattaRaoSmith2005, YoshimoriLahiri2014}, we group the small areas according to specific patterns of $D_i$'s, i.e., $m = 15$ and $30$ small areas are clustered into $5$ groups of $3$ and $5$ areas, respectively, having same $D_i$ values in each group. We consider the following two patterns of $D_i$'s, keeping in mind that the ratio of $D_i$'s and $A$ should be small so that M$1$ and regularity conditions hold. 
    \begin{enumerate}
        \item[\textbf{S11.}] The unique $D_i$'s for groups G1 - G5 are $(0.3, 0.4, 0.5, 0.6, 0.7)$ and $A=1$. 
        \item[\textbf{S12.}] The unique $D_i$'s for groups G1 - G5 are $(0.1, 0.4, 0.5, 0.6, 3)$ and $A = 1$. 
    \end{enumerate}
Note that, setup S11 is a nearly balanced case with very similar $D_i$'s, whereas S12 has more variability with minimum $D_i$ as $0.1$ in G1 and maximum $D_i$ as $3$ in G5.

\paragraph{Balanced case - S2:} For the balanced case, we consider the same values of $p, \; x$ and $\beta$ as in S1, with the following patterns of $D$ and $A$ values:
\begin{enumerate}
        \item[\textbf{S21.}] $A = 1$ and $D = \{2,1.5,1,0.5\}$,
        \item[\textbf{S22.}] $A = 2$ and $D = \{4,3,2,1\}$.
    \end{enumerate} 
    Note that in setup S22, the values of $A$ and $D$ are double as that of S21.

\subsection{Simulation results}

\paragraph{EPC}  -- We first report EPC values of the two intervals $I_i^\text{N}$ and $I_i^{\text{YL}}$. 
We use \texttt{Rstan} to define our matching prior and in each MC iteration, we draw a total of $2,000$ samples from the posterior distribution of $\theta_i$, with $500$ samples in each of four chains, after discarding a burn-in of size $500$. We then check the percentage of such samples falling inside each EBCI, which gives the posterior coverage (PC). To calculate the EPC, keeping $(\beta, A)$ as fixed, we generate $100$ independent $(\theta,y)$ values, where for each dataset we perform the posterior sample draw and PC calculation for $M$ MC iterations. Finally, averaging over $100$ independent datasets and $200$ MC iterations in each dataset, we get the EPC.
    
    
In Table \ref{tab: 11}, we provide EPC for five groups G1 - G5, averaged over the unique $D_i$'s in each group for $m=15$ areas. For the simulation settings S11 and S12, we report EPC of both EBCIs for $\alpha = \{0.05, 0.1\}$.
In Table \ref{tab4}, EPC value for same setups are provided for $m=30$.
We observe that robustness in EPC is  achieved for $m=15$ with improvement for $m=30$, specially noticeable in EPC for $I_i^{YL}$ in G5 for S12 (where $D_i$'s are maximum). 
    For the balanced cases S21 and S22, since there is no group structure in the $D_i$'s, we provide a single value in each setting, averaged over all areas. These are given in Table \ref{tab: 12} for $m=15$ and Table \ref{tab5} for $m=30$. 
    We observe that in case of both S21 and S22 as the $D/A$ ratio increases, for example in the case of $D=2$ and $A=1$ in S21 and $D=4$ and $A=2$ in S22, the EPC values show under-coverage, especially for $I_i^{\text{YL}}$.
    Finally, we conclude that in both unbalanced and balanced cases, reported in Tables \ref{tab: 11} to \ref{tab5}, EPC of both $I_i^\text{N}$ and $I_i^{\text{YL}}$ are close to nominal levels, for appropriate ratio of $A$ and $D$. This corroborates with our theoretical finding that under the chosen matching prior, both EBCIs can be characterized as approximate credible intervals with expected posterior coverage close to the nominal value.

    \begin{table}[!htb]
\caption[Expected posterior coverage, empirical Bayes coverage and length for S11 and S12 and $m=15$]{Expected posterior coverage (EPC) and empirical Bayes coverage (EBC) in \% and length in parenthesis of $I_i^\text{N}$ and $I_i^{\text{YL}}$ for the simulation settings S11 and S12, taking number of areas $m=15$.}
\label{tab: 11}
\footnotesize
\begin{center}
\begin{tabular}{|c|c|c|c|c|c|c|}
\hline
& &  & \multicolumn{2}{c|}{S11} & \multicolumn{2}{c|}{S12} \\
\hline
$\alpha$ & Measure & Group & $I_i^\text{N}$ & $I_i^{\text{YL}}$ & $I_i^\text{N}$ & $I_i^{\text{YL}}$ \\
\hline 
0.05 & EPC & G1 & 94.8 & 94.5 &94.8 & 94.7\\
 &  & G2 & 94.7 & 94.4& 94.8 & 94.4\\
 &  & G3 & 94.8 & 94.4 & 94.9 & 94.4\\
 &  & G4 & 94.9 & 94.4 & 94.9 & 94.3\\
 &  & G5 & 94.8 & 94.1 & 95.2 & 93.4\\
\cline{2-7}
 & EBC & G1 & 95.7 (2.0) & 95.8 (2.0)  & 95.7 (1.2) & 95.7 (1.2)\\
 & (Length) & G2 &  94.2 (2.2) & 94.7 (2.2)  & 94.0 (2.2) & 94.7 (2.3)\\
 &  & G3 &  95.8 (2.4) & 96.2 (2.5)  & 95.8 (2.4) & 95.8 (2.5)\\
 &  & G4 &  94.8 (2.6) & 95.3 (2.6)  & 94.5 (2.6) & 95.8 (2.7)\\
 &  & G5 &  92.5 (2.7) & 93.8 (2.8)  & 91.3 (3.9) & 96.0 (4.4)\\
\hline 
0.1 & EPC & G1 & 89.7 & 89.4 & 89.7 & 89.6\\ 
  && G2 & 89.6 & 89.0 & 89.7 & 89.1\\ 
  && G3 & 89.8 & 89.1 & 89.9 & 89.2\\ 
  && G4 & 90.0 & 89.0 & 90.1 & 89.0\\ 
  && G5 & 89.7 & 88.6 & 90.5 & 87.2\\
\cline{2-7}
 & EBC & G1 &  90.2 (1.6) & 90.8 (1.7)  & 90.3 (1.0) & 90.3 (1.0) \\
 & (Length) & G2 &  88.8 (1.8) & 89.5 (1.9)  & 89.0 (1.9) & 89.5 (1.9)\\
 &  & G3 &  91.0 (2.0) & 92.0 (2.1)  & 90.5 (2.0) & 91.8 (2.1)\\
 &  & G4 &  88.7 (2.1) & 89.5 (2.2)  & 88.5 (2.2) & 89.7 (2.2)\\
 &  & G5 &  86.0 (2.3) & 87.8 (2.3)  & 84.5 (3.2) & 90.0 (3.7)\\
 \hline
\end{tabular}
\end{center}
\end{table}

\begin{table}[!htb]
\caption[Expected posterior coverage, empirical Bayes coverage and length for S11 and S22 and $m=30$]{EPC and EBC in \% and length in parenthesis of $I_i^\text{N}$ and $I_i^{\text{YL}}$ for the simulation settings S11 and S12, taking number of areas $m=30$.}
\label{tab4}
\footnotesize
\begin{center}
\begin{tabular}{|c|c|c|c|c|c|c|}
\hline
& &  & \multicolumn{2}{c|}{S11} & \multicolumn{2}{c|}{S12} \\
\hline
$\alpha$ & Measure & Group & $I_i^\text{N}$ & $I_i^{\text{YL}}$ & $I_i^\text{N}$ & $I_i^{\text{YL}}$ \\
\hline 
0.05 & EPC & G1 & 94.8 & 94.4 & 94.9 & 94.7 \\ 
 &  & G2 & 94.9 & 94.4 & 94.9 & 94.5 \\ 
 &  & G3 & 94.9 & 94.5 & 95.1 & 94.7 \\ 
 &  & G4 & 94.9 & 94.4 & 95.0 & 94.5 \\
 &  & G5 & 95.0 & 94.5  & 95.3 & 94.2 \\
\cline{2-7}
 & EBC & G1 &  94.5 (1.9) & 94.7 (1.9)  & 95.2 (1.2) & 95.4 (1.2)\\
 & (Length) & G2 &  94.1 (2.1) & 94.7 (2.2)  &94.2 (2.1) & 94.8 (2.2)\\
 &  & G3 &  94.4 (2.3) & 94.9 (2.3)  &94.3 (2.3) & 94.9 (2.4)\\
 &  & G4 &  94.5 (2.5) & 95.1 (2.5)  & 94.8 (2.5) & 95.2 (2.5) \\
 &  & G5 &  93.9 (2.6) & 94.4 (2.6)  & 93.2 (3.6) & 94.9 (3.8) \\
\hline 
0.1 & EPC & G1 & 89.7 & 89.1 & 89.8 & 89.6 \\  
  && G2 & 89.8 & 89.1 & 90.0 & 89.2 \\ 
  && G3 & 89.9 & 89.3 & 90.2 & 89.5 \\ 
  && G4 & 90.9 & 89.1 & 90.0 & 89.1 \\
  && G5 & 90.0 & 88.2 & 90.7 & 88.7 \\
\cline{2-7}
 & EBC & G1 &  89.0 (1.6) & 89.9 (1.6)  & 89.4 (1.0) & 89.7 (1.0)\\
 & (Length) & G2 &  89.3 (1.8) & 89.9 (1.8)  & 88.9 (1.8) & 89.7 (1.8)\\
 &  & G3 &  88.2 (1.9) & 88.7 (2.0)  & 88.9 (1.9) & 89.4 (2.0)\\
 &  & G4 &  89.7 (2.1) & 90.5 (2.1)  & 89.5 (2.1) & 90.6 (2.1)\\
 &  & G5 &  88.5 (2.2) & 89.4 (2.2)  & 87.8 (3.0) & 90.1 (3.2) \\
 \hline
\end{tabular}
\end{center}
\end{table}

\begin{table}[!htb]
\caption[Expected posterior coverage, empirical Bayes coverage and length for S21 and S22 and $m=15$]{EPC and EBC in \% and length in parenthesis of $I_i^\text{N}$ and $I_i^{\text{YL}}$ for simulation settings S21 and S22, taking number of areas $m=15$.}
\label{tab: 12}
\footnotesize
\begin{center}
\begin{tabular}{|c|c|c|c|c|c|c|c|c|c|c|}
\hline
 \multirow{8}{*}{S21} & &  & \multicolumn{2}{c|}{$D = 0.5$} & \multicolumn{2}{c|}{$D = 1$} & \multicolumn{2}{c|}{$D = 1.5$}& \multicolumn{2}{c|}{$D = 2$} \\
\cline{2-11}
& $\alpha$ & Measure & $I_i^\text{N}$ & $I_i^{\text{YL}}$ & $I_i^\text{N}$ & $I_i^{\text{YL}}$ & $I_i^\text{N}$ & $I_i^{\text{YL}}$ & $I_i^\text{N}$ & $I_i^{\text{YL}}$ \\
\cline{2-11} \cline{2-11}
& 0.05 & EPC & 94.8 & 94.1 & 94.4 & 93.6 & 94.1 & 93.2 & 93.9 & 93.0\\
& & EBC   & 94.2 & 95.3  &  94.5 & 95.5  &  95.1 & 95.6  & 95.2 & 95.8 \\ 
& & (Length) &  (2.4) & (2.5)  &  (3.1) & (3.3)  &  (3.7) & (3.9)  &  (4.2) & (4.4) \\ 
\cline{2-11}
& 0.1 & EPC & 89.6 & 88.4 & 89.0 & 87.5 & 88.4 & 86.8 & 88.0 & 86.5\\
&& EBC & 88.4 &  90.0 &  88.6 & 90.2  &  89.0 & 90.2  & 89.5 & 90.7 \\ 
 && (Length) & (2.0) & (2.1) &  (2.6) & (2.7)  &  (3.1) & (3.2)  & (3.5) & (3.6) \\ 
 \hline 
 \multirow{8}{*}{S22} & &  & \multicolumn{2}{c|}{$D = 1$} & \multicolumn{2}{c|}{$D = 2$} & \multicolumn{2}{c|}{$D = 3$}& \multicolumn{2}{c|}{$D = 4$} \\
\cline{2-11}
& $\alpha$ & Measure & $I_i^\text{N}$ & $I_i^{\text{YL}}$ & $I_i^\text{N}$ & $I_i^{\text{YL}}$ & $I_i^\text{N}$ & $I_i^{\text{YL}}$ & $I_i^\text{N}$ & $I_i^{\text{YL}}$ \\
\cline{2-11}\cline{2-11}
& 0.05 & EPC & 94.6 & 94.2 & 94.2 & 93.6 & 93.9 & 93.3 & 93.6 & 93.1\\
& & EBC  &  94.3 & 94.9  &  94.7 & 95.2  &  95.1 & 95.4  &  95.3 & 95.6 \\ 
& & (Length) &   (3.4) & (3.4)   &  (4.5) & (4.6)  &  (5.3) & (5.4)  &  (6.0) & (6.1) \\ 
\cline{2-11}
& 0.1 & EPC & 89.4 & 88.7 & 88.6 & 87.8 & 88.1 & 87.2 & 87.6 & 86.8\\
&& EBC &  88.7 & 89.6  &  88.9 & 89.7  &  89.3 & 90.0  &  89.8 & 90.4 \\ 
 && (Length) &  (2.8) & (2.9)  &  (3.7) & (3.8)  &  (4.4) & (4.5)  &  (5.0) & (5.0) \\ 
 \hline
\end{tabular}
\end{center}
\end{table}

\begin{table}[!htb]
\caption[Expected posterior coverage, empirical Bayes coverage and length for S21 and S22 and $m=30$]{EPC and EBC in \% and length in parenthesis of $I_i^\text{N}$ and $I_i^{\text{YL}}$ for simulation settings S21 and S22, taking number of areas $m=30$.}
\label{tab5}
\footnotesize
\begin{center}
\begin{tabular}{|c|c|c|c|c|c|c|c|c|c|c|}
\hline
 \multirow{8}{*}{S21} & &  & \multicolumn{2}{c|}{$D = 0.5$} & \multicolumn{2}{c|}{$D = 1$} & \multicolumn{2}{c|}{$D = 1.5$}& \multicolumn{2}{c|}{$D = 2$} \\
\cline{2-11}
& $\alpha$ & Measure & $I_i^\text{N}$ & $I_i^{\text{YL}}$ & $I_i^\text{N}$ & $I_i^{\text{YL}}$ & $I_i^\text{N}$ & $I_i^{\text{YL}}$ & $I_i^\text{N}$ & $I_i^{\text{YL}}$ \\
\cline{2-11} \cline{2-11}
& 0.05 & EPC & 94.9  & 94.3 & 94.4 & 93.5 & 94.1 & 93.2& 93.7 & 92.8\\
& & EBC &  94.2 & 95.0  &  94.4 & 95.1  &  94.8 & 95.4  &  95.2 & 95.7 \\ 
& & (Length) &  (2.3) & (2.4)  &  (2.9) & (3.0)  &  (3.4) & (3.5)  &  (3.8) & (3.9) \\ 
\cline{2-11}
& 0.1 & EPC & 89.8 & 88.7 & 89.0 & 87.6 & 88.4 & 86.8 & 87.8 & 86.3 \\
&& EBC &  88.5 & 89.5  &  88.5 & 89.6  &  88.7 & 90.0  &  89.3 & 90.4 \\ 
 && (Length) &  (1.9) & (2.0)  &  (2.4) & (2.5)  &  (2.8) & (2.9)  &  (3.1) & (3.2)\\ 
 \hline
 \multirow{8}{*}{S22} & &  & \multicolumn{2}{c|}{$D = 1$} & \multicolumn{2}{c|}{$D = 2$} & \multicolumn{2}{c|}{$D = 3$}& \multicolumn{2}{c|}{$D = 4$} \\
\cline{2-11}
& $\alpha$ & Measure & $I_i^\text{N}$ & $I_i^{\text{YL}}$ & $I_i^\text{N}$ & $I_i^{\text{YL}}$ & $I_i^\text{N}$ & $I_i^{\text{YL}}$ & $I_i^\text{N}$ & $I_i^{\text{YL}}$ \\
\cline{2-11}\cline{2-11}
& 0.05 & EPC & 94.6 & 94.2 & 94.1 & 93.6 & 93.8 & 93.2 & 93.3 & 92.8 \\
& & EBC &  94.2 & 94.7  &  94.4 & 94.8  &  94.9 & 95.2   &  95.2 & 95.5 \\ 
& & (Length) &  (3.3) & (3.3)  &  (4.2) & (4.2)  &  (4.8) & (4.9)  &  (5.3) & (5.4) \\ 
\cline{2-11}
& 0.1 & EPC & 89.4 & 88.8 & 88.5 & 87.7 & 87.9 & 87.1 & 87.2 & 86.4 \\
&& EBC &  88.5 & 89.2  &  88.6 & 89.2  &  89.0 & 89.6  &  89.4 & 90.0 \\ 
 && (Length) &  (2.7) & (2.8)  &  (3.5) & (3.5)  &  (4.0) & (4.0)  &  (4.4) & (4.5) \\ 
 \hline
\end{tabular}
\end{center}
\end{table}

\paragraph{EBC} -- In this section, we provide EBC values and lengths of $I_i^\text{N}$ and $I_i^{\text{YL}}$ in a similar fashion as in EPC. Tables \ref{tab: 11} and \ref{tab4} report values for unbalanced setups S11 and S12 for $m=15$ and $m=30$ respectively. Tables \ref{tab: 12} and \ref{tab5} report values for balanced setups S21 and S22 for $m=15$ and $m=30$ respectively. In this case, keeping $(\beta, A)$ as fixed, we generate the $(\theta,y)$ in 100 independent iterations and report EBC as the percentage of cases when an EBCI contains the true $\theta$. In the unbalanced cases, referring to values in Tables \ref{tab: 11} and \ref{tab4}, we observe that overall both intervals $I_i^\text{N}$ and $I_i^{\text{YL}}$ achieve EBC close to the nominal value, except for under-coverage for $I_i^\text{N}$ in G5 where $D_i$'s are highest compared to the other groups. This phenomenon improves as number of areas increases from $m=15$ to $m=30$.
    In the balanced cases, referring to values in Tables \ref{tab: 12} and \ref{tab5}, we observe that both EBCIs have EBC close to nominal levels. We also calculate the lengths of the two EBCIs in this section, which are provided in parenthesis alongside EBC values in Tables \ref{tab: 11} to \ref{tab5}. We observe from the empirical results that the lengths of $I_i^\text{N}$ and $I_i^{\text{YL}}$ are equivalent, with very little or no difference. This supports our theoretical finding that the EBCIs are similar in terms of length. 
    Hence, overall our simulation results corroborate with our theoretical findings.

\section{Data Analysis}
\label{realdata}
In this section, we analyze two real datasets, which we discuss in the following subsections.

\subsection{Baseball data analysis}

The baseball data on performance of baseball players was originally studied by \citet{efron1975data} in the context of Stein's estimator and later by \citet{MorrisJASA} in an EBCI context. This data, collected from the New York Times, contains batting averages of $18$ major league players through their first $45$ official at-bats of the $1970$ season. The highest batting average in this dataset is of an unusually good hitter named Clemente, whose batting average is $0.4$, and the least batting average is that of Alvis, whose value is $0.2$. We are interested in interval estimation of each player's batting average ($\theta_i$'s) for the remainder of the season. The true season batting average of a player for the 1970 season is also available in this dataset. In order to consider the batting averages as direct estimates, we use the arc-sin transformation as in \citet{efron1975data}, and take $y_i$ as the transformed batting average of Player $i; \; i = 1,\cdots, 18$. We demonstrate the balanced common mean model in this data analysis given by \eqref{DGM1} - \eqref{DGM3} in Section \ref{sec:mat}, and compute $I_i^\text{YL}$ only, as $I_i^\text{N}$ does not exist in this case. For brevity we denote this EBCI as $I_i^\text{YL}$ (M3). Following Section 2 of \citet{efron1975data}, we take $D=1$ in Level 1 and estimate the unknown common mean in Level 2 by the sample average value: $\hat \mu = -3.275$.  

\begin{figure}[!htb]
  \centering
  \begin{minipage}{0.45\textwidth}
    \caption{Plot of true batting averages ($\theta_i$) and EBCI $I_i^{\text{YL}}$ (M3) for $m=18$ players.}
\label{baseball_pic}
\centering
\includegraphics[width=\linewidth]{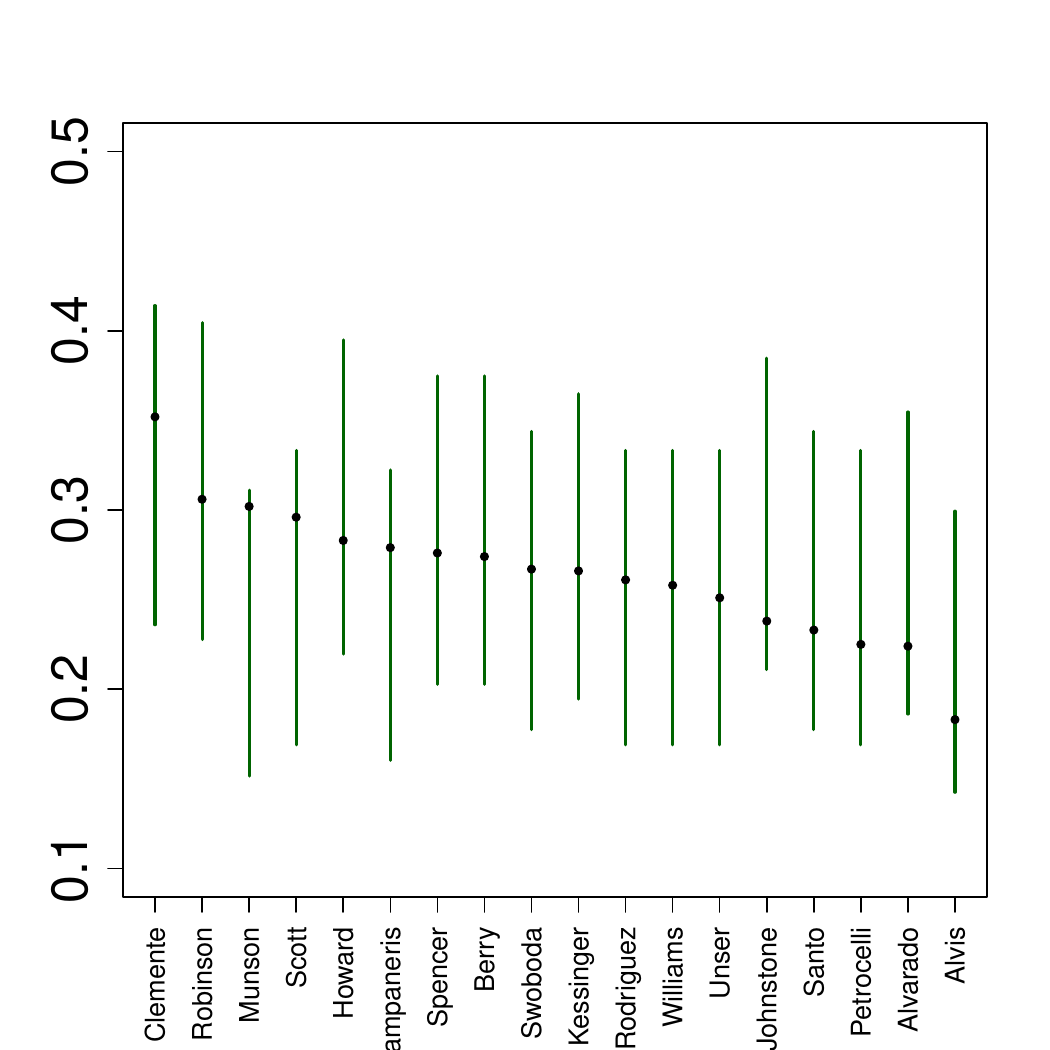}
  \end{minipage}
  \hspace{0.05\textwidth}
  \begin{minipage}{0.45\textwidth}
    \caption{Plot of PC and Monte Carlo error for $I_i^{\text{YL}}$ (M3) for $m=18$ players.}
\label{baseball_pic_2}
\centering
\includegraphics[width=\linewidth]{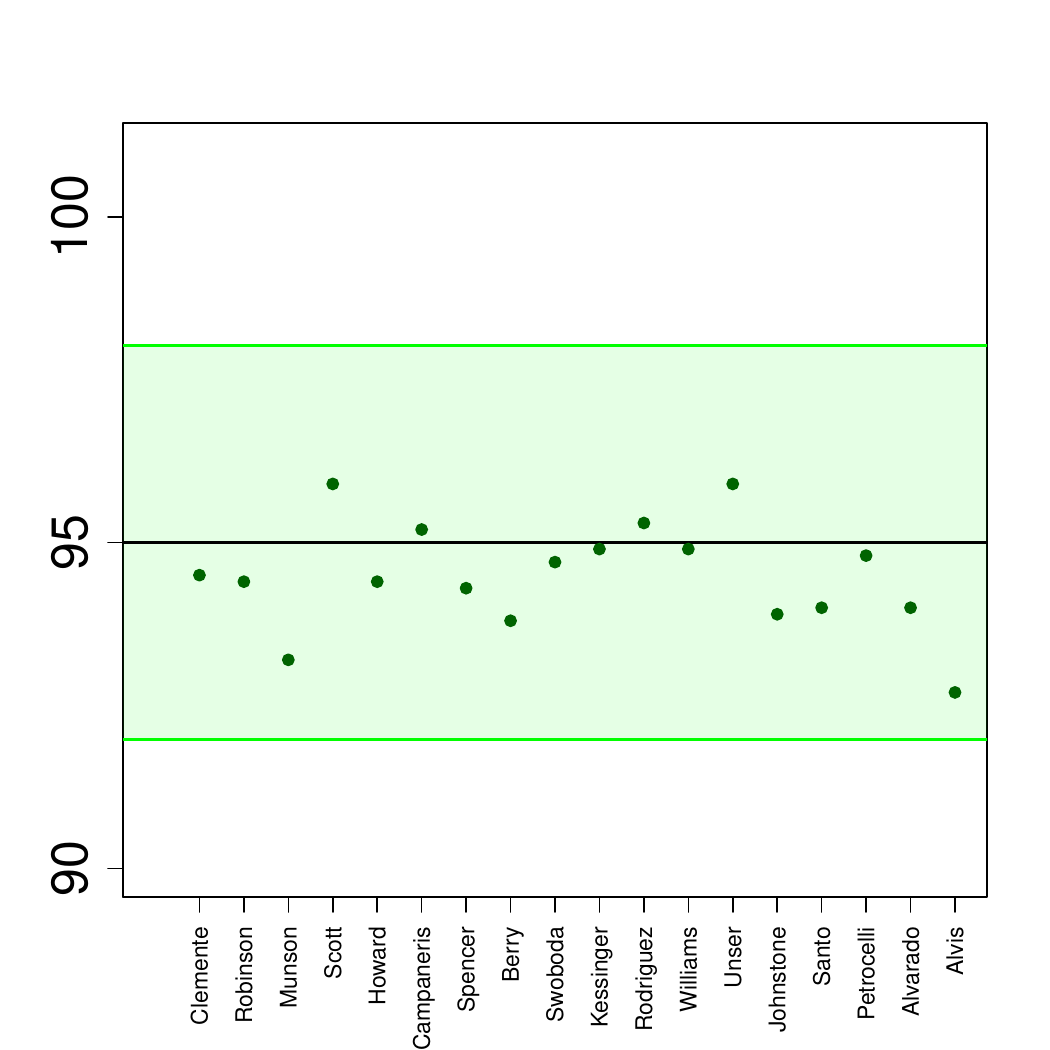}
  \end{minipage}
\end{figure}

In Figure \ref{baseball_pic}, we plot  $I_i^{\text{YL}}$ (M3) in color dark green with the true batting average values in black dots and observe that the intervals contain true values for all players. In Figure \ref{baseball_pic_2}, we provide PC of $I_i^{\text{YL}}$ (M3) in dark green dots along with the MC errors in light green band. PC is calculated similarly as done in the simulation section, averaging from $M =200$ MC iterations taking $2,000$ posterior samples in each iteration. Note that, since this is a real data where $y$ is fixed, we cannot calculate EPC. We observe in Figure \ref{baseball_pic_2} that PCs of $I_i^\text{YL}$ (M3) are all close to $95\%$ and within the Monte Carlo Error (MCE) bands. In conclusion, the baseball data analysis demonstrates that the proposed method is practically implementable on real data, and that the second-order correct EBCI achieves posterior coverage close to the nominal level.

\subsection{BRFSS data analysis}

In this section, we analyze survey data from the Behavioral Risk Factor Surveillance System (BRFSS), a nationwide health surveillance program administered by the Centers for Disease Control and Prevention (CDC) in partnership with state health agencies. Established in 1984, BRFSS collects information from adult U.S. residents through monthly telephone interviews conducted via both landline and cellular telephones. The survey gathers information on health-related behaviors (such as cigarette smoking, alcohol use), chronic health conditions (such as diabetes, heart diseases, asthma, arthritis), and preventive health practices (such as physical activity) that are important for monitoring population health and informing public health policy. Responses collected by individual states are compiled and processed by the CDC, which provides standardized survey methodology and technical support. With more than 400,000 adult interviews conducted annually from all 50 United States (U.S.) states, the District of Columbia (DC), and three U.S. territories, BRFSS is recognized as the world’s largest continuously operating health survey system. 
For a brief overview of BRFSS, refer to the documentation \url{https://www.cdc.gov/brfss/data_documentation/pdf/UserguideJune2013.pdf}.

For our analysis, we consider BRFSS data from 50 U.S. states and DC on five health and preventive outcomes for the year 2024, obtained from the \cite{BRFSS2024}. These outcomes include percentages (\%) of people who ever had diabetes, \% of people who had flu shot within a year, \% of people who visited a doctor for a routine checkup within a year, \% of people who have ever experienced a depressive disorder and \% of people who have smoked at least 100 cigarettes in their entire lives. The state Tennessee (TN) is excluded from the analysis due to unavailable or incomplete data for the selected outcomes and auxiliary variables. Following Model M2, for area $i$, $y_i$ denotes the direct survey-weighted estimate of the corresponding health outcome prevalence, obtained using \texttt{survey} R package. The sampling variance $D_i$ is obtained using a generalized variance function (GVF) smoothing approach applied from \cite{otto1995sampling} to the state-level direct variance estimates, where $i=1,\ldots,m$ and $m=50$.

To improve estimation efficiency, we incorporate auxiliary information from a larger 
survey, the American Community Survey (ACS) conducted by the U.S. Census Bureau, sampling approximately 3.5 million addresses annually to produce detailed insights on various characteristics of the U.S. population. (\cite{ACS2024}). 
Auxiliary variables from the ACS 1-year 2024 Data Profile tables were obtained using the \texttt{tidycensus} R package (\cite{tidycensus}) and include demographic, socioeconomic, and healthcare-access related variables. Different sets of auxiliary variables are considered for different outcomes based on their socio-economic and healthcare relevance. The corresponding covariates include age composition, sex composition, race composition, poverty, Supplemental Nutrition Assistance Program (SNAP) participation, insurance coverage, educational attainment, disability prevalence, unemployment, and housing-related variables. A complete list of outcomes and auxiliary variables used in each model M2 is provided in Table~\ref{tab:brfss_acs_models_2024}.

\begin{table}[ht]
\centering
\caption{BRFSS outcomes and ACS covariates used in M2 for 2024.}
\label{tab:brfss_acs_models_2024}
\begin{tabular}{|l|l|}
\hline
Outcome & Covariates (in percentages) \\
\hline
Diabetes 
& Age 65+, Age 25--34, Female, White, Insured \\

Smoker 
& Age 65+, Age 25--34, Female, White, Poverty, SNAP \\

Flu Shot 
& Age 65+, Age 25--34, Female, White, Low education, Insured \\

Checkup 
& Age 65+, Age 25--34, Female, White, Low education, Insured \\

Depression 
& \begin{tabular}{l}
     Age 65+, Age 25--34, Female, White, Disability, Unemployed, Poverty\\
     Percentage of gross rent of the monthly income
\end{tabular} \\
\hline
\end{tabular}
\end{table}

\begin{figure}
    \centering
    \includegraphics[width=0.9\linewidth, page = 1]{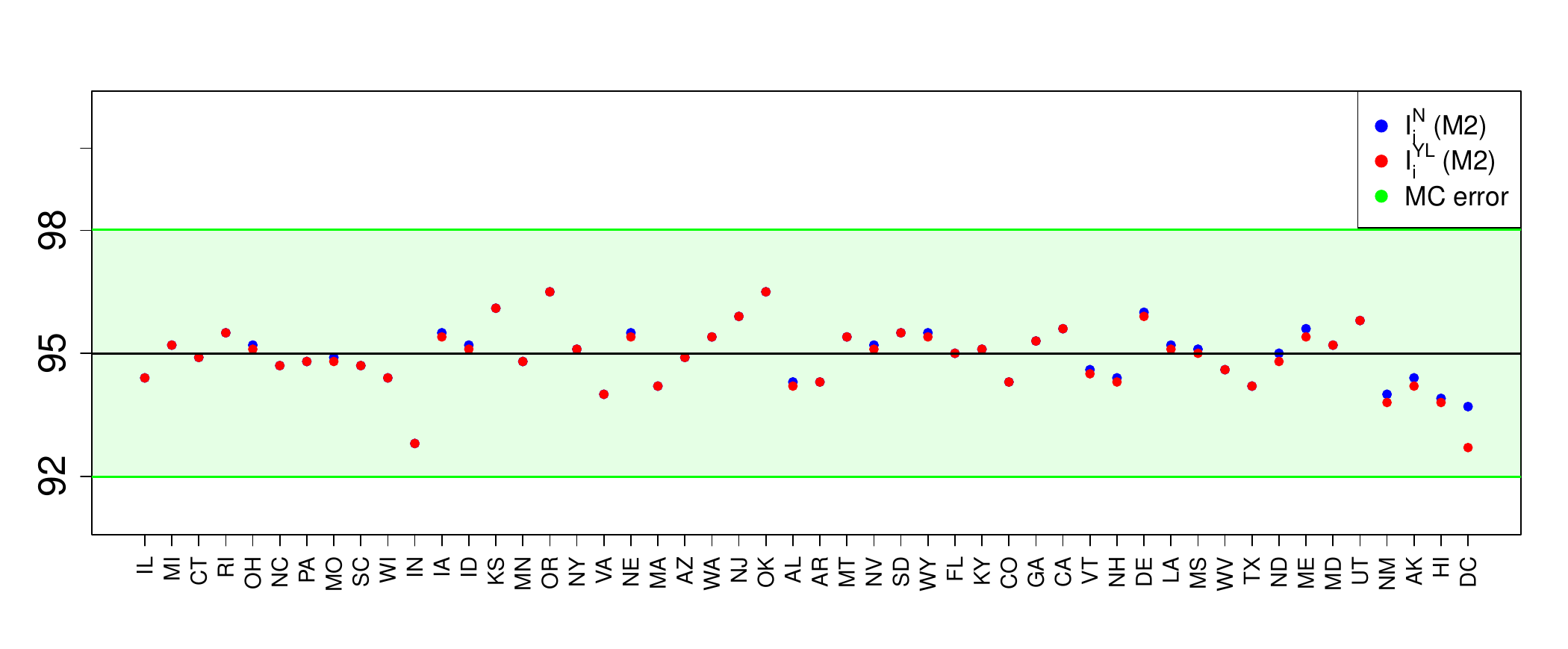}
    \caption{Plot of PC and Monte Carlo error for two EBCIs:  $I_i^\text{N}$ (M2) and $I_i^{\text{YL}}$ (M2) of percentage of diabetes patients for 49 states and DC from the BRFSS data year 2024.}
    \label{brfss_pc_dia}
\end{figure}

\begin{figure}
    \centering
    \includegraphics[width=0.9\linewidth, page = 2]{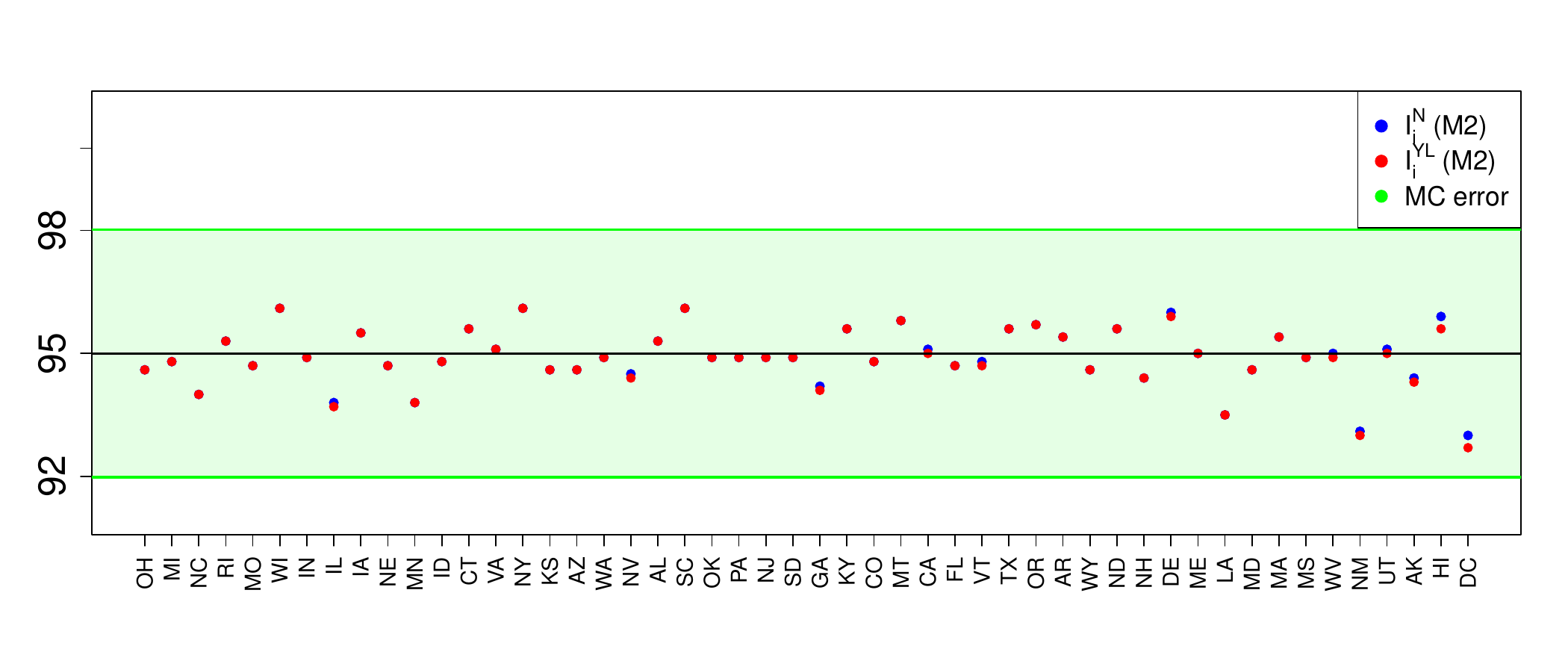}
    \caption{Plot of PC and Monte Carlo error for two EBCIs:  $I_i^\text{N}$ (M2) and $I_i^{\text{YL}}$ (M2) of percentage of smokers for 49 states and DC from the BRFSS data year 2024.}
    \label{brfss_pc_smk}
\end{figure}

We construct two 95\% EBCIs for the true small area parameter $\theta_i$, representing the BRFSS outcome in area $i$ and plot the PC of $I_i^{\text{N}}$ (M2) in blue colored dots and that of $I_i^{\text{YL}}$ (M2) in red dots. 
We present the PC for percentage of diabetes patients and smokers in Figures \ref{brfss_pc_dia} and \ref{brfss_pc_smk}, respectively, with the 49 states (excluding TN) and DC arranged in an ascending order of the leverage values of the model M2. 
The results for other outcomes are deferred to Appendix \ref{appD} for interest of space. The empirical results indicate that PC of both EBCIs remain close to the nominal level for the majority of the states, with all coverage values lying within the MCE bounds. Overall, the two EBCIs demonstrate comparable performance and provide reliable uncertainty quantification for state-level prevalence estimation. This study highlights the usefulness of small area estimation techniques for analyzing large-scale public health survey data. By combining direct BRFSS survey estimates with auxiliary behavioral and demographic information from ACS, the proposed framework yields more stable and accurate state-level estimates while maintaining approximately nominal posterior coverage. The analysis also demonstrates the practical value of empirical Bayesian methods in supporting public health monitoring and policy-making applications.


\section{Data availability}
\label{sec:code}
The codes for simulation study (Section \ref{simstudy}), real data analysis and dataset (Section \ref{realdata}) are available at 
\url{https://github.com/asen123/matching_prior}
in GitHub. 

\section{Discussion}
\label{sec:dis}

In this article, we propose a framework that facilitates a hierarchical Bayesian justification for second-order efficient EBCIs for random effects under a two-level normal model. Within this framework, an EBCI is shown to achieve expected posterior coverage that closely aligns with the nominal level for a specific prior on the hyperparameter. Consequently, a chosen EBCI can also be interpreted as an approximate credible interval. The matching prior is area-specific, following earlier literature 
on matching prior, which depend on individual sampling variance for a mean $\theta_i$. The proposed prior differs slightly from those in earlier literature -- derived by matching the posterior variance of the small area mean with the estimated mean squared error or with the integrated Bayes risk of the empirical Bayes estimator. Moreover, through real data analysis of two datasets, including an important socio-economic one, we illustrate that our method is  practically implementable and computationally efficient, requiring minimal processing time. Through the  real data investigations, we observe how the presence of large leverage values affect expected posterior coverage performances. Essentially, this study sheds light on the importance of regularity conditions in achieving good posterior coverage. 

Overall, the proposed approach addresses an important gap in the matching prior literature related to interval estimation of random effects in two-level normal models. In comparison to previous matching prior literature in two-level normal models, this work provides a novel perspective of providing Bayesian justification to EB confidence intervals. The implications and advantages of this analysis are as follows. This investigation provides a useful tool to Bayesians, as instead of calculating credible intervals through complicated procedures, EB confidence intervals can be directly used for interval estimation. In order to find EB coverage of credible intervals, one has to perform integration multiple times. An EB confidence interval uses maximization instead of integration, and are easier to calculate and time-saving as compared to the fully Bayesian route. Moreover, Markov Chain Monte Carlo procedures may give different answers while computing credible intervals multiple times. We would like to acknowledge that these arguments in favor of our contributions, are similar to the appeal of ADM of \cite{morristang} in Bayesian scenarios, when a procedure is to be used repeatedly for model selection and model checking.

An interesting avenue for future investigation would be to establish dual justifications for interval estimation across broader modeling frameworks beyond the two-level normal model discussed in this article, such as multiple variance components in linear mixed models -- specifically, the nested-error regression model 
or non-normal and other hierarchical models.

\section{Acknowledgments}
The second author’s research was partially supported by Kyushu University Diversity and Super Global Training Program for Female and Young Faculty (SENTAN Q) and JSPS KAKENHI grant number 22K01426. The research of the third author was supported, in part, by a sabbatical leave granted by the University of Maryland College Park, USA.

\bibliographystyle{chicago}
\bibliography{sample.bib}

\newpage
\appendix
\renewcommand{\thefigure}{\thesection.\arabic{figure}}
\setcounter{figure}{0}
\renewcommand{\thetable}{\thesection.\arabic{table}}
\setcounter{table}{0}

\section{Appendix}
\label{appA}
In this section, we provide the proofs of theorems in Appendix \ref{sec:app_th} and those of propositions and corollaries in Appendix \ref{sec:cor_prop}. 

\subsection{Proof of theorems}
\label{sec:app_th}

Before proving Theorem \ref{Th1}, we state the following lemmas, which provide detailed expressions of the posterior coverage of $I_i^\text{N}$. Proofs of these lemmas are deferred to Appendix \ref{applemmapf}. 

\begin{lemma} \textbf{[First level expansion of $I_{i}^\text{N}$]} \label{lem-1} Under model M$2$, we have the posterior coverage probability of the upper and lower limits of $I_{i}^\text{N}$ as
\begin{align}  \label{prob}
    P^{\pi}(E_i^\text{U}|y) = E[\Phi\{Q_i(z) + z\}|y] \mbox{ and }
    P^{\pi}(E_i^\text{L}|y) = E[\Phi\{Q_i(-z) - z\}|y],
\end{align}
where the expectations in 
\eqref{prob} are with respect to the posterior distribution $\pi_i(A|y)$ and the sets $E_i^\text{U},\;E_i^\text{L}$ and functions $Q_i(z)$ and $Q_i(-z)$ are given by
\begin{align*}
    &E_i^\text{U} \equiv \lbrace \theta_i < \Tilde{\theta_i} + z \; \Tilde{\delta}_i \rbrace, \\
    &E_i^\text{L} \equiv \lbrace \theta_i < \Tilde{\theta_i} - z \; \Tilde{\delta}_i \rbrace, \\
    & Q_i(z) \equiv Q_i(z,y,A,D_i)
    = \frac{k_i}{\delta_i} + z \Big (\frac{\Tilde{\delta}_i}{\delta_i} - 1 \Big ),\\ & Q_i(-z) \equiv Q_i(-z,y,A,D_i)
    = \frac{k_i}{\delta_i} - z \Big (\frac{\Tilde{\delta}_i}{\delta_i} - 1 \Big ),\\
    &k_i = (B_i - \Tilde{B}_i)y_i + \Tilde{B}_i x_i'\Tilde{\beta} - B_i x_i'\Bar{\beta}.
\end{align*}
\end{lemma}

\begin{lemma} \textbf{[Second level expansion of $I_{i}^\text{N}$]} \label{lem-2}
Under model M$2$, the posterior coverage probability of $I^\text{N}_i$ yields the expansion
\begin{align} \label{lem-2-2}
    P^{\pi}(E_i|y) 
    &= 1 - \alpha + \phi(z) \Big [ \Big (T_1 - T_4 \Big ) - \frac{z}{2} \Big (T_2 + T_5 \Big ) \Big ] + (T_3 - T_6),
\end{align}
where $E_i \equiv \lbrace \theta_i \in I_i^\text{N} \rbrace$ and $T_1, T_2, T_3, T_4, T_5, T_6$ are all functions of $(y,z)$, defined as
\begin{align*}
   &  \; T_1 =
   E(Q_i(z)|y), \quad 
   T_2 = 
   E(Q_i^{2}(z)|y), \\ 
   & T_4 
   = E(Q_i(-z)|y), \quad T_5 
   = E(Q_i^{2}(-z)|y),\\
    & T_3 
    = \frac{1}{2} E \Big [ \Big \lbrace \int_{z}^{z + |Q_i(z)|} (z + |Q_i(z)| - x)^2 (x^2 - 1) \phi(x) dx \Big \rbrace \Big |y \Big ], 
    \allowdisplaybreaks\\
    & T_6 
    = \frac{1}{2} E \Big [ \Big \lbrace \int_{-z}^{- z + |Q_i(-z)|} ( - z + |Q_i(-z)| - x)^2 (x^2 - 1) \phi(x) dx \Big \rbrace \Big |y \Big ].
\end{align*}
\end{lemma}

In the next three lemmas, we provide simplified expressions of $(T_1 - T_4)$ and $(T_2 + T_5)$, and the asymptotic order of the remainder term $(T_3 - T_6)$. The proof of these lemmas are deferred to Appendix \ref{applemmapf}, where we use some results from \citet{DattaRaoSmith2005}, stated in Lemma \ref{aux-lem-1} and Lemma \ref{aux-lem-2} in Appendix \ref{auxlem}.

\begin{lemma} \textbf{[Expansion of $(T_1 - T_4)$]} \label{lem-3}
The expression $(T_1 - T_4)$ in \eqref{lem-2-2} yields upto $O_p(m^{-1})$
\begin{align} 
    T_1 - T_4 
    &= \frac{2z \; \hat{B}_i}{\hat{A}\;\mbox{tr}[\hat{V}^{-2}]} \Bigg [ \hat{l}^{(1)}_{i} - 
    \frac{1}{2}\Bigg \lbrace \Big (\frac{\hat{B}_i}{2} - 2 \Big ) \frac{1}{\hat{A}} + 4 \frac{\mbox{tr}[\hat{V}^{-3}]}{\mbox{tr}[\hat{V}^{-2}]} \Bigg \rbrace - \hat{\rho}_1 \Bigg ].
\end{align}
\end{lemma}

\begin{lemma} \textbf{[Expansion of $(T_2 + T_5)$]} \label{lem-4}
The expression $(T_2 + T_5)$ in \eqref{lem-2-2} yields upto $O_p(m^{-1})$
\begin{align}
    T_2+T_5
    &= \frac{4 \hat{B}_i}{\hat{A}\;\mbox{tr}[\hat{V}^{-2}]} \Bigg[ \Bigg ( \frac{y_i - x_i'\hat{\beta}}{\hat{A} + D_i} \Bigg )^2 + \frac{z^2 \hat{B_i}}{4 \hat{A}} \Bigg ].  
\end{align}
\end{lemma}

\begin{lemma} \textbf{[Remainder terms]} \label{lem-5}
The order of the expression $T_3 - T_6$ in \eqref{lem-2-2} is $o_p(m^{-1})$.
\end{lemma}

We now discuss the proof of the Theorem \ref{Th1} as follows. 

\begin{proof}[\textbf{Proof of Theorem \ref{Th1}}] 
\textbf{[Posterior expansion of $I_i^\text{N}$]}

The posterior coverage probability of $I_{i}^\text{N}$ under model M$2$ can be expressed as
\begin{align} \label{cov}
    &P^{\pi}(\theta_i \in I_{i}^\text{N} |y) \equiv P^{\pi}(E_i |y) = P^{\pi}(E_i^\text{U}|y) - P^{\pi}(E_i^\text{L}|y),
\end{align}
whereby using Lemma \ref{lem-1} and Lemma \ref{lem-2} subsequently, we further simplify \eqref{cov} in terms of conditional expectations of $y$, i.e., $T_1, T_2, T_3, T_4, T_5$ and $T_6$.

For the posterior coverage in \eqref{lem-2-2} to be approximately $1-\alpha$, under regularity conditions R$1$-R$5$
in Section \ref{sec:not_reg} and model M$2$, we need to find a prior $\pi(A)$ such that the following holds
\begin{align} \label{stoch} 
\Big (T_1 - T_4 \Big ) - \frac{z}{2} \Big (T_2 + T_5 \Big )  = o_p(m^{-1}).
\end{align}
Combining the expressions of $(T_1 - T_4)$ and $(T_2+T_5)$ upto $O_p(m^{-1})$ from Lemma \ref{lem-3} and Lemma \ref{lem-4} respectively, the left hand side of \eqref{stoch} yields the expression of $\hat c_i$ as follows
\begin{align*}
        \hat c_i &=
        \frac{2 z \hat{B}_i}{\hat{A} \; \mbox{tr}[\hat{V}^{-2}] } \Bigg [ 
        \hat{l}^{(1)}_{i} - 
        \frac{1}{2}\Bigg \lbrace \Big (\frac{\hat{B}_i}{2} - 2 \Big ) \frac{1}{\hat{A}} + 4 \frac{\mbox{tr}[\hat{V}^{-3}]}{\mbox{tr}[\hat{V}^{-2}]} \Bigg \rbrace - \hat{\rho}_1
        - \Bigg ( \frac{y_i - x_i'\hat{\beta}}{\hat{A} + D_i} \Bigg )^2 - \frac{z^2 \hat{B_i}}{4 \hat{A}} \Bigg]. 
    \end{align*}
This concludes the proof.
\end{proof}

\begin{proof}[\textbf{Proof of Theorem \ref{Th2}}]

\textbf{[Matching prior of $I_i^\text{N}$ and $I_i^{\text{YL}}$]}

The expression of the matching prior $\pi(A)$ for $I_i^\text{N}$ is obtained by solving the differential equation that arises from { $E[\hat c_i] = 0$ upto $O(m^{-1})$}. Towards that, using the expression of $E[\hat c_i]$ from \eqref{Eci} we obtain the expression of $\rho_1$ as
\begin{align*}
    \rho_1 
    &= l^{(1)}_{i} - 
        \frac{1}{2}\Bigg \lbrace \Big (\frac{B_i}{2} - 2 \Big ) \frac{1}{A} + 4 \frac{\mbox{tr}[V^{-3}]}{\mbox{tr}[V^{-2}]}\Bigg \rbrace 
        {- \frac{1}{A+D_i}}
        - \frac{z^2 B_i}{4 A} \\
        &= \frac{\partial}{\partial A} \Big[ \log( h_{i}) + \log (\mbox{tr} [V^{-2}]) + \Big (\frac{z^2+1}{4} \Big )\log \Big (\frac{A+D_i}{A} \Big) + \log A 
        { - \log (A+D_i)} \Big ]. \labthis \label{eq:rho}
\end{align*}
Moreover, since $\rho_1 = \dfrac{\partial}{\partial A} \lbrace \log \pi(A) \rbrace$, we get the following differential equation
\begin{align}\label{eq:diff}
    \frac{\partial}{\partial A} \lbrace \log \pi(A) \rbrace
        &= \frac{\partial}{\partial A} \log \Bigg[  h_{i} \; \mbox{tr} [V^{-2}] \Big (\frac{A+D_i}{A} \Big)^{{(z^2+1)}/{4}}  \frac{A}{A+D_i} 
        \Bigg ].
\end{align}
Finally, using the expression of $h_{i}$ given in \eqref{adj_term} in Remark \ref{rem:post}, we obtain the solution as 
\begin{align*} 
    \pi_i(A) 
    & \propto A^{{(1+z^2)}/{4}}(A+D_i)^{{(7-z^2)}/{4}} 
    \Big (\frac{A+D_i}{A} \Big)^{{(z^2+1)}/{4}} \; \mbox{tr}[V^{-2}] \; \frac{A}{A+D_i}  
    \notag \\
    & = \mbox{tr} [V^{-2}] (A+D_i) A. 
    \labthis \label{prN}
\end{align*}

For obtaining the prior for $I_i^{\text{YL}}$, we solve the differential equation arising from { $E[\hat c_{i;\text{YL}}] = 0$ upto $O(m^{-1})$}, where the expression of $\hat c_{i;\text{YL}}$ is noted in \eqref{adj_lYL} in Remark \ref{rem:post}. In this case, we similarly get
\begin{equation*}
\begin{split}
    \rho_{1} 
    &= l^{(1)}_{i;\text{YL}} - \frac{\mbox{tr}[V^{-2}]}{2} r_i -
        \frac{1}{2}\Bigg \lbrace \Big (\frac{B_i}{2} - 2 \Big ) \frac{1}{A} + 4 \frac{\mbox{tr}[V^{-3}]}{\mbox{tr}[V^{-2}]}\Bigg \rbrace 
        - { \frac{1}{A+D_i}}
        - \frac{z^2 B_i}{4 A} \\
        &= \frac{\partial}{\partial A} \Bigg[ \log( h_{i;\text{YL}}) + \log (\mbox{tr} [V^{-2}]) + \Big (\frac{z^2+1}{4} \Big )\log \Big (\frac{A+D_i}{A} \Big) + \log A - \int \frac{\mbox{tr}[V^{-2}]}{2} r_i dA \\ 
        &\qquad \qquad { - \log (A+D_i)} 
        \Bigg ]. 
\end{split}
\end{equation*}
Finally, using the value of $h_{i;\text{YL}}$ from \eqref{adj_term}, we obtain the solution to the above differential equation as 
\begin{align*} 
    \pi_i(A) 
    & \propto A^{{(1+z^2)}/{4}}(A+D_i)^{{(7-z^2)}/{4}} \exp \Bigg \lbrace \int \frac{\mbox{tr}[V^{-2}]}{2} \; r_i \; dA \Bigg \rbrace \; \mbox{tr} [V^{-2}]
    \Big (\frac{A+D_i}{A} \Big)^{{(z^2+1)}/{4}} \frac{A}{A+D_i}  
    \\ & \qquad \times \exp \Bigg \lbrace - \int \frac{\mbox{tr}[V^{-2}]}{2} \; r_i \; dA  
    \Bigg \rbrace  \notag \\
    & = \mbox{tr} [V^{-2}] (A+D_i) A. 
    \labthis \label{prYL}
\end{align*}
Thus, the prior for $I_i^\text{YL}$ derived in \eqref{prYL} is same as that of $I_i^\text{N}$ derived in \eqref{prN}. 
\end{proof}

\begin{proof}[\textbf{Proof of Theorem \ref{Th3}}]\label{app:prop}
\textbf{[Propriety of matching prior]}

To this end, we use the definition of propriety provided in \eqref{post}, where we plug in the expressions of $h_i$ from \eqref{adj_term} and $L_\text{RE}$ as follows
\begin{align*} 
    \int_{0}^{\infty} \pi_i(A|y) dA 
    &\propto \int_{0}^{\infty} \Bigg [
    \mbox{tr} [V^{-2}] (A+D_i) A \; 
     |V|^{-1/2}|X'V^{-1}X|^{-1/2} \exp \left \lbrace - \frac{1}{2} (y'Py) \right \rbrace \Bigg ] dA, \labthis \label{prop1}
\end{align*}
where $P \equiv P(A) = V^{-1}-V^{-1}X ({X'V^{-1}X})^{-1}X'V^{-1}$. 
From \eqref{prop1}, we observe that the exponential term satisfies the following condition 
\begin{align}
  \exp \Big \lbrace - \frac{1}{2} (y'Py) \Big \rbrace \le 1.  \label{expineq}
\end{align}
Now let, $ c_1 = \mbox{inf}_{j \ge 1}  D_j$ and $c_2 = \mbox{sup}_{j \ge 1}  D_j$.
Then, the following relations hold
\begin{align}
    & |V|^{-1/2} = \prod_{j=1}^{m} (A+D_j)^{-1/2} \le (A + c_1)^{-m/2}, \label{IC7}\\ 
    & |X'V^{-1}X|^{-1/2} \le C (A + c_2)^{p/2} \label{IC8}\\
    &\mbox{tr} [V^{-2}] = \sum_{j=1}^{m} (A+D_j)^{-2} \le m (A + c_1)^{-2}, \label{IC9} \\
    &(A + D_i) \le (A + c_2), \label{IC10}
\end{align}
where $C$ is a generic positive constant used in equalities, which does not depend on $A$. Thus, using \eqref{expineq} - \eqref{IC10}, the expression \eqref{prop1} reduces to
\begin{align}  \label{post1}
   0 < \int_0^\infty \pi_i(A|y) \; dA 
    & \le C \int_{0}^{\infty} (A + c_1)^{-m/2 - 2} (A+c_2)^{p/2 + 1}\; A
    \; dA \notag \\
    & \stackrel{(a)}\le C \int_{0}^{\infty} \frac{A}{(A+c_2)^{m/2 - p/2 + 1}} \; dA \notag \\ &\stackrel{(b)}\le C \int_{0}^{\infty} \frac{1}{(A+c_2)^{m/2 - p/2 }} \; dA,
\end{align}
where the inequality $(a)$ holds since $(A+c_2)/(A+c_1)$ is a decreasing function of $A$ and therefore $(A+c_2)/(A+c_1) < c_2/c_1.$ The inequality $(b)$ holds since $A/(A+c_2) < 1,\; \forall A \in (0,\infty)$. 
Finally, the conclusion holds by noting the fact that the last expression in \eqref{post1} is integrable when { $m > p + 2$.}
\end{proof}

\subsection{Proofs of propositions and corollaries}
\label{sec:app_cr}
\begin{proof}
[\textbf{Proof of Proposition \ref{propEci}}]
We recall from \eqref{eq:ci} that 
\begin{align*}
        \hat c_i &=
        \frac{2 z \hat{B}_i}{\hat{A} \; \mbox{tr}[\hat{V}^{-2}] } \Bigg [ \hat{l}^{(1)}_{i}
        - 
        \frac{1}{2}\Bigg \lbrace \Big (\frac{\hat{B}_i}{2} - 2 \Big ) \frac{1}{\hat{A}} + 4 \frac{\mbox{tr}[\hat{V}^{-3}]}{\mbox{tr}[\hat{V}^{-2}]} \Bigg \rbrace - \hat{\rho}_1
        - \Bigg ( \frac{y_i - x_i'\hat{\beta}}{\hat{A} + D_i} \Bigg )^2 - \frac{z^2 \hat{B_i}}{4 \hat{A}} \Bigg] \\
        &= 2z \times [\hat c_{i1} - \hat c_{i2} + \hat c_{i3} + \hat c_{i4} - \hat c_{i5} - \hat c_{i6} - \hat c_{i7}]
    \end{align*}
where 
\begin{align*}
    &\hat c_{i1} = \frac{\hat{B}_i \hat{l}^{(1)}_{i}}{\hat{A} \; \mbox{tr}[\hat{V}^{-2}] }, \; 
    \hat c_{i2} = \frac{\hat{B}_i^2}{4\hat{A}^2 \; \mbox{tr}[\hat{V}^{-2}] }, \; \hat c_{i3} = \frac{\hat{B}_i}{\hat{A}^2 \; \mbox{tr}[\hat{V}^{-2}] }, \; 
    \hat c_{i4} =  \frac{\hat{B}_i \mbox{tr}[\hat{V}^{-3}]}{\hat{A} \; \mbox{tr}[(\hat{V}^{-2}])^2}, \\ 
    &\hat c_{i5} = \frac{\hat{B}_i \hat{\rho_1}}{\hat{A} \; \mbox{tr}[\hat{V}^{-2}]}, 
    \; \hat c_{i6} = \frac{\hat{B}_i}{\hat{A} \; \mbox{tr}[\hat{V}^{-2}]} \Bigg ( \frac{y_i - x_i'\hat{\beta}}{\hat{A} + D_i} \Bigg )^2 , 
    \; \hat c_{i7} = \frac{z^2\hat{B}_i^2}{4\hat{A}^2 \; \mbox{tr}[\hat{V}^{-2}]}. 
\end{align*}
Note that from the terms $\hat c_{ij}, \; j=1,\cdots,7$, defined above, only $\hat c_{i6}$ involves $y_i$ explicitly. Hence, we compute $E[\hat c_{i1}]$ and comment that the other terms, except $\hat c_{i5}$ and $\hat c_{i6}$, can be handled similarly. Due to the involvement of $\rho$ in $\hat c_{i5}$ and $y_i$ in $\hat c_{i6}$, we compute $E[\hat c_{i5}]$ and $E[\hat c_{i6}]$ separately.

Letting $\hat c_{1i} = f_m(\hat A)$ and using Taylor expansion on $f_m(\hat A)$ around $A$, we get
\begin{align}
    f_m(\hat A) = f_m(A) + f_m'(A^*)(\hat A - A), \label{c1i}
\end{align}
where $ f_m(A) = \dfrac{{B}_i {l}^{(1)}_{i}}{{A} \; \mbox{tr}[{V}^{-2}]}$ and $A^* = \gamma_1 \hat A + (1 - \gamma_1) A $, for some $\gamma_1 \in [0,1]$. Taking expectation of \eqref{c1i} with respect to M1, we get 
\begin{align*}
    E[f_m(\hat A)] &= f_m(A) + E[f_m'(A^*)(\hat A - A)] = O(m^{-1}) + E[r_1], \labthis \label{fm1} 
\end{align*}
where we denote the remainder term as $r_1 = f_m'(A^*)(\hat A - A)$ and use the fact that $f_m(A)$ is of the order $O(m^{-1})$. We next find the order of $f_m'(A)$ in terms of $A$. Let $ B_i l_i^{(1)}/A = \zeta_i(A) \equiv \zeta_i$ and $\mbox{tr}[{V}^{-2}] = s_m(A) \equiv s_m$. Then, using chain rule of derivatives, we have
\begin{align*}
    f_m'(A) &= \frac{\partial}{\partial A} \Big [ \zeta_i \{s_m ^{-1} \} \Big ] = \{s_m ^{-1} \} \frac{\partial}{\partial A} \zeta_i + \zeta_i \frac{\partial}{\partial A} \Big[ \{s_m ^{-1} \} \Big]. \labthis \label{er_prop3}
\end{align*}
Note that
\begin{align*}
&\frac{\partial}{\partial A}\zeta_i\\
&= -\frac{1}{A^2}B_il_i^{(1)}
-\frac{l_i^{(1)}}{A}\frac{D_i}{(A+D_i)^2}
+\frac{B_i}{A} \frac{\partial}{\partial A}
\left[ \frac{2}{A+D_i} + \frac{(1+z^2)D_i}{4A(A+D_i)} \right] \\
&= -\frac{B_il_i^{(1)}}{A^2}-\frac{l_i^{(1)}D_i}{A(A+D_i)^2} - \frac{B_i}{A}\left[\frac{2}{(A+D_i)^2} + \frac{(1+z^2)}{4} \left\{ \frac{D_i}{A(A+D_i)^2} + \frac{D_i}{A^2(A+D_i)}\right\} \right].
\end{align*}
We take $ c_1 = \mbox{inf}_{j \ge 1}  D_j$ and $c_2 = \mbox{sup}_{j \ge 1}  D_j$, as in the proof of Theorem \ref{Th3}.
Now using 
the fact that $s_m^{-1}
= 1/\sum_{j=1}^m (A+D_j)^{-2} \le (A+c_2)^2/m$, $1/(A+D_i) \le 1/(A+c_1) \le 1/c_1$, $(A+c_2)/(A+c_1) < c_2/c_1$ and $B_i <1$, we find the orders of all terms involved in the first term of \eqref{er_prop3} as follows
\begin{align*}
& \frac{B_il_i^{(1)}}{s_m} \le \frac{D_i}{A+D_i}
\left[ \frac{2}{A+D_i} + \frac{1+z^2}{4A(A+D_i)}
\right] \frac{(A+c_2)^2}{m} \le \frac{c_2}{m} \left( 2+\frac{1+z^2}{4A}
\right) \left( \frac{A+c_2}{A+c_1} \right)^2 \lesssim \frac{1}{Am},
\\
& \frac{B_il_i^{(1)}}{A^2}s_m^{-1}
\lesssim \frac{1}{A^3m},
\\
& \frac{l_i^{(1)}D_i}{A(A+D_i)^2}s_m^{-1}
= \frac{B_il_i^{(1)}}{A(A+D_i)}s_m^{-1} \lesssim \frac{1}{A^2m},
\\
& \frac{B_i}{A}\frac{2}{(A+D_i)^2}s_m^{-1}
\le \frac{B_i}{A} \frac{2}{(A+c_1)^2}
\frac{(A+c_2)^2}{m} \lesssim \frac{1}{Am},
\\
&\frac{B_i}{A}\frac{(1+z^2)}{4} 
\frac{D_i}{A(A+D_i)^2}
s_m^{-1}
\le \frac{B_i}{A} \frac{(1+z^2)}{4} \frac{c_2}{A(A+c_1)^2} \frac{(A+c_2)^2}{m} \lesssim \frac{1}{A^2m},\\
& \frac{B_i}{A} \frac{(1+z^2)}{4}
\frac{D_i}{A^2(A+D_i)}s_m^{-1}
\lesssim \frac{1}{A^3}\frac{(1+z^2)}{4}
\frac{c_2^2}{(A+c_1)^2}
\frac{(A+c_2)^2}{m} \lesssim \frac{1}{A^3m}.
\end{align*}
Therefore, the first term in \eqref{er_prop3} reduces to
\begin{align*}
\left|s_m^{-1} \frac{\partial}{\partial A}\zeta_i \right|
\lesssim \frac{1}{m} \left( \frac{1}{A^3} + \frac{1}{A^2} + \frac{1}{A} \right).\labthis \label{er_prop4}
\end{align*}
For the second term in \eqref{er_prop3} we have
\begin{align*}
\left|\zeta_i \frac{\partial}{\partial A} \{s_m^{-1} \} \right|
& \lesssim \left|
\frac{D_i}{A+D_i} \cdot \frac{1}{A(A+D_i)}
\left[ 2+\frac{(1+z^2)D_i}{4A}
\right] \frac{\sum_{j=1}^m (A+D_j)^{-3}}
{\left\{\sum_{j=1}^m (A+D_j)^{-2}\right\}^2} \right| \\
& \le \frac{D_i}{A+D_i} \cdot \frac{1}{A(A+D_i)}
\left[ 2+\frac{(1+z^2)D_i}{4A} \right]
\frac{1}{m} \frac{(A+c_2)^4}{(A+c_1)^3} \\
&\le \frac{1}{mA}\left[ 2+\frac{(1+z^2)D_i}{4A} \right] \left( \frac{A+c_2}{A+c_1} \right)^4 \\
&\le \frac{1}{mA} 
\left[ 2+\frac{(1+z^2)c_2}{4A} \right] \left( \frac{c_2}{c_1} \right)^4 \\
&\lesssim \frac{1}{m} \left( \frac{1}{A^2} + \frac{1}{A}\right). \labthis \label{er_prop5}
\end{align*}
Next, we show the remainder term in \eqref{fm1} to be ignorable after expectation, i.e., the order of $E[r_1]$ is $o(m^{-1})$.
Using \eqref{er_prop4} and \eqref{er_prop5}, we compute as follows
\begin{align*}
    &E\Big[|f'_m(A^*)(\hat A - A)|\Big] \\
    & \stackrel{(i)}{\le} \sqrt{E[f'_m(A^*)]^2} \sqrt{E(\hat A - A)^2}\\
    & \lesssim \frac{1}{m} \sqrt{E\left( \frac{1}{{A^*}} + \frac{1}{{A^*}^2} + \frac{1}{{A^*}^3} \right)^2} \sqrt{E(\hat A - A)^2}\\
    & \stackrel{(ii)}{\lesssim} \frac{1}{m} \sqrt{E\left( \frac{1}{{A^*}^2} + \frac{1}{{A^*}^4} + \frac{1}{{A^*}^6} \right)} \sqrt{E(\hat A - A)^2}\\
    & \stackrel{(iii)}{\lesssim} \frac{1}{m} \min\left\{ \sqrt{\frac{1}{A^2}+\frac{1}{A^4}+ \frac{1}{A^6}}, \sqrt{E\left( \frac{1}{{\hat A}^2} + \frac{1}{{\hat A}^4} + \frac{1}{{\hat A}^6} \right)} \right\} \sqrt{E(\hat A - A)^2}\\
    & \stackrel{(iv)}{\lesssim} \frac{1}{m} \min\left\{ \frac{1}{A}+ \frac{1}{A^2} + \frac{1}{A^3}, \sqrt{E\left(\frac{1}{\hat A^2}\right)} + \sqrt{E\left(\frac{1}{\hat A^4}\right)} + \sqrt{E\left(\frac{1}{\hat A^6}\right)} \right\} O(m^{-{1/2}})\\
    & \stackrel{(v)}{=} O(m^{-3/2}),
\end{align*}
where $(i)$ uses Cauchy Schwartz (CS) inequality, $(ii)$ follows form $C_r$-inequality, $(iii)$ uses the fact that $A^* \ge \min\{A,\hat A\}>0$, $(iv)$ uses the facts that for any positive real numbers $x_1,x_2,x_3 >0$, $\sqrt{x_1+x_2+x_3} \le \sqrt{x_1}+ \sqrt{x_2} + \sqrt{x_3}$ and $E(\hat A - A)^2 = O(m^{-{1}})$ and lastly, $(v)$ uses the regularity condition R6.
Finally, we have 
\begin{align}
    E[\hat c_{i1}] = \frac{{B}_i {l}^{(1)}_{i}}{{A} \; \mbox{tr}[{V}^{-2}]} + o(m^{-1}). \label{ci1}
\end{align}

We next perform the calculation for $\hat c_{i5} = \dfrac{\hat{B}_i \hat{\rho_1}}{\hat{A} \; \mbox{tr}[\hat{V}^{-2}]} \equiv f_m(\hat A)$. Now, the order of $f'_m(A)$ in terms of $A$ is given by
\begin{align*}
    \left|\frac{\partial}{\partial A} \left( \frac{B_i \rho_1}{A \;\text{tr}[V^{-2}]} \right) \right| 
    &= \left| - \frac{B_i \rho_1}{A(A+D_i) \;\text{tr} [V^{-2}]} + \frac{B_i \rho_2}{A \;\text{tr}[V^{-2}]} - \frac{B_i \rho_1}{A^2 \;\text{tr}[V^{-2}]} + \frac{B_i \rho_1\; \text{tr}[V^{-3}]}{A \;(\text{tr}[V^{-2}])^2} \right|\\
    & \le \frac{c_2 (A+c_2)^2\rho_1}{mA(A+c_1)^2} + \frac{c_2 (A+c_2)^2\rho_2}{mA(A+c_1)} + \frac{c_2 (A+c_2)^2\rho_1}{mA^2(A+c_1)} + \frac{2c_2 (A+c_2)^4\rho_1}{mA(A+c_1)^4}\\
    & \lesssim \frac{\rho_1}{m} \left[ \frac{1}{A} + \frac{1}{A^2} \right] + \frac{\rho_2}{mA}, \labthis \label{er_prop6}
\end{align*}
where the last inequality follows from the facts that $(A+c_2)/(A+c_1) < c_2/c_1$ and $1/(A+c_1)\le 1/c_1$.
Therefore, using \eqref{er_prop6},
we calculate the order of the remainder term as follows
\begin{align*}
    & E \Big[\left|f_m'(A^*) (\hat A - A)\right| \Big] \\&\lesssim E \left| \left[  \frac{\rho_1^*}{m} \left( \frac{1}{A^*} + \frac{1}{{A^*}^2} \right) + \frac{\rho_2^*}{m A^*} \right] (\hat A - A)\right|\\
    &\le \frac{1}{m} E \left| \left[  \left\{\sup_{A\in(0,\infty)} \rho_1(A)\right\} \left( \frac{1}{A^*} + \frac{1}{{A^*}^2} \right) + \left\{\sup_{A\in(0,\infty)} \rho_2(A)\right\} \frac{1}{A^*} \right] (\hat A - A)\right|\\
    & \stackrel{(i)}{\lesssim} \frac{1}{m} E \left| \left[  \left( \frac{1}{A^*} + \frac{1}{{A^*}^2}  \right)\right] (\hat A - A)\right| \\
    & \stackrel{(ii)}{\lesssim} \frac{1}{m} \sqrt{ E \left[ \frac{1}{{A^*}^2} + \frac{1}{{A^*}^4}\right] } \sqrt{ E \left(\hat A - A\right)^2 } \\
    & \stackrel{(iii)}{\lesssim} \frac{1}{m} \min \left\{ \frac{1}{A} + \frac{1}{A^2}, \sqrt{ E \left[ \frac{1}{{\hat A}^2} + \frac{1}{{\hat A}^4}\right] } \right\}  O(m^{-1/2})\\
    & \stackrel{(iv)}{=} O(m^{-3/2}),
\end{align*}
where $(i)$ uses regularity condition R7, i.e., $\rho_k^* = \rho_k(A^*) \le \sup_{A \in (0, \infty)} \rho_k(A) < \infty$ for $k=1,2$, $(ii)$ uses CS inequality and $C_r$-inequality, $(iii)$ uses the facts that $A^* \ge \min\{\hat A, A\}>0$, for any positive real numbers $x_1,x_2>0$, $\sqrt{x_1+x_2} \le \sqrt{x_1}+ \sqrt{x_2}$ and $E(\hat A - A)^2 = O(m^{-{1}})$, lastly, $(iv)$ uses the regularity condition R6.
Thus, noting the fact that $\dfrac{B_i\rho_1}{A \; \text{tr}[V^{-2}]} = O(m^{-1})$, we get
\begin{align}
    E[\hat c_{i5}] = \dfrac{B_i\rho_1}{A \; \text{tr}[V^{-2}]} + o(m^{-1}).\label{ci5}
\end{align}
We next work out $E[\hat c_{i6}]$ by considering $f_m(\hat A) = \hat c_{i6}$.
Using Taylor expansion and taking expectation as before, we get 
\begin{align*}
    E[f_m(\hat A)] &= E[f_m(A)] + E[f_m'(A^*)(\hat A - A)],
\end{align*}
where 
\begin{align*}
    &f_m(A) = \frac{{B}_i}{{A} \; \mbox{tr}[{V}^{-2}]} \Big ( \frac{y_i - x_i'{\beta}}{{A} + D_i} \Big )^2,\\
    &f_m'(A^*) = (y_i - x_i'{\beta})^2 \frac{\partial}{\partial A} \Bigg \lbrace \frac{{B}_i}{{A} \; \mbox{tr}[{V}^{-2}] (A + D_i)^2} \Bigg \rbrace \Bigg |_{A = A^*}.
\end{align*}
Note that marginally, $(y_i - x_i'{\beta}) \sim N(0, A+D_i)$, which implies $E[(y_i - x_i'{\beta})^2] = (A+D_i)$. Thus, we have $E[f_m(A)] = \dfrac{{B}_i}{{A} \; \mbox{tr}[{V}^{-2}]} \dfrac{1}{({A} + D_i)}  = O(m^{-1}) $. 
Next we show that $E[f_m'(A^*)(\hat A - A)] = o(m^{-1})$. Note that 
\begin{align*}
&E\left[\left|f_m'(A^*)(\hat A-A)\right|\right]\\
&= E\left[\left|(y_i-x_i'\beta)^2 \times \left.
\frac{\partial}{\partial A}
\left\{
\frac{B_i} {A\,\text{tr}[V^{-2}](A+D_i)^2}
\right\} \right|_{A=A^*}
(\hat A-A) \right|
\right] \\
&\le \sqrt{ E\!\left[
(y_i-x_i'\beta)^4 \left( \left.
\frac{\partial}{\partial A}
\left\{\frac{B_i} {A\,\text{tr}[V^{-2}](A+D_i)^2}
\right\} \right|_{A=A^*}
\right)^2 \right]
} \times \sqrt{E(\hat A-A)^2} \,
, \labthis \label{eq_prop2}
\end{align*}
from CS inequality. For the second term in \eqref{eq_prop2} we know that $E(\hat A-A)^2 = O(m^{-1})$. For the first term in \eqref{eq_prop2}, using H\"older's inequality (with components $p = 3$ and $q = 3/2$) we get that
\begin{align*}
&E\!\left[ (y_i-x_i'\beta)^4
\left( \left.
\frac{\partial}{\partial A}
\left\{ \frac{B_i}  {A\,\text{tr}[V^{-2}](A+D_i)^2}
\right\}
\right|_{A = A^*} \right)^2
\right]
\\
&\le \Bigl(
E\bigl[|y_i-x_i'\beta|^{4 \times 3}\bigr]\Bigr)^{1/3}
\Biggl(E\!\left[ \Bigg| \left.
\frac{\partial}{\partial A}
\left\{
\frac{B_i}
     {A\,\text{tr}[V^{-2}](A+D_i)^2}
\right\}
\right|_{A = A^*}
\Bigg|^{2\times \frac{3}{2}}
\right]
\Biggr)^{2/3}
\\
&=
O(1)\,
\Biggl( E\!\left[ \Bigg| \left.
\frac{\partial}{\partial A}
\left\{
\frac{B_i}
     {A\,\text{tr}[V^{-2}](A+D_i)^2}
\right\} \right|_{A = A^*} \Bigg|^3
\right] \Biggr)^{2/3},  \labthis \label{eq_prop1} 
\end{align*}
where the last inequality in \eqref{eq_prop1} follows due to the fact that $E[|y_i - x_i'{\beta}|^{12}]  \propto (A+D_i)^{6}$ which in turn implies $E[|y_i - x_i'{\beta}|^{12}]= O(1)$. Now, we calculate the derivative as follows
\begin{align*}
& \left|\frac{\partial}{\partial A}
\left\{
\frac{B_i}
     {A\,\text{tr}[V^{-2}](A+D_i)^2}
\right\} \right|\\
&=
\left|\frac{\partial}{\partial A}
\left\{
\frac{D_i}
     {A\,\text{tr}[V^{-2}](A+D_i)^3}
\right\} \right|\\
& = \left|
-\frac{D_i}
       {A^2\text{tr}[V^{-2}](A+D_i)^3}
-\frac{3D_i}
       {A\text{tr}[V^{-2}](A+D_i)^4}
+ \frac{2D_i}
       {A(A+D_i)^3}
 \frac{\text{tr}[V^{-3}]}
      {\{\text{tr}[V^{-2})\}^2} \right|
\\ 
&  \lesssim
\frac{c_2(A+c_2)^2}
     {A^2\,m\,(A+c_1)^3} +
\frac{3c_2(A+c_2)^2}
     {A\,m\,(A+c_1)^4}+
\frac{2c_2}{A(A+c_1)^3}
\frac{(A+c_2)^4}
     {m(A+c_1)^3}\\
& \lesssim \frac{1}{m} \left( \frac{1}{A^2} + \frac{1}{A} \right), \labthis \label{eq_prop3}
\end{align*}
where the last inequality follows from the facts that $(A+c_2)/(A+c_1) < c_2/c_1$ and $1/(A+c_1)\le 1/c_1$.
Using \eqref{eq_prop3}, we have
\begin{align*}
E\!\left[
\Bigg|
\left.
\frac{\partial}{\partial A}
\left\{
\frac{B_i}
     {A\,\text{tr}[V^{-2}](A+D_i)^2}
\right\}
\Bigg|_{A = A^*}
\right|^3
\right]
&\stackrel{(i)}{\lesssim} \frac{1}{m^3 } E\left[ \frac{1}{{A^*}^6} + \frac{1}{{A^*}^3} \right]\\
& \stackrel{(ii)}{\lesssim} \frac{1}{m^3}\min\left\{ \frac{1}{A^6} + \frac{1}{A^3}, E\left(\frac{1}{\hat A^6} + \frac{1}{\hat A^3}\right)\right\}\\
& \stackrel{(iii)}{=} O(m^{-3}), \labthis \label{eq_prop3.1}
\end{align*}
where $(i)$ follows from $C_r$-inequality, $(ii)$ follows from the fact that $ A^* \ge \min\{A, \hat A\} >0$ and $(iii)$ uses the regularity condition R6. Therefore, combining \eqref{eq_prop1} and \eqref{eq_prop3.1}, we obtain the order of \eqref{eq_prop1} to be $O(m^{-2})$.
Hence, from \eqref{eq_prop2} we get
\begin{align*}
    E[f_m'(A^*)(\hat A - A)] = O(m^{-1/2}) O(m^{-1}) = O(m^{-3/2})
\end{align*}
and
\begin{align}
    E[\hat c_{i6}] = \frac{{B}_i}{{A} \; \mbox{tr}[{V}^{-2}]} \frac{1}{({A} + D_i)}  + o(m^{-1}).\label{ci6}
\end{align}
The other expectations can be computed similar to \eqref{ci1}, \eqref{ci5} and  \eqref{ci6}. Finally, combining all the terms, we arrive at the expression of $E[\hat c_i]$ given in \eqref{Eci}.
\end{proof}

\begin{proof}
[\textbf{Proof of Proposition \ref{cor5}}]
\label{sec:cor_prop}

We first provide the proof of (i), followed by (ii).

\textbf{[(i) Bias comparison]}
From Theorem 3 in \citet{Hirose2017} we get the bias of $\tilde{A}_i$ as
\begin{align}
    E[\tilde{A}_i - A] = \frac{2}{\mbox{tr}[V^{-2}]} l^{(1)}_{i}  + o(m^{-1}), \label{biasN}
\end{align}
where the expression of $l^{(1)}_{i}$ at $A=\hat A$ is provided in \eqref{eq:ci}.
Again, from the result (ii) of Corollary to Theorem 4 of \citet{YoshimoriLahiri2014} we get the bias of $\hat{A}_{i;\text{gls}}$ as 
\begin{align}
    E[\hat{A}_{i;\text{gls}} - A] = \frac{2}{\mbox{tr}[V^{-2}]} l^{(1)}_{i;\text{YL}}  + O(m^{-3/2}), \label{biasYL}
\end{align}
where the expression of $l^{(1)}_{i;\text{YL}}$ at $A=\hat A$ is provided in \eqref{adj_lYL}. From the expressions for $l^{(1)}_{i}$ and $l^{(1)}_{i;\text{YL}}$, we see that 
\begin{align}
    l^{(1)}_{i;\text{YL}} - l^{(1)}_i =  \frac{\mbox{tr}[V^{-2}]}{2} \; r_i + o(m^{-1}).\label{ldiff}
\end{align}
Thus using \eqref{ldiff} we get the difference of bias from \eqref{biasN} and \eqref{biasYL}, upto the order $O(m^{-1})$ as 
\begin{align} \label{bias-exp}
    E[\hat{A}_{i;\text{gls}} - A] - E[\tilde{A}_i - A] = r_i.
\end{align}
Since $(X'V^{-1}X)^{-1}$ is a positive definite matrix, the quadratic from $r_i = x_i'(X'V^{-1}X)^{-1} x_i$ is positive. Hence from \eqref{bias-exp} we have the result.

\textbf{[(ii) Length comparison]} We first note that the lengths of $I^\text{N}_i$ and $I^\text{YL}_i$ are given by Length $[I^\text{N}_i] = 2 z \Tilde{\delta}_i \equiv 2z\delta_i(\tilde{A}_i)$ and Length $[I^\text{YL}_i] = 2 z \hat \sigma_i \equiv 2 z \sigma_i(\hat{A}_{i;\text{gls}})$, respectively. To compare these lengths, we use Taylor expansions on $\tilde \delta_i$ and $\hat \sigma_i$ around $A$ respectively, and get 
\begin{align*}
    & \Tilde \delta_i - \hat \sigma_i\\
    & \equiv \delta_i(\tilde{A}_i) - \sigma_i(\hat{A}_{i;\text{gls}}) \\
    &= (\delta_i - \sigma_i) + \Big \lbrace (\tilde{A}_i - A) \delta_i' - (\hat{A}_{i;\text{gls}} - A) \sigma_i' \Big \rbrace 
     + \Big \lbrace (\tilde{A}_i - A)^2 \frac{\delta_i''}{2} - (\hat{A}_{i;\text{gls}} - A)^2 \frac{\sigma_i''}{2} \Big \rbrace \\
     & \qquad + R_{\delta} -  R_{\sigma} \\
    & = \Big \lbrace \frac{g_{2i}}{2\sqrt{g_{1i}}} + O(m^{-2}) \Big \rbrace  + \Big [ (\tilde{A}_i  - A) \Big \lbrace \sigma_i' + \frac{\partial}{\partial A} \Big ( \frac{g_{2i}}{2\sqrt{g_{1i}}} \Big ) + O(m^{-2}) \Big \rbrace - (\hat{A}_{i;\text{gls}} - A) \sigma_i' \Big ] \\ 
    &  \qquad + \Big [ (\tilde{A}_i - A)^2 \frac{1}{2} \Big \lbrace \sigma_i'' + \frac{\partial^2}{\partial A^2} \Big ( \frac{g_{2i}}{2\sqrt{g_{1i}}} \Big ) + O(m^{-2}) \Big \rbrace - (\hat{A}_{i;\text{gls}} - A)^2 \frac{\sigma_i''}{2} \Big ] \\
    & \qquad + R_{\delta} - R_{\sigma}, 
    \labthis \label{length1}
\end{align*}
where the remainder term $R_{\delta}$ is defined as 
$R_{\delta} = (\tilde A_{i} - A)^3 \tilde \delta_{i0}'''$ with $\tilde \delta_{i0}''' \equiv \delta_i'''(A)|_{A = \tilde A_0}$ and $\tilde A_0 = \gamma_1 \tilde A_{i} + (1 - \gamma_1) A $, for some $\gamma_1 \in [0,1]$. Similarly, we have $R_{\sigma} = (\hat A_{i;\text{gls}} - A)^3 \hat \sigma_{i0}'''$ with $\hat \sigma_{i0}''' \equiv \sigma_i'''(A)|_{A = \hat A_0}$, and $\hat A_0 = \gamma_2 \hat A_{i;\text{gls}} + (1 - \gamma_2) A $, for some $\gamma_2 \in [0,1]$. For more details on the calculations in \eqref{length1}, we refer the readers to similar calculations in Appendix \ref{b1011} (cf. \eqref{sigmatilde}
- \eqref{ahat_deltaprime}). We note that in this length difference proof we intend to keep terms upto the order $O_p(m^{-3/2})$, which is why we show that the difference in the remainder terms $ R_{\delta}$ and $R_{\sigma}$ is of the order $o_p(m^{-3/2})$. 
In other derivations finer calculations than the order $O_p(m^{-1})$ are not needed, so while using results from Appendix \ref{b1011}, we keep the difference of remainder terms as $O_p(m^{-3/2})$. 

Next, using Mean Value Theorem (MVT) and the relationships between $\sigma$ and $\delta$ and their derivatives (cf. Appendix \ref{b1011}), we have
\begin{align*}
    & R_{\delta} - R_{\sigma} \\
    & \equiv (\tilde A_{i} - A)^3 \tilde \delta_{i0}''' - (\hat A_{i;\text{gls}} - A)^3 \hat \sigma_{i0}'''\\
    & = (\tilde A_{i} - A)^3 \Big [\sigma_{i}'''(\tilde{A}_0) + O_p(m^{-1}) \Big ] - (\hat A_{i;\text{gls}} - A)^3 \hat \sigma_{i0}'''\\
    &= (\tilde A_{i} - A)^3 \Big [\hat \sigma_{i0}''' + O_p(m^{-1/2}) \Big ] + O_p((\tilde A_{i} - A)^3 \; m^{-1}) - (\hat A_{i;\text{gls}} - A)^3 \hat \sigma_{i0}'''\\
    & = \hat \sigma_{i0}''' \; \Big [(\tilde A_{i} - A)^3 - (\hat A_{i;\text{gls}} - A)^3 \Big ] + O_p((\tilde A_{i} - A)^3 \; m^{-1/2})  + O_p(m^{-5/2})\\
    & = \hat \sigma'''_{i0} ( \tilde A_{i} - \hat A_{i;\text{gls}}) \; \Big [ (\tilde A_{i} - A)^2  + (\tilde A_{i} - A) (\hat A_{i;\text{gls}} - A)  + (\hat A_{i;\text{gls}} - A)^2 \Big ] + O_p(m^{-2})\\ 
    & \qquad + O_p(m^{-5/2})\\
    & = O_p(m^{-2}) + O_p(m^{-5/2}), \labthis \label{remterm}
\end{align*}
where we have $\sigma'''' = O(1)$ in MVT, and  $\hat \sigma_{i0}''' = O_p(1)$, and finally all these four terms $ ( \tilde A_{i} - \hat A_{i;\text{gls}}), \; (\tilde A_i - A)^2, \; (\hat A_{i;\text{gls}} - A)^2,\; (\hat A_{i;\text{gls}} - A) (\tilde A_{i} - A)$ are of the order $O_p(m^{-1})$. Now using \eqref{remterm} in \eqref{length1}, 
we get 
\begin{align*}
 \Tilde \delta_i - \hat \sigma_i 
    &= \frac{g_{2i}}{2\sqrt{g_{1i}}} + \lbrace (\tilde{A}_i - \hat{A}) - (\hat{A}_{i;\text{gls}} - \hat{A}) \rbrace \sigma_i' + \Big \lbrace (\tilde{A}_i - A)^2 - (\hat{A}_{i;\text{gls}} - A)^2 \Big \rbrace \frac{\sigma_i''}{2} \\ 
    & \quad \quad + (\tilde{A}_i - A) \frac{\partial}{\partial A} \Big ( \frac{g_{2i}}{2\sqrt{g_{1i}}} \Big ) + o_p(m^{-3/2}). \labthis \label{length2}
\end{align*}

Next leveraging some tools in the proof of Theorem 1 of \citet{HiroseLahiri2021}, provided in Lemmas \ref{aux-lem-3} and \ref{aux-lem-4} in Appendix \ref{auxlem}, and using the relationship between ${l}^{(1)}_i$ and ${l}^{(1)}_{i;\text{YL}}$ provided in \eqref{adj_lYL}, from \eqref{length2} we have
\begin{align*}
    & \Tilde \delta_i  - \hat \sigma_i \\
    & = \frac{g_{2i}}{2\sqrt{g_{1i}}} + \frac{2}{\text{tr}[V^{-2}]} ( {l}^{(1)}_i - {l}^{(1)}_{i;\text{YL}}) \sigma_i' + \Big \lbrace \frac{2}{\text{tr}[V^{-2}]} \Big \rbrace ^2 ({l}^{(1)}_i - {l}^{(1)}_{i;\text{YL}} ) \Big \lbrace l^{(2)}_\text{RE}(A) - E [l^{(2)}_\text{RE}(A) ] \Big \rbrace \sigma_i' \\
    & \qquad + \Big \lbrace (\tilde{A}_i - A)^2 - (\hat{A}_{i;\text{gls}} - A)^2 \Big \rbrace \frac{\sigma_i''}{2} + (\tilde{A}_i - A) \frac{\partial}{\partial A} \Big ( \frac{g_{2i}}{2\sqrt{g_{1i}}} \Big ) + o_p(m^{-3/2}) \\
    &= r_i \sigma_i' - r_i \sigma_i' - \frac{2}{\text{tr}[V^{-2}]} r_i \sigma_i' \Big \lbrace - y'P^3 y + \text{tr}(P^2) \Big \rbrace + \Big \lbrace (\tilde{A}_i - A)^2 - (\hat{A}_{i;\text{gls}} - A)^2 \Big \rbrace \frac{\sigma_i''}{2} \\ 
    & \qquad + (\tilde{A}_i - A) \frac{\partial}{\partial A} \Big ( \frac{g_{2i}}{2\sqrt{g_{1i}}} \Big ) + o_p(m^{-3/2})\\
    &= \frac{2}{\text{tr}[V^{-2}]} r_i \sigma_i' \Big \lbrace y'P^3 y - \text{tr}(P^2) \Big \rbrace + \Big \lbrace  (\tilde{A}_i - A)^2 - (\hat{A}_{i;\text{gls}} - A)^2 \Big \rbrace \frac{\sigma_i''}{2} + (\tilde{A}_i - A) \frac{\partial}{\partial A} \Big ( \frac{g_{2i}}{2\sqrt{g_{1i}}} \Big ) \\ 
    & \qquad + o_p(m^{-3/2}) \labthis \label{length3},
\end{align*}
since $\frac{g_{2i}}{2\sqrt{g_{1i}}} = r_i \sigma_i'$, $\frac{\partial}{\partial A} \Big ( \frac{g_{2i}}{2\sqrt{g_{1i}}} \Big )  = O(m^{-1})$,  $ l^{(2)}_\text{RE}(A) = - y'P^3y + \frac{1}{2} \text{tr}(P^2)$ and $E [l^{(2)}_\text{RE}(A) ] = - \frac{1}{2} \text{tr}(P^2) \rbrace$. We also note that $ \Big \lbrace l^{(2)}_\text{RE}(A) - E [l^{(2)}_\text{RE}(A) ] \Big \rbrace  = O_p(\sqrt{m})$. Finally due to the fact that all the three terms in \eqref{length3} are $O_p(m^{-3/2})$, without needing to further simplify, we get our result.
\end{proof}




\newpage
\begin{proof}[\textbf{Proof of Corollary \ref{cor1}}]
\textbf{[Balanced case of the two-level model]}

We note that in the balanced case of model M2, we have $D_i = D, \; i = 1, \cdots, m$. Using this in \eqref{finalprior} and the value of $h_{i}$ for the balanced case from \eqref{adj_term}, we obtain the prior for $I_i^\text{N}$ for the balanced case as     
\begin{align*}
    \pi_{i;\text{bal}}(A) & \propto \frac{1}{(A+D)^2} (A+D) \; A 
    = \frac{A}{A+D}, 
\end{align*}
since $\mbox{tr} [V^{-2}] = m/(A+D)^2$. 
The propriety condition can be obtained following similar calculations as in the proof of Theorem \ref{Th3}.
\end{proof}

\section{Appendix}
\label{appB}

In this section, we provide the proofs of lemmas in Appendix \ref{applemmapf}, state auxiliary lemmas in Appendix \ref{auxlem} and provide some detailed derivations in Appendix \ref{app:eq}. 

\subsection{Proof of lemmas} 
\label{applemmapf}

\begin{proof}[\textbf{Proof of Lemma \ref{lem-1}}] 
To discuss the proof of \eqref{prob}, we first observe that
\begin{align*} 
    P^{\pi}(E_i^\text{U}|y)
    &= P^{\pi} \Bigg \lbrace \frac{\theta_i - \hat{\theta}^\text{BLUP}_i}{\delta_i} < \frac{\Tilde{\theta}_i - \hat{\theta}^\text{BLUP}_i}{\delta_i} + z \Big (\frac{\Tilde{\delta}_i}{\delta_i} - 1 \Big ) + z \Big |y \Bigg \rbrace \\
    &= P^{\pi} \Bigg \lbrace \frac{\theta_i - \hat{\theta}^\text{BLUP}_i}{\delta_i} < \frac{(B_i - \Tilde{B}_i)y_i + \Tilde{B}_i x_i'\Tilde{\beta} - B_i x_i'\Bar{\beta}}{\delta_i} + z \Big (\frac{\Tilde{\delta}_i}{\delta_i} - 1 \Big ) + z \Big |y \Bigg \rbrace.
\end{align*}    
Since, $\theta_i | A, y \sim N(\hat{\theta}^\text{BLUP}_i, \delta_i^2)$, we define $Z_i = (\theta_i - \hat{\theta}^\text{BLUP}_i)/\delta_i$ and get
\begin{align} 
    & P^{\pi}(E_i^\text{U}|y) = E \left[ P^{\pi} \Big ( Z_i < Q_i(z) +  z |A,y \Big ) \Big  |y \right] = E [\Phi\{Q_i(z) + z\}|y], \label{phiup}\\
    & P^{\pi}(E_i^\text{L}|y) = E \left[P^{\pi} \Big ( Z_i < Q_i(-z) - z \Big |A,y \Big ) \Big |y\right] = E\left[\Phi\{Q_i(-z) - z\}|y\right].\label{philow}
\end{align}
\end{proof}

\begin{proof}[\textbf{Proof of Lemma \ref{lem-2}}]
From Lemma \ref{lem-1}, using Taylor series expansion, we expand $\Phi(Q_i(z) + z)$ around $Q_i(z) + z = z$ or equivalently $Q_i(z) = 0$ and get 
\begin{align*}
    &\Phi \Big \lbrace Q_i(z) + z \Big \rbrace\\  
    &= \Phi(z) + \Phi'(z) Q_i(z) + \frac{1}{2} \Phi''(z) Q_i^2(z)  + \frac{1}{2} \int_{z}^{z + |Q_i(z)|} \Phi'''(x) (|Q_i(z)| + z - x)^2 dx \\
    &= 1 - \alpha/2 + \phi(z) Q_i(z) - \frac{z}{2} \phi(z) Q_i^{2}(z) + \frac{1}{2} \int_{z}^{z + |Q_i(z)|} (|Q_i(z)| + z - x)^2 (x^2 - 1) \phi(x) dx,
\end{align*}
where $\Phi(z) = 1 - \alpha/2$, $\phi'(z) = -z\phi(z)$ and $\phi''(x) = (x^2 -1 )\phi(x)$. 

Similarly, the other expansion yields 
\begin{align*}
    &\Phi \Big \lbrace Q_i(-z) - z \Big \rbrace  
    \\ 
    &= \alpha/2 + \phi(z) Q_i(-z) + \frac{z}{2} \phi(z) Q_i^{2}(-z) + \frac{1}{2} \int_{-z}^{- z + |Q_i(-z)|} (|Q_i(-z)| - z - x)^2 (x^2 - 1) \phi(x) dx, 
    \end{align*}
using  
$\Phi(-z) = \alpha/2$, $\phi(-z) = \phi(z)$ and $\phi''(-x) = (x^2 -1 )\phi(x)$. 
Thus, from \eqref{phiup}
\begin{align}\label{eu} 
    P^{\pi}(E_i^\text{U}|y) 
    &= 1 - \alpha/2 + \phi(z) \Big (T_1 - \frac{z}{2} \times T_2 \Big ) + T_3,
\end{align}
where $T_1 \equiv T_{1i} (z) = E(Q_i(z)|y), \; T_2 \equiv T_{2i}(z) = E(Q_i^{2}(z)|y)$ and the remainder term is $$T_3 \equiv T_{3i} (z) = \frac{1}{2} E \Big [ \Big \lbrace \int_{z}^{z + |Q_i(z)|} (z + |Q_i(z)| - x)^2 (x^2 - 1) \phi(x) dx \Big \rbrace |y \Big ].$$
Similarly, \eqref{philow} entails 
\begin{equation}
\begin{split} \label{el}
     P^{\pi}(E_i^\text{L}|y) 
    &= \alpha/2 + \phi(z) \Big (T_4 + \frac{z}{2} \times T_5 \Big ) + T_6,
\end{split}
\end{equation}
where $T_4 \equiv T_{4i}(z) = E(Q_i(-z)|y), \; T_5 \equiv T_{5i}(z) = E(Q_i^{2}(-z)|y)$ and $T_6$ is the remainder term defined similarly as $T_3$. Finally, combining \eqref{eu} and \eqref{el} we get
\begin{align*}
    P^{\pi}(E_i|y) 
    &= 1 - \alpha + \phi(z) \Bigg [ \Big (T_1 - T_4 \Big ) - \frac{z}{2} \Big (T_2 + T_5 \Big ) \Bigg ] + (T_3 - T_6)  .
\end{align*}
\end{proof}

\begin{proof}[\textbf{Proof of Lemma \ref{lem-3}}]
Note that, the expression $(T_1 - T_4)$ can be simplified as follows
\begin{align} \label{t1t4}
    T_1 - T_4
    &= E \Bigg [ \Bigg \lbrace \frac{k_i}{\delta_i} + z \Big (\frac{\Tilde{\delta}_i}{\delta_i} - 1 \Big ) - \frac{k_i}{\delta_i} + z \Big (\frac{\Tilde{\delta}_i}{\delta_i} - 1 \Big ) \Bigg \rbrace \Bigg | y \Bigg ] \notag \\
    &= 2z \times u,
\end{align}
where the conditional expectation $u$ is given by $u = E \Big [ \Big ({\Tilde{\delta}_i}/{\delta_i} - 1 \Big ) \Big |y \Big ]$. To evaluate conditional expectations in the following sections, we repeatedly use an important result from \cite{DattaRaoSmith2005}, given in Lemma \ref{aux-lem-1} (cf. \eqref{DRS}).  
We note that, from the expansions of $(T_1 - T_4)$ and $(T_2 + T_5)$ we intend to keep the terms of order $O_p(m^{-1})$ for $\hat{b}$ and those of order $O_p(1)$ for $\hat{b}_1$ and $\hat{b}_2$ (cf. \eqref{DRS}), as $m$ is already present in the denominator.

From \eqref{t1t4}, taking $b^u \equiv b^{u}(A) = (\Tilde{\delta}_i/\delta_i - 1 )$ in \eqref{DRS}  
and using Taylor series expansion on $\Tilde{\delta}_i$ and $\hat{\delta}_i$ around $A$ we get 
\begin{equation} \label{bhatu}
    \hat{b}^{u} \equiv b^{u}(A)|_{A=\hat{A}} = \Big (\frac{\Tilde{\delta}_i}{\hat{\delta}_i} - 1 \Big ) 
    =  
    (\Tilde{A}_i - \hat{A}) \frac{\sigma'_i}{\sigma_i} + o_p(m^{-1}),
\end{equation}
which is a function of both $\Tilde{A}_i$ and $\hat{A}$. Finally, using the asymptotic properties of $\Tilde{A}_i,\; \hat{A}$ and the relationship between these given in Lemma \ref{aux-lem-3}, we get $\hat{b}^{u}$ upto the order $O_p(m^{-1})$ as
\begin{align} 
    \hat{b}^{u} = 
    \frac{\hat l^{(1)}_{i}}{\mbox{tr}[\hat{V}^{-2}]} \times \frac{\hat{B}_i}{\hat{A}}. \label{bbhatu}   
\end{align}
For more detailed calculation regarding equations \eqref{bhatu} and \eqref{bbhatu}, refer to Appendix \ref{b1011}.  

We next evaluate $\hat{b}^{u}_1$ upto the order $O_p(1)$ and get
\begin{equation*}
    \hat{b}^{u}_1 \equiv \frac{\partial}{\partial A}{b}^{u}(A)|_{A=\hat{A}} = - \frac{\Tilde{\sigma}_i \hat{\sigma}_i}{2\hat{A}^2}. 
\end{equation*}
In the above $\Tilde{\sigma}_i \hat{\sigma}_i$ can be replaced with $\hat{\sigma}_i^2$ using Taylor series expansion (cf. \eqref{tildehat}) and thus we get $\hat{b}^{u}_1$, upto the order $O_p(1)$, as follows
\begin{equation} \label{b1u}
    \hat{b}^{u}_1 =  - \frac{\hat{\sigma}_i^2}{2\hat{A}^2} 
    = - \frac{\hat{B}_i}{2 \hat{A}}.
\end{equation}
Similarly, for the second derivative $\hat{b}^{u}_2$, upto the order $O_p(1)$ we get
\begin{equation} \label{b2u}
    \hat{b}^{u}_2 = - \frac{\hat{B}_i}{2\hat{A}^2}  \Big ( \frac{\hat{B}_i}{2} - 2\Big ).
\end{equation}
For the preceding calculations, details including derivatives are deferred to Appendix \ref{b1213}. Note that, for obtaining a solution for $\pi(A)$, we need to express everything in terms of $\hat{A}$, but we write the intermediate steps in terms of $\hat{B}_i$ only for brevity. 
Thus, using \eqref{b1u} and \eqref{b2u}  and \eqref{k23} from Lemma \ref{aux-lem-2}, the second term in \eqref{DRS} 
can be worked out, upto $O_p(m^{-1})$, as
\begin{align} \label{bhatu1}
    \frac{1}{2m\hat{k}_2} \left(\hat{b}^u_2 - \frac{\hat{k}_3}{\hat{k}_2} \hat{b}^u_1\right)
    &= -\frac{\hat{B}_i}{2 \hat{A} \; \mbox{tr}[\hat{V}^{-2}]}\Bigg \lbrace \Big (\frac{\hat{B}_i}{2} - 2 \Big ) \frac{1}{\hat{A}} + 4 \frac{\mbox{tr}[\hat{V}^{-3}]}{\mbox{tr}[\hat{V}^{-2}]} \Bigg \rbrace. 
\end{align}
The third term in \eqref{DRS}, which involves the prior $\pi$, upto the order $O_p(m^{-1})$ is
\begin{align} \label{bhatu2}
\frac{1}{m\hat{k}_2} \hat{b}^u_1 \hat{\rho}_1 =
\frac{2}{\mbox{tr}[\hat{V}^{-2}]} \times \Bigg (- \frac{\hat{B}_i}{2\hat{A}} \Bigg )\hat{\rho}_1 =-\frac{1}{\mbox{tr}[\hat{V}^{-2}]} \frac{\hat{B}_i}{\hat{A}} \hat{\rho}_1.
\end{align}
Finally, from \eqref{bbhatu}, \eqref{bhatu1} and \eqref{bhatu2}, we get $(T_1 - T_4)$, upto the order $O_p(m^{-1})$, as follows
\begin{align*}
    T_1 - T_4 &= 2z \times u  \\
    &= 2z \frac{\hat{B}_i}{\hat{A}\mbox{tr}[\hat{V}^{-2}]} \Bigg [ \hat{l}^{(1)}_{i} - 
    \frac{1}{2}\Bigg \lbrace \Big (\frac{\hat{B}_i}{2} - 2 \Big ) \frac{1}{\hat{A}} + 4 \frac{\mbox{tr}[\hat{V}^{-3}]}{\mbox{tr}[\hat{V}^{-2}]} \Bigg \rbrace - \hat{\rho}_1 \Bigg ]. 
\end{align*}
\end{proof}

\begin{proof}[\textbf{Proof of Lemma \ref{lem-4}}] 

$(T_2 + T_5)$ can be further simplified as follows
\begin{align*} \label{t2t5}
    T_2 + T_5
    &= E [Q^2_i(z) |y] + E [Q^2_i(-z)|y] \notag \\
    &= E \Bigg [ \Bigg \lbrace {\frac{k_i}{\delta_i} + z \Big (\frac{\Tilde{\delta}_i}{\delta_i} - 1 \Big )\Bigg \rbrace }^2 
    \notag + 
    \Bigg \lbrace  {\frac{k_i}{\delta_i} - z \Big (\frac{\Tilde{\delta}_i}{\delta_i} - 1 \Big ) \Bigg \rbrace}^2 \Bigg | y \Bigg ] \notag \\
    &= 2v + 2w, \labthis
\end{align*}
where the conditional expectations are $v = E \Big [ \left\{ k_i/\delta_i 
\right\}^2 \Big |y  \Big ] $ and $w = E \Big [ \left\{ z\Big (\Tilde{\delta}_i/\delta_i - 1 \Big )\right\}^2 \Big |y \Big ]$. 

Firstly, we work with the conditional expectation $w$ in \eqref{t2t5}, as we observe that $w$ involves the square of $b^u = (\Tilde{\delta}_i/\delta_i - 1 )$, which we already worked out in the proof of Lemma \ref{lem-3}. Hence, we have $ b^w \equiv b^w(A) = (\Tilde{\delta}_i/\delta_i - 1 )^2 $ and use the following relationship of derivatives
\begin{gather} \label{der}
    b^w = ({b^u})^2, \quad 
    b^w_1 = 2 b^u b^u_1, \quad
    b^w_2 = 2 [b^u b^u_2 + {(b^u_1)}^2].
\end{gather}
Using \eqref{der} and applying the ${C_r}$-inequality, we conclude that $\hat{b}^w$ and $\hat{b}^w_1$ are ignorable in the context of using in \eqref{DRS}, as described in the proof of Lemma \ref{lem-3}, i.e.,
\begin{equation}
\begin{split}
    \hat{b}^w = o_p(m^{-1}) \quad \text{and} \quad 
    \hat{b}^w_1 = O_p(m^{-1}).
\end{split}\label{bhatw}
\end{equation}
Detailed derivations of \eqref{bhatw} are deferred to Appendix \ref{b1920}. Again, we see that for $\hat{b}^w_2$ the contributing term to the order $O_p(1)$ is ${(\hat{b}^u_1)}^2$. Therefore, using \eqref{b1u}, $\hat{b}^w_2$ approximated upto the order $O_p(1)$ is
\begin{equation*}
    \hat{b}^w_2 = 2 {(\hat{b}^u_1)}^2 = \frac{\hat{\sigma}^4_i}{2 \hat{A}^4} = \frac{\hat{B}^2_i}{2 \hat{A}^2} = \frac{D_i^2}{2 \hat{A}^2 (\hat{A} + D_i)^2}.
\end{equation*}
Hence, the contribution of $w$ in \eqref{t2t5}, upto the order $O_p(m^{-1})$, involves only part of the second term of \eqref{DRS}. To finally use in \eqref{stoch}, we multiply by $-\frac{z}{2}(2z^2) = -z^3$ and get
\begin{equation} \label{wexp}
     -z^3 \frac{1}{2m\hat{k}_2} \hat{b}^w_2 =  -z^3\frac{1}{\mbox{tr}[\hat{V}^{-2}]} \frac{\hat{B}^2_i}{2 \hat{A}^2}.
\end{equation}

Secondly, we work with the conditional expectation $v$ in \eqref{t2t5}. To begin with we take $b^{v}(A) = k_i/\delta_i$, i.e., at first we do not consider the square and use the relation provided in \eqref{der}. We compute, one by one, $\hat{b}^v, \hat{b}_1^v, \hat{b}_2^v$. Now, $\hat{b}^v$ is
\begin{align} \label{bvhat}
    \hat{b}^{v} = 
    \frac{(\hat{B}_i - \Tilde{B}_i)}{\hat{\delta}_i}y_i + \frac{\Tilde{B}_i x_i'\Tilde{\beta} - \hat{B}_i x_i'\hat{\beta}}{\hat{\delta}_i}.
\end{align}
From the first part of $\hat{b}^{v}$, using Taylor series expansion and similar logic as $\hat{b}^u$, we get that the order of $(\hat{B}_i - \tilde{B}_i)$ is $O_p(m^{-1})$ and that of $1/\hat{\delta}_i$ is $O_p(1)$. Thus, $\dfrac{(\hat{B}_i - \Tilde{B}_i)}{\hat{\delta}_i}y_i$ is $O_p(m^{-1})$. We work with the second part of $\hat{b}^{v}$, for which we expand $\Tilde{B}_i x_i'\Tilde{\beta}$ and $\hat{B}_i x_i'\hat{\beta}$, respectively around $A$, using Taylor series, and subtract to get
\begin{align}
    &\Tilde{B}_i x_i'\Tilde{\beta} - \hat{B}_i x_i'\hat{\beta} \notag \\
    &= (\tilde{A}_i - \hat{A}) \frac{\partial (B_i x_i'\bar{\beta})}{\partial A}  + \Big \lbrace  \frac{(\tilde{A}_i - A)^2 - (\hat{A} - A)^2}{2} \Big \rbrace \frac{\partial^2 (B_i x_i'\bar{\beta})}{\partial A^2}
    + o_p(m^{-\frac{3}{2}}).
\end{align}
Note that, $\frac{\partial}{\partial A} (B_i x_i'\bar{\beta}) = O(m^{-1/2})$. Hence, using earlier results, we get the order of $(\Tilde{B}_i x_i'\Tilde{\beta} - \hat{B}_i x_i'\hat{\beta})$ as $O_p(m^{-3/2})$. Thus, $\hat{b}^{v}$ is of order $O_p(m^{-1})$. We then consider the square, i.e., 
\begin{align*}
    b^{s} \equiv (b^{v})^2 = \Big \lbrace \frac{(B_i - \Tilde{B}_i)y_i + \Tilde{B}_i x_i'\Tilde{\beta} - B_i x_i'\Bar{\beta}}{\delta_i} \Big \rbrace^2.
\end{align*} 
We need $\hat b^s \equiv (\hat{b}^{v})^2$, which is of smaller order, so in the sense of \eqref{DRS}, $\hat b^s$ becomes ignorable. After this we compute the first and second derivatives of $b^s$. From \eqref{DRS} we see that for the first derivative, we need $\hat{b}^{s}_1 = 2 \hat{b}^{v} \hat{b}^{v}_1$ upto the order $O_p(1)$. We compute the first derivative as follows
\begin{align} \label{bv1}
b^{v}_1 &= \frac{(y_i - x_i'\bar{\beta}) }{\delta_i} \frac{\partial B_i }{\partial A} + \Big \lbrace (\Tilde{B}_i x_i'\Tilde{\beta} - B_i x_i'\bar{\beta}) + (B_i - \Tilde{B}_i) y_i \Big \rbrace \frac{\partial}{\partial A} \Big (\frac{1}{\delta_i} \Big ) - \frac{B_i}{\delta_i}\frac{\partial (x_i'\bar{\beta})}{\partial A}  
\end{align}
and $\hat{b}^{v}_1$, upto $O_p(1)$ is 
\begin{align} \label{bvhat1}
\hat{b}^{v}_1
    &= - (y_i - x_i'\hat{\beta}) \hat{\sigma}_i \frac{\hat{B}_i}{\hat{A}D_i}.
\end{align}
The detailed calculations of (\ref{bv1}), (\ref{bvhat1}) and (\ref{b2vhat}) are given in Appendix \ref{b242526}. Then, using \eqref{bvhat} and \eqref{bvhat1}, we get $\hat{b}^{s}_1$ as 
\begin{align*}
    \hat{b}^{s}_1 &\equiv \frac{\partial}{\partial A} (b^{v}(A))^2\Big |_{ A = \hat{A}} = 2 \hat{b}^{v} \hat{b}^{v}_1 = O_p(m^{-1}).
\end{align*} 
Thus, with similar arguments $\hat{b}^{s}_1$ becomes ignorable, as we are interested upto the order $O_p(1)$. Coming to the second derivative, we see that 
\begin{align} \label{b2vhat}
    \hat{b}_2^v = O_p(1),
\end{align} 
so that from $\hat{b}^s_2 = 2[\hat{b}^v \hat{b}_2^v + (\hat{b}^{v}_1)^2]$, only $2(\hat{b}^{v}_1)^2$ contributes to the order $O_p(1)$. Hence, we ignore $2 \hat{b}^{v} \hat{b}^{v}_2$ and using \eqref{bvhat1}, we get $\hat{b}^s_2$ upto the order $O_p(1)$ as
\begin{align} \label{B.21}
    2(\hat{b}^{v}_1)^2 = 
    2 (y_i - x_i'\hat{\beta}) ^2 \frac{\hat{\sigma}_i^2 \hat{B}_i^2}{\hat{A}^2 D_i^2}.
\end{align}
Putting \eqref{B.21} in the form of \eqref{DRS}, we thus need $ \Big \lbrace 2(\hat{b}^{v}_1)^2 \Big \rbrace / 2m\hat{k}_2$. Detailed derivations of (\ref{bv1}), (\ref{bvhat1}) and (\ref{b2vhat}) are deferred to Appendix \ref{b242526}.

Going back to \eqref{stoch}, we need 
$-z (T_2 + T_5)/2  
= -z (v + w)$ and using \eqref{B.21} get the first term upto order $O_p(m^{-1})$ as follows 
\begin{equation} \label{vexp}
     - z v = -z \frac{1}{2m\hat{k}_2} \Big \lbrace 2(\hat{b}^{v}_1)^2 \Big \rbrace = 
     - \frac{2z}{\mbox{tr}[\hat{V}^{-2}]} (y_i - x_i'\hat{\beta}) ^2 \frac{\hat{\sigma}_i^2 \hat{B}_i^2}{\hat{A}^2 D_i^2}.
\end{equation}
Hence, the final contribution of $ -\frac{z}{2} (T_2 + T_5)$ upto $O_p(m^{-1})$ in \eqref{stoch}, combining \eqref{wexp} and \eqref{vexp}, is
\begin{align*}
    -\frac{z}{2} (T_2+T_5) 
     &= \frac{-2 z \hat{B}_i}{\mbox{tr}[\hat{V}^{-2}]\hat{A}} \times \frac{\hat{B}_i}{2\hat{A}} \times \Bigg[ 2 (y_i - x_i'\hat{\beta})^2 \frac{\hat{\sigma}_i^2}{D_i^2} + \frac{z^2}{2} \Bigg ] \notag \\
     &= \frac{-2 z \hat{B}_i}{\hat{A} \; \mbox{tr}[\hat{V}^{-2}]} \Bigg[ \Bigg ( \frac{y_i - x_i'\hat{\beta}}{\hat{A} + D_i} \Bigg )^2 + \frac{z^2 \hat{B_i}}{4 \hat{A}} \Bigg ]. 
\end{align*}
\end{proof}

\begin{proof}[\textbf{Proof of Lemma \ref{lem-5}}] 
We need to show that the remainder terms $T_3$ and $T_6$ are of the order $o_p(m^{-1})$ and hence ignorable. Note that, similar result involving remainder terms are shown in \cite{HLthesis} [Proof of Theorem 4.1.], 
\cite{CLL2008} [Proof of Theorem 3.1] and \cite{YoshimoriLahiri2014} [Proof of Theorem 1], the difference in our case being that the expectation in our remainder terms $T_3$ and $T_6$ are conditional on $y$. We prove the result for the first term $T_3$, i.e., 
$$ T_3 = E \Big [ \Big \lbrace \int_{z}^{z + |Q_i(z)|} (z + |Q_i(z)| - x)^2 (x^2 - 1) \phi(x) dx \Big \rbrace |y \Big ] = o_p(m^{-1})$$ 
and the other follows similarly. Let us denote $Q_i \equiv Q_i(z) $ for notational simplicity. For, $x \in (z,z+ |Q_i|)$, we have 
\begin{align}\label{in1}
0 < |z+|Q_i|-x| < |Q_i|
\end{align}
We also note that, as the function $(x^2 -1)\phi(x)$ is maximized at $x = \sqrt{3}$, we have 
\begin{align}\label{in2}
    |(x^2 -1)\phi(x)| \le 2\phi(\sqrt{3}).
\end{align}
Using inequalities in \eqref{in1} and \eqref{in2} along with triangle inequality we get that
\begin{align*}
    0 \; \le \; & T_3 
     \le 2 \phi(\sqrt{3}) \; E \Big [ Q_i^2 \int_{z}^{z + |Q_i|} dx \Big | y \Big ] = C E [ |Q_i|^3 | y ] ,
\end{align*}
where $C$ stands for a generic positive constant free of $A$, as noted earlier. 
Hence, to achieve our goal, it will be enough to find the order of $E[|Q_i|^3|y]$. To this end, we first show that $E[Q_i^4|y] = O_p(m^{-2})$, then using Lyapunov's inequality we get that $E[|Q_i|^3|y] = O_p(m^{-3/2})$. 

Recall that, $Q_i = k_i/\delta_i + z \Big (\Tilde{\delta}_i/\delta_i - 1 \Big )$. From $C_r$-inequality we know that 
\begin{align*}
    0 \le Q_i^4 &\le 2^{4-1} \Bigg [ \Big \lbrace \frac{k_i}{\delta_i} \Big \rbrace^4 + \Big \lbrace z \Big (\frac{\Tilde{\delta}_i}{\delta_i} - 1 \Big ) \Big \rbrace^4 \Bigg ] = 8 [b^x + b^y], 
\end{align*}
where $b^x =  \lbrace k_i/\delta_i \rbrace^4$ and $b^y = \Big \lbrace z \Big (\Tilde{\delta}_i/\delta_i - 1 \Big ) \Big \rbrace^4$. In the proofs of Lemmas \ref{lem-3} and \ref{lem-4}, assuming $b^u =  ( \Tilde{\delta}_i/\delta_i - 1 )$ and $b^v = k_i/\delta_i$, we get the following results
\begin{align*}
    &\hat{b}^u = O_p(m^{-1}), \; \hat{b}^u_1 = O_p(1) \text{ and } \hat{b}^u_2 = O_p(1), \\
    &\hat{b}^v = O_p(m^{-1}), \; \hat{b}^v_1 = O_p(1) \text{ and } \hat{b}^v_2 = O_p(1).
\end{align*}

Using the relationship of derivatives given in \eqref{der} (but for power $4$), we also get the following
\begin{gather*}
    \hat{b}^x = O_p(m^{-4}), \; \hat{b}^x_1 = O_p(m^{-3}) \text{ and } \hat{b}^x_2 = O_p(m^{-2}),\\
    \hat{b}^y = O_p(m^{-4}), \; \hat{b}^y_1 = O_p(m^{-3}) \text{ and } \hat{b}^y_2 = O_p(m^{-2}).
\end{gather*} 
Using the result given in \eqref{DRS} in Lemma \ref{aux-lem-1}, we get that $E[Q_i^4|y] = O_p(m^{-2})$. Thus, from Lyapunov's inequality we conclude that $E[|Q_i|^3|y] \le \lbrace E[Q_i^4|y] \rbrace ^{3/4}$,
which implies $T_3 = O_p(m^{-3/2})$.
Similarly, it can be shown that $T_6 = O_p(m^{-3/2})$. Hence, the difference of the remainder terms $T_3$ and $T_6$ is of the order $o_p(m^{-1})$. 
\end{proof}

\subsection{Auxiliary lemmas}
\label{auxlem}
\begin{lemma} \textbf{[\cite{DattaRaoSmith2005}, Equation (21)]} \label{aux-lem-1}
For any function $b(A)$ of $A$, for which the below expectation exists,
we have
\begin{equation} \label{DRS}
    E \lbrace b(A)|y \rbrace = \hat{b} + \frac{1}{2m\hat{k}_2} \left(\hat{b}_2 - \frac{\hat{k}_3}{\hat{k}_2}\hat{b}_1\right) + \frac{\hat{b}_1}{m\hat{k}_2}\hat{\rho}_1 + O_p(m^{-2}),
\end{equation}
where ($m\times k$) is the negative logarithm of the REML likelihood $L_\text{RE}$ which yields
\begin{align*}
    k &\equiv k_\text{RE}(A)\\ &= \frac{1}{2m} \log(|V|) + \frac{1}{2m} \log(|X'V^{-1}X|) + \frac{1}{2m}y'Py \\ &= - \frac{1}{m} \log \lbrace L_\text{RE}(A) \rbrace.
\end{align*}
Here $\hat A$ is the REML estimator of $A$ that minimizes $k$ and with $D^j = \dfrac{\partial^j}{\partial A^j}$, we have $\hat{b} \equiv b(\hat{A}), \; \hat{b}_j \equiv D^{j}b(A)|_{A = \hat{A}}, j =1,2$ and $\hat{k}_j \equiv D^{j}k(A)|_{A = \hat{A}}, j = 2,3$ and $\hat \rho _1 = D \rho(A)|_{A = \hat A}$, where $\rho(A) = \log \pi (A)$. 
For more details on higher-order expansions of conditional expectations, which follow from expansion of ratio of integrations, we refer the readers to \cite{lindley1980approximate} and \cite{kass1989approximate}.
\end{lemma}

\begin{lemma} \textbf{[\cite{DattaRaoSmith2005}, Page 193, Section 7]} \label{aux-lem-2}
$\hat{k}_2, \hat{k}_3$ can be approximated, upto the order of $O_p(1)$ as 
$$\hat{k}_2 = \frac{\mbox{tr}[\hat{V}^{-2}]}{2m} + o_p(1) \quad \mbox{and} \quad
     \hat{k}_3 = - \frac{2\mbox{tr}[\hat{V}^{-3}]}{m} + o_p(1),$$
which yields, 
\begin{gather}
     \frac{\hat{k}_3}{\hat{k}_2} = - 4 \frac{\mbox{tr}[\hat{V}^{-3}]}{\mbox{tr}[\hat{V}^{-2}]} + o_p(1). \label{k23}
\end{gather}
\end{lemma}

\begin{lemma} \textbf{[\cite{HiroseLahiri2021}, Theorem 1 and \cite{YoshimoriLahiri2014} Theorem 2]} \label{aux-lem-3}
Under regularity conditions, relationship between $\Tilde{A}_i$ and $\hat{A}$ is given by
\begin{equation} \label{eq:Arel}
    \Tilde{A}_i - \hat{A} = \frac{2 {l}^{(1)}_{i;\text{YL}}}{\mbox{tr}[V^{-2}]} + o_p(m^{-1}),
\end{equation}
where $\mbox{tr}[V^{-2}] = \sum_{u=1}^{m} (A+D_u)^{-2}$ and ${l}_{i;\text{YL}}^{(1)} = \dfrac{\partial}{\partial A} \log h_{i;YL}(A) = O(1)$. 
Thus, we have $\Tilde{A}_i - \hat{A} = O_p(m^{-1})$. 
\end{lemma}

\begin{lemma} \textbf{
[\cite{HiroseLahiri2021}, Appendix A.1: Proof of Theorem 1]}  \label{aux-lem-4}
   For a general class of estimator $\hat{A}_{i;G}$ of $A$ based on adjusted likelihood, it holds that $\hat{A}_{i;G} - A =  O_p(m^{-1/2})$. As $\Tilde{A}_i$ is a special case of $\hat{A}_{i;G}$, it follows that $\Tilde{A}_{i} - A =  O_p(m^{-1/2})$. Similarly, as $\hat{A}$ is a consistent estimator of $A$, we have $\hat{A} - A = O_p(m^{-1/2})$. 
   Again, $\hat{A}_{i;G} - \hat{A}$ yields
\begin{align*}
    \hat{A}_{i;\text{G}} - \hat{A}
    &= \frac{2}{\text{tr}[V^{-2}]} l^{(1)}_{i;\text{G}} + \Big \lbrace \frac{2}{\text{tr}[V^{-2}]} \; \Big \rbrace ^2 l^{(1)}_{i;\text{G}} (A) \; \Big \lbrace l^{(2)}_\text{RE}(A) - E [l^{(2)}_\text{RE}(A) ] \Big \rbrace \\
    & \qquad + \frac{1}{2}\Big \lbrace \frac{2}{\text{tr}[V^{-2}]} \Big \rbrace ^3 \; \Big \lbrace l^{(1)}_{i;\text{G}} (A) \; (l^{(1)}_{i;\text{G}}(A) + 2 l^{(1)}_\text{RE}(A)) \Big \rbrace \Big \lbrace l^{(3)}_{i;\text{G}} + o_p(m) \Big \rbrace,
\end{align*}
where $l^{(1)}_\text{RE}$ and $l^{(2)}_\text{RE}$ are the first and second derivatives of the REML log-likelihood, and $l^{(1)}_{i;\text{G}}$ and $l^{(3)}_{i;\text{G}}$ are the first and third derivatives of the logarithm of a general adjustment factor $h_{i;G}(A)$, respectively.
\end{lemma}

\subsection{Detailed calculation}
\label{app:eq}
In this section, we provide detailed derivations necessary to prove the main results in the paper.

\subsubsection{Equations (\ref{bhatu}) and (\ref{bbhatu}).} \label{b1011}

We use Taylor series expansion on $\Tilde{\delta}_i$ and $\hat{\delta}_i$, respectively, around $A$ to get the following:
\begin{align}
    &\Tilde{\delta}_i
    = \delta_i + (\Tilde{A}_i - A) \delta'_i + (\Tilde{A}_i - A)^2 \frac{\delta''_i}{2} + O_p(m^{-3/2}), \label{sigmatilde} \\
    &\hat{\delta}_i
    = \delta_i + (\hat{A} - A) \delta'_i + (\hat{A} - A)^2 \frac{\delta''_i}{2} + O_p(m^{-3/2}) \label{sigmahat},
\end{align}
where $\delta'_i$ and $\delta''_i$ are respectively the first and second order derivatives of $\delta_i$ with respect to $A$.
Dividing \eqref{sigmatilde} by \eqref{sigmahat} and subtracting $1$ from the ratio, we get 
\begin{equation*}
\begin{split}
    \hat{b}^{u} &= \Big ( \frac{\Tilde{\delta}_i}{ \hat{\delta}_i} - 1 \Big ) \\
    &= \frac{ 
    (\Tilde{A}_i - \hat{A})\delta'_i +
    (\Tilde{A}_i - A)^2 \frac{\delta''_i}{2} - (\hat{A} - A)^2 \frac{\delta''_i}{2} + O_p(m^{-3/2})}{\delta_i + (\hat{A} - A) \delta'_i + o_p(m^{-1/2})}.
\end{split}
\end{equation*}
Now using negative binominal expansion, we express $\delta_i$ in terms of $\sigma_i$ as
\begin{align*}
    \delta_i 
    &= \sqrt{g_{1i}} \Big (1 + \frac{g_{2i}}{g_{1i}} \Big )^{1/2} \\
    &=  \sqrt{g_{1i}} \Big \lbrace 1+\frac{g_{2i}}{2 g_{1i}} + O(m^{-2}) \Big \rbrace \\
    &= \sigma_i + \frac{g_{2i}}{2 \sqrt{g_{1i}}}  + O(m^{-2}) \\
    &= \sigma_i + O(m^{-1}),\labthis \label{sig_del}
\end{align*}
which yields $\delta_i'$ and $\delta_i''$ as 
\begin{equation}
    \begin{split}
    & \delta'_i = \frac{\partial \sqrt{g_{1i}}}{\partial A} + \frac{1}{2}    \frac{\partial}{\partial A} \Big ( \frac{g_{2i}}{\sqrt{g_{1i}}} \Big ) + O(m^{-2}) = \sigma_i' + O(m^{-1}),\\
    & \delta''_i = \frac{\partial^2 \sqrt{g_{1i}}}{{\partial A^2}}
    + \frac{1}{2}    \frac{\partial^2}{\partial A^2}
    \Big ( \frac{g_{2i}}{\sqrt{g_{1i}}} \Big ) + O(m^{-2}) = \sigma_i'' + O(m^{-1}),
    \end{split}\label{sig_del_pr}
\end{equation}
since $\dfrac{\partial}{\partial A} \Big ( \dfrac{g_{2i}}{\sqrt{g_{1i}}} \Big )$ and $ \dfrac{\partial ^2}{{\partial A^2}}
    \Big ( \dfrac{g_{2i}}{\sqrt{g_{1i}}} \Big )$ are $O(m^{-1})$. Note that $(\hat{A} - A)= O_p(m^{-1/2})$, hence we get 
\begin{align*}
    &(\hat{A} - A) \delta'_i 
    = (\hat{A} - A) \sigma_i' + O_p(m^{-3/2}), \labthis \label{ahat_delta}\\
    &(\hat{A} - A)^2 \delta''_i 
    = (\hat{A} - A)^2\sigma_i'' + O_p(m^{-2}).\labthis \label{ahat_deltaprime}
\end{align*}

Thus, $\hat{b}^{u}$ 
simplifies to
{\allowdisplaybreaks
\begin{align*}
\hat{b}^{u}
&= \frac{
 (\Tilde{A}_i - \hat{A})\sigma_i'
+ \Big \lbrace (\Tilde{A}_i - A)^2 \frac{\sigma_i''}{2} - (\hat{A} - A)^2 \frac{\sigma_i''}{2} \Big \rbrace + O_p(m^{-3/2}) 
}{\delta_i +(\hat{A} - A)\delta_i' + o_p(m^{-1/2})}\\
&= \frac{\sqrt{g_{1i}}}{\sqrt{g_{1i} + g_{2i}}} \times
\frac{\Big \lbrace \Big 
(\Tilde{A}_i - \hat{A}) \frac{\partial \sqrt{g_{1i}}}{\partial A} \frac{1}{\sqrt{g_{1i}}} + O_p(m^{-3/2}) 
\Big \rbrace}{\Big \lbrace 1 +(\hat{A} - A)\frac{\delta_i'}{\delta_i} + o_p(m^{-1/2}) \Big \rbrace} \\
&= \Big \lbrace 1 - \frac{g_{2i}}{2g_{1i}}  + o_p(m^{-1}) \Big \rbrace \Big \lbrace 
(\Tilde{A}_i - \hat{A}) \frac{\sigma_i'}{\sigma_i} + O_p(m^{-\frac{3}{2}}) + O_p(m^{-2}) \Big \rbrace  \\
& \qquad \times \Big \lbrace 1 - (\hat{A} - A) \frac{\delta_i'}{\delta_i} + o_p(m^{-\frac{1}{2}}) \Big \rbrace\\
&= \lbrace 1 + O_p(m^{-1}) \rbrace \Big \lbrace 
(\Tilde{A}_i - \hat{A}) \frac{\sigma_i'}{\sigma_i} + O_p(m^{-3/2}) + O_p(m^{-2}) \Big\rbrace \lbrace 1 + O_p(m^{-1/2}) \rbrace ,
\end{align*}}
using negative binomial expansion and the facts that $g_{1i},\delta_i,\;\delta_i',\; \delta_i'/\delta_i, \;\sigma_i,\;\sigma_i',\;\sigma_i'/\sigma_i,\; \sigma_i''$ are all $O(1)$, $g_{2i}= O(m^{-1})$, $(\Tilde{A}_i - \hat{A}) = O_p(m^{-1})$ and $(\Tilde{A}_i - A), \; (\hat{A} - A)$ are $O_p(m^{-1/2})$. 
Thus, ignoring lower order terms we get $\hat{b}^{u}$ upto $O_p(m^{-1})$ as
$\hat{b}^{u}= 
(\Tilde{A}_i - \hat{A}) {\sigma_i'}/{\sigma_i}$.

We then compute $\sigma_i'/\sigma_i$ as follows
\begin{align*}
    \frac{\sigma_i'(A)}{\sigma_i(A)} 
    &= \sqrt{ \frac{A + D_i}{AD_i}} \times \frac{\partial}{\partial A} \Big \lbrace \sqrt{ \frac{AD_i}{(A+D_i)}} \Big \rbrace \\
     &= \sqrt{ \frac{A + D_i}{AD_i}} \times \frac{D_i \sqrt{D_i}}{2 \sqrt{A} \sqrt{A+D_i} (A+D_i)} \\
     &= \frac{B_i}{2A}.
\end{align*}
Thus, $\hat{b}^u$ upto order $O_p(m^{-1})$ is
\begin{align*}
    \hat{b}^u &= 
    (\tilde{A}_i - \hat{A}) \frac{\hat{B}_i}{2\hat{A}} = 
    2 \frac{l^{(1)}_{i}}{\mbox{tr}[\hat{V}^{-2}]} 
 \frac{\hat{B}_i}{2\hat{A}} 
 = 
 \frac{l^{(1)}_{i}}{\mbox{tr}[\hat{V}^{-2}]} 
 \times \frac{D_i}{\hat{A}(\hat{A}+D_i)}.
\end{align*}

\subsubsection{Equations (\ref{b1u}) and (\ref{b2u}).} \label{b1213}

Note that, 
\begin{align*}
    b_1^u 
    &= - \frac{\Tilde{\delta_i}}{2} (g_{1i} + g_{2i})^{-3/2} \Big \lbrace B_i^2 - 2 \frac{B_i^2}{A+D_i} r_i + B_i^2 x_i'(X'V^{-1}X)^{-1} X'V^{-2}X (X'V^{-1}X)^{-1} x_i \Big \rbrace.
\end{align*}
From regularity condition R4, in \eqref{b1u}, ignoring lower order terms, we obtain the first derivative upto $O(1)$ as 
\begin{equation*}
    b_1^u = - \frac{\Tilde{\sigma}_i}{2} (g_{1i})^{-3/2} B_i^2
    = - \frac{\Tilde{\sigma}_i \sigma_i}{2A^2}.
\end{equation*}
For the approximation in \eqref{b1u}, to replace $\tilde{\sigma}_i \hat{\sigma}_i$ with $\hat{\sigma}^2_i$, we use the following argument
\begin{align*}
     \Tilde{\sigma}_i \hat{\sigma}_i &= 
     \Big [\sigma_i + (\Tilde{A}_i - A)\sigma'_i + (\Tilde{A}_i - A)^2 \frac{\sigma''_i}{2} + o_p(m^{-1})  \Big ] \hat{\sigma}_i \\ 
     &= \sigma_i \hat{\sigma}_i + O_p(m^{-1/2}) \\
     &= \hat{\sigma}^2_i + o_p(1), \labthis \label{tildehat}
\end{align*}
using the fact that $\hat{\sigma}_i = \sigma_i + o_p(1)$. Thus, we get $\hat{b}^{u}_1$ upto the order $O_p(1)$ as:
\begin{align*}
    \hat{b}^{u}_1 =  - \frac{\hat{\sigma}_i^2}{2\hat{A}^2} 
    = - \frac{\hat{B}_i}{2 \hat{A}} = - \frac{D_i}{2\hat{A}(\hat{A} + D_i)}.
\end{align*}

For \eqref{b2u}, after computing the second derivative and using previous argument we get 
\begin{equation*}
    \hat{b}^{u}_2 \equiv \frac{\partial}{\partial A}{b}^{u}_1(A)|_{A=\hat{A}}
    = - \frac{\hat{\sigma}_i^2}{2\hat{A}^3}  \Big ( \frac{\hat{B}_i}{2} - 2\Big ) + o_p(1).
\end{equation*}
Finally, we get $\hat{b}^{u}_2$ upto the order $O_p(1)$ as
\begin{align*}
    \hat{b}^{u}_2 = - \frac{\hat{B}_i}{2\hat{A}^2}  \Big ( \frac{\hat{B}_i}{2} - 2\Big )
    = \frac{D_i (4\hat{A} + 3D_i)}{4\hat{A}^2(\hat{A}+D_i)^2}.
\end{align*}

\subsubsection{Equation (\ref{bhatw}).} \label{b1920}
Since $\hat{b}^u = O_p(m^{-1})$ we get $\hat{b}^w = (\hat{b}^u)^2 = O_p(m^{-2})$
and hence, ignorable in this context. For $\hat{b}^w_1$, we use \eqref{b1u} and \eqref{der} and from similar arguments as before, we get that $\hat{b}^w_1$ is ignorable, i.e., $\hat{b}^w_1 = O_p(m^{-1})O_p(1) = O_p(m^{-1})$.

\subsubsection{Equations (\ref{bv1}), (\ref{bvhat1}) and (\ref{b2vhat}).} \label{b242526}

For \eqref{bv1}, we compute the first derivative $b^{v}_1$ as
{\allowdisplaybreaks
\begin{align*}
    b^{v}_1 &= 
    \frac{\partial}{\partial A} \Bigg \lbrace
    \frac{(B_i - \Tilde{B}_i)}{\delta_i} y_i + \frac{\Tilde{B}_i x_i'\Tilde{\beta} - B_i x_i'\bar{\beta}} {\delta_i} \Bigg \rbrace \notag \\
    &= y_i \Bigg \lbrace \frac{1}{\delta_i} \frac{\partial B_i}{\partial A}  + (B_i - \tilde{B}_i) \frac{\partial}{\partial A} \Big (\frac{1}{\delta_i} \Big ) \Bigg \rbrace + \Tilde{B}_i x_i'\Tilde{\beta} \frac{\partial}{\partial A} \Big (\frac{1}{\delta_i} \Big ) - \Bigg \lbrace x_i'\bar{\beta} \frac{\partial B_i}{\partial A} + B_i \frac{\partial}{\partial A} (x_i'\bar{\beta}) \Bigg \rbrace \frac{1}{\delta_i} \notag \\ 
    & \qquad - B_i x_i'\bar{\beta} \frac{\partial}{\partial A} \Big (\frac{1}{\delta_i} \Big ) \notag \\
    &= \frac{(y_i - x_i'\bar{\beta}) }{\delta_i} \frac{\partial B_i }{\partial A} + \Big \lbrace (\Tilde{B}_i x_i'\Tilde{\beta} - B_i x_i'\bar{\beta}) + (B_i - \Tilde{B}_i) y_i \Big \rbrace \frac{\partial}{\partial A} \Big (\frac{1}{\delta_i} \Big ) - \frac{B_i}{\delta_i}\frac{\partial (x_i'\bar{\beta})}{\partial A} . 
\end{align*}}
Therefore, upto $O_p(1)$, for \eqref{bvhat1} we have 
{\allowdisplaybreaks
\begin{align*}
    \hat{b}^{v}_1 &= 
    \frac{\partial}{\partial A} b^v (A) |_{A=\hat{A}} \\ 
    &= \frac{(y_i - x_i'\hat{\beta}) }{\hat{\sigma}_i} \frac{\partial B_i }{\partial A} \Big |_{A=\hat{A}} + \lbrace (\Tilde{B}_i x_i'\Tilde{\beta} - \hat{B}_i x_i'\hat{\beta}) + (\hat{B}_i - \Tilde{B}_i) y_i \rbrace  \frac{\partial}{\partial A} \Big (\frac{1}{\sigma_i} \Big ) \Big |_{A=\hat{A}} \\
    &= \frac{(y_i - x_i'\hat{\beta}) }{\hat{\sigma}_i} \frac{\partial B_i }{\partial A} \Big |_{A=\hat{A}} \\ 
    &= (y_i - x_i'\hat{\beta}) \frac{\sqrt{\hat{A}+D_i}}{\sqrt{\hat{A}D_i}} \times \Bigg \lbrace \frac{- D_i}{(\hat{A}+D_i)^2} \Bigg \rbrace \\
    &= - (y_i - x_i'\hat{\beta}) \hat{\sigma}_i \frac{\hat{B}_i}{\hat{A}D_i},
\end{align*}}
since, $(\Tilde{B}_i x_i'\Tilde{\beta} - \hat{B}_i x_i'\hat{\beta}) = \Tilde{B}_i ( x_i'\Tilde{\beta} - x_i'\hat{\beta}) + (\Tilde{B}_i - \hat{B}_i) x_i'\hat{\beta} = O_p(m^{-1})$. We get this by expanding $x_i'\Tilde{\beta}$ and $x_i'\hat{\beta}$, respectively, around $A$ and subtracting, i.e., $x_i'\Tilde{\beta} - x_i'\hat{\beta} = O_p(m^{-1})$ and we already know that, $\Tilde{B}_i - \hat{B}_i = O_p(m^{-1})$,  $\Tilde{B}_i = O_p(1)$ and $x_i'\hat{\beta} = O_p(1)$. 

We note that 
\begin{align*}
     \frac{\partial B_i}{\partial A} = - \frac{D_i}{(A+D_i)^2} \quad \text{ and } \quad \frac{\partial}{\partial A} \Big (\frac{1}{\sigma_i} \Big) = - \frac{\sqrt{D_i}}{2A^{3/2} \sqrt{A+D_i}}.
\end{align*}
Then from \eqref{bv1}, we have upto the order $O(1)$
\begin{align*}
    b_2^v &= \frac{\partial b_1^v (A)}{\partial A} \\
    &= \frac{\partial}{\partial A} \Bigg [ (y_i - x_i'\bar{\beta}) \frac{1}{\sigma_i} \frac{\partial B_i }{\partial A} + \lbrace (\Tilde{B}_i x_i'\Tilde{\beta} - B_i x_i'\bar{\beta}) + (B_i - \Tilde{B}_i) y_i \rbrace \frac{\partial}{\partial A} \Big (\frac{1}{\sigma_i} \Big ) \Bigg ] \\
    &= \frac{\partial}{\partial A} \Bigg [ - (y_i - x_i'\bar{\beta}) \frac{\sqrt{A+D_i}}{\sqrt{AD_i}} \frac{D_i}{(A+D_i)^2} \\
    & \qquad \qquad - \Big \lbrace (\Tilde{B}_i x_i'\Tilde{\beta} - B_i x_i'\bar{\beta}) + (B_i - \Tilde{B}_i) y_i \Big \rbrace \frac{\sqrt{D_i}}{2A^{3/2} \sqrt{A+D_i}} \Bigg ] \\
    &= - y_i \sqrt{D_i} \frac{\partial}{\partial A}  \Big  \lbrace \frac{1}{\sqrt{A} (A+D_i)^{3/2}} \Big \rbrace + \frac{\sqrt{D_i}}{\sqrt{A} (A+D_i)^{3/2} } \frac{\partial}{\partial A} (x_i'\bar{\beta}) \\ 
    & \quad \quad  + \sqrt{D_i} x_i'\bar{\beta} \frac{\partial}{\partial A}  \Bigg \lbrace \frac{1}{\sqrt{A} (A+D_i)^{3/2}} \Bigg \rbrace 
    - \frac{\sqrt{D_i}}{2A^{3/2} \sqrt{A+D_i}} \frac{\partial}{\partial A} \Big \lbrace - B_i x_i'\bar{\beta} + B_i y_i \Big \rbrace \\
    & \quad \quad - \Big \lbrace (\Tilde{B}_i x_i'\Tilde{\beta} - B_i x_i'\bar{\beta}) + (B_i - \Tilde{B}_i) y_i \Big \rbrace \frac{\sqrt{D_i}}{2} \frac{\partial}{\partial A} \Big \lbrace \frac{1}{A^{3/2} \sqrt{A+D_i}} \Big \rbrace \\
    &= - \sqrt{D_i} (y_i  -  x_i'\bar{\beta}) \frac{\partial}{\partial A}  \Big  \lbrace \frac{1}{\sqrt{A} (A+D_i)^{3/2}} \Big \rbrace - \frac{\sqrt{D_i}}{2A^{3/2} \sqrt{A+D_i}} \Big \lbrace - x_i'\bar{\beta} \frac{\partial B_i}{\partial A}  + y_i \frac{\partial B_i}{\partial A}  \Big \rbrace \\
    & \quad \quad- \lbrace (\Tilde{B}_i x_i'\Tilde{\beta} - B_i x_i'\bar{\beta}) + (B_i - \Tilde{B}_i) y_i \rbrace \frac{\sqrt{D_i}}{2} \frac{\partial}{\partial A} \Big \lbrace \frac{1}{A^{3/2} \sqrt{A+D_i}} \Big \rbrace,
\end{align*}
ignoring the term $\frac{\partial}{\partial A} (x_i'\bar{\beta})$ that is of lower order $O(m^{-1/2})$.

We then calculate the following derivatives
\begin{align*}
    & \frac{\partial}{\partial A} \Big \lbrace \frac{1}{A^{3/2} \sqrt{A+D_i}} \Big \rbrace = -\frac{1}{2} \frac{(4A+3D_i)}{A^{5/2}(A+D_i)^{3/2}},\\
    & \frac{\partial}{\partial A} \Big \lbrace \frac{1}{\sqrt{A} (A+D_i)^{3/2} } \Big \rbrace = -\frac{(4A+D_i)}{2A(A+D_i)^{5/2}}
\end{align*}
and the fact that $(\Tilde{B}_i x_i'\Tilde{\beta} - \hat{B}_i x_i'\hat{\beta}) + (\hat{B}_i - \Tilde{B}_i) y_i = O_p(m^{-1})$. Finally we get $\hat{b}_2^v$, upto $O_p(1)$, as follows
\begin{align*}
    \hat{b}_2^v &= - \sqrt{D_i} (y_i - x_i'\hat{\beta}) \Bigg \lbrace - \frac{(4\hat{A}+D_i)}{2\hat{A}(\hat{A}+D_i)^{5/2}} \Bigg \rbrace - \frac{\sqrt{D_i}}{2\hat{A}^{3/2} \sqrt{ \hat{A}+D_i}} (y_i  -  x_i'\hat{\beta}) \Bigg \lbrace - \frac{D_i}{(\hat{A}+D_i)^2} \Bigg \rbrace \\
    &= (y_i - x_i'\hat{\beta}) \frac{\sqrt{D_i}}{2\hat{A}(\hat{A}+D_i)^{5/2}} \Bigg [ (4\hat{A}+D_i) + \frac{D_i}{\sqrt{\hat{A}}} \Bigg ].
\end{align*}
Thus, we have \eqref{b2vhat}, i.e, $\hat{b}_2^v = O_p(1)$.

\FloatBarrier
\section{Appendix : Additional simulation results}
\label{appC}

In this section, we provide additional results form the baseball analysis. We consider two other simulation setups mimicking the baseball data discussed in Section \ref{realdata}. In these simulations, we consider the model without covariates and simulate the direct estimates as follows
\paragraph{S3}: We simulate values {$ \hat{s}_i \sim \text{Bin}(45,p_i)$}, where $p_i$ is the known true season batting average of player $i$, and then we apply arc-sin transformation to get {$y_i = \sqrt{45}\arcsin(2\hat s_i/45 -1)$ and $\theta_i = \sqrt{45}\arcsin(2p_i -1)$} for all $i = 1, \cdots, 18$.
    \paragraph{S4}: We simulate $y_i$ directly from $N(\theta_i, 1)$, where {$\theta_i = \sqrt{45}\arcsin(2p_i -1)$} is the arc-sin transformed true batting average of player $i$ for all $ i = 1, \cdots, 18$.

{We subtract $\hat \mu = -3.275$ from $y_i$ and $\theta_i$ to construct a mean zero process.}
We provide EPC of $I_i^\text{YL}$ for setups S4 and S5 for 18 players in Table \ref{tab6}, which are again close to the nominal levels of $95\%$ and $90\%$.

\begin{table}[h]
\caption[Expected posterior coverage of $I_i^\text{YL}$ for {S3 and S4} for 18 baseball players]{EPC of $I_i^\text{YL}$ for setups S4 and S5 for 18 players.}
\label{tab6}
\footnotesize
\begin{center}
\begin{tabular}{|c|c|c|c|c|}
\hline
&  \multicolumn{2}{c|}{$\alpha = 0.05$}& \multicolumn{2}{c|}{$\alpha = 0.1$} \\
\hline
Player & S3 & S4 & S3 & S4 \\
\hline \hline
Clemente (Pitts, NL) & 94.7 & 94.7 & 89.0 & 89.2 \\ 
F. Robinson (Balt, AL) & 94.5 & 94.7 & 88.7 & 89.3 \\ 
Munson (NY, AL) & 94.2 &94.2 & 88.3 & 88.8 \\ 
Scott (Bos, AL) & 94.9 & 95.0 & 89.5 & 89.4 \\ 
F. Howard (Wash, AL) & 94.7 &94.6 & 89.0 & 88.9 \\ 
Campaneris (Oak, AL) & 94.7 & 94.4 & 89.4 & 88.9 \\ 
Spencer (Cal, AL) & 94.7 & 94.7 & 89.1 & 89.1 \\ 
Berry (Chi, AL) & 94.8 & 94.9 & 89.3 & 89.4 \\ 
Swoboda (NY, NL) & 94.7 &94.7 & 89.2 &89.4 \\ 
Kessinger (Chi, NL) & 94.8 &94.4 & 89.3 & 89.0\\
E. Rodriguez (KC, AL) & 94.7 & 94.9 & 89.1 & 89.7\\
Williams (Chi, AL) & 95.1 &94.8 & 89.8 & 89.4\\
Unser (Wash, AL) & 94.7 & 94.9 & 89.2 & 89.3\\
Johnstone (Cal, AL) & 94.6 & 94.9 & 89.0 & 89.7\\
Santo (Chi, NL) & 94.5 & 94.9 & 89.0 & 89.5\\
Petrocelli (Bos, AL) & 94.6 & 94.6 & 89.1 & 89.0\\
L. Alvarado (Bos, AL) & 94.7 & 94.4 & 89.3 & 88.8\\
Alvis (Mil, NY) & 94.1 & 93.9 & 88.2 & 87.9\\
\hline
\end{tabular}
\end{center}
\end{table}

\FloatBarrier
\section{Appendix : Additional BRFSS data analysis}
\label{appD}

\begin{figure}
    \centering
    \includegraphics[width=0.9\linewidth, page = 3]{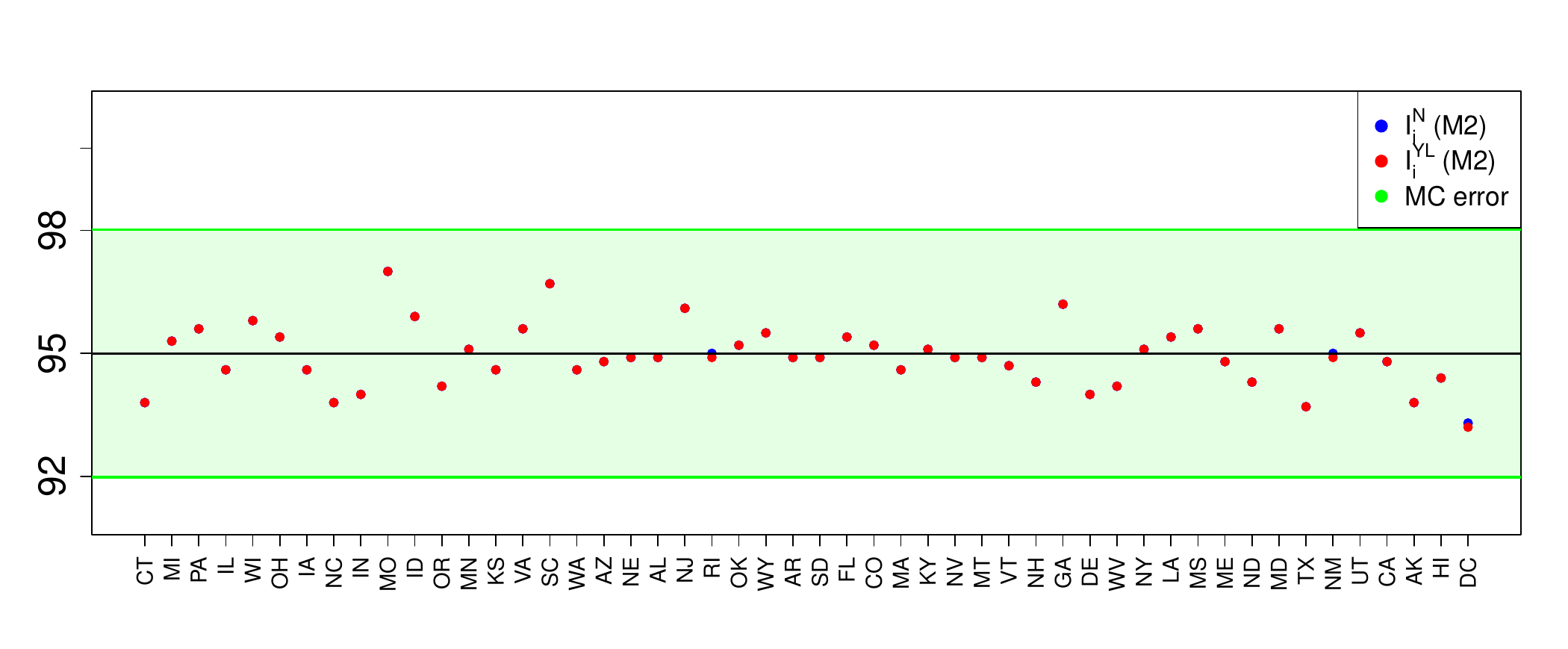}
    \caption{Plot of PC and MCE for two EBCIs: $I_i^\text{N}$ (M2) and $I_i^{\text{YL}}$ (M2) of percentage of people had flu shot within a year for 49 states and DC from the BRFSS data year 2024.}
    \label{brfss_pc_flu}
\end{figure}

\begin{figure}
    \centering
    \includegraphics[width=\linewidth, page = 4]{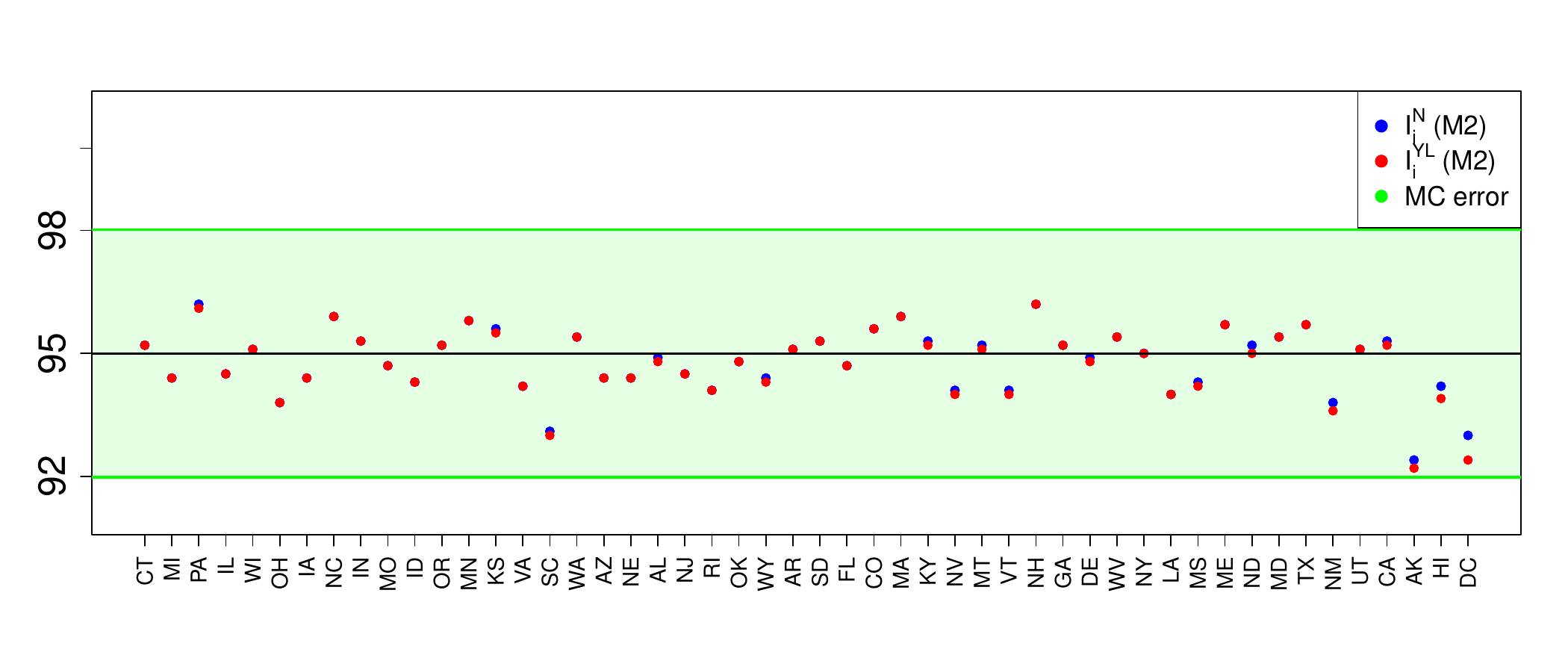}
    \caption{Plot of PC and MCE for two EBCIs: $I_i^\text{N}$ (M2) and $I_i^{\text{YL}}$ (M2) of percentage of people visited a doctor for regular checkup within a year for 49 states and DC from the BRFSS data year 2024.}
    \label{brfss_pc_chk}
\end{figure}

\begin{figure}
    \centering
    \includegraphics[width=\linewidth, page = 5]{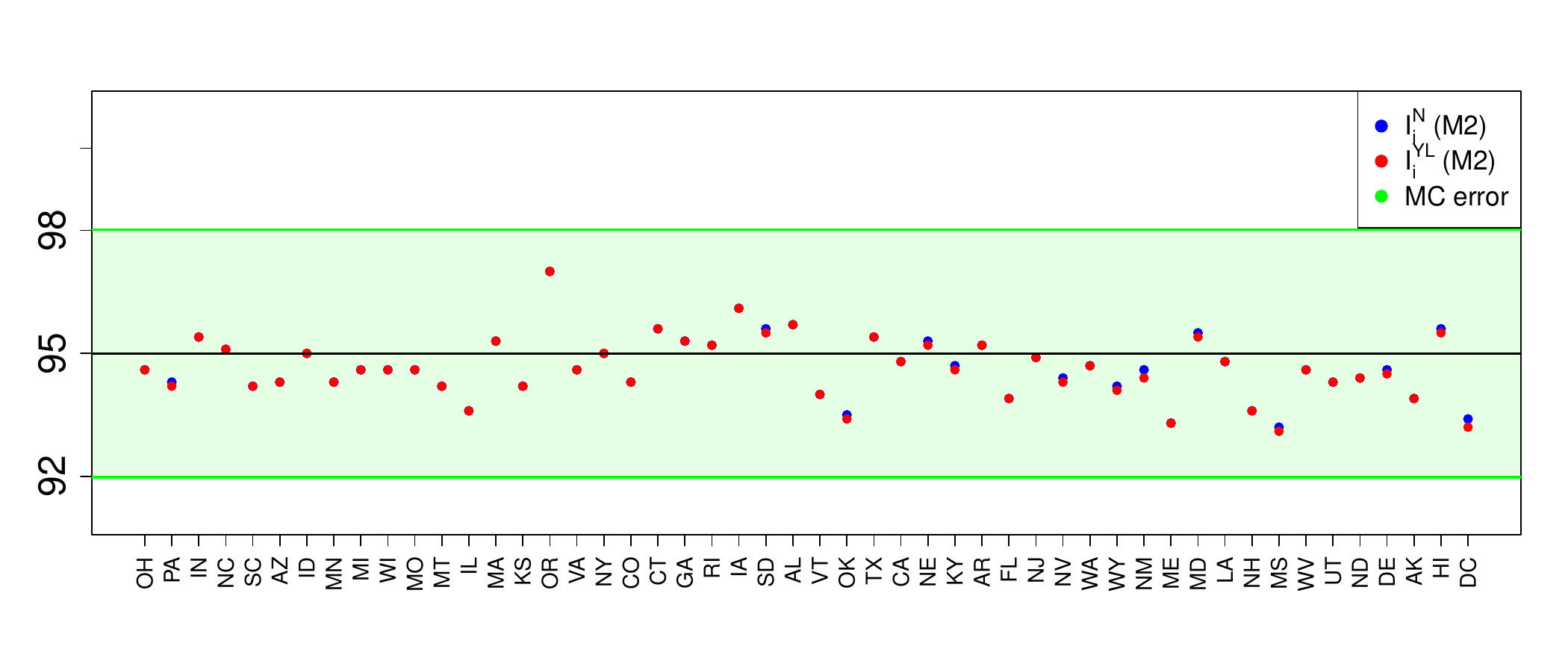}
    \caption{Plot of PC and MCE for two EBCIs: $I_i^\text{N}$ (M2) and $I_i^{\text{YL}}$ (M2) of percentage of people with depressive disorders for 49 states and DC from the BRFSS data year 2024.}
    \label{brfss_pc_dpr}
\end{figure}

\begin{figure}
    \centering
\includegraphics[width=\linewidth, page = 1]{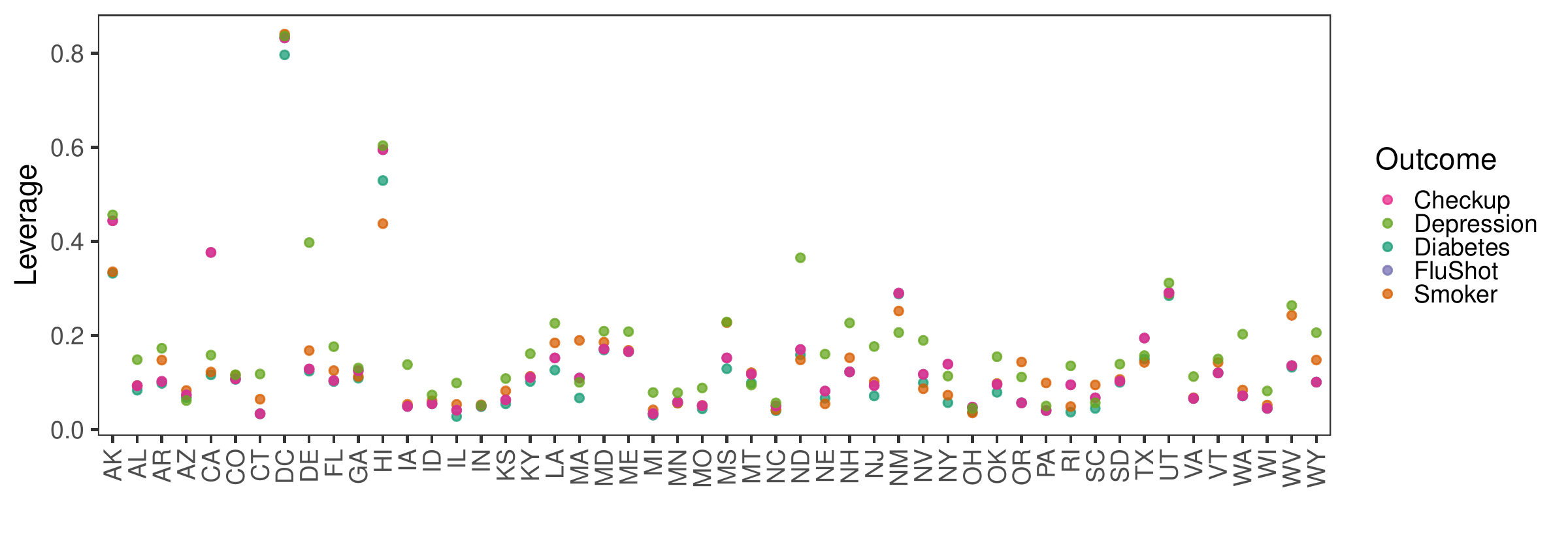}
    \caption{Leverage values from M2 for 49 states and DC from the ACS data years 2024.}
    \label{brfss_lev}
\end{figure}

In this section, we provide additional results from the BRFSS data analysis. Figures \ref{brfss_pc_flu}--\ref{brfss_pc_dpr} show PC of $I_i^{\text{N}}$ (M2) in blue colored dots and $I_i^{\text{YL}}$ (M2) in red dots and MCE in green bands, for the outcomes: \% of people who had flu shot within a year, \% of people who visited a doctor for a regular checkup within a year and \% of people with andepressive disorder, respectively. Figure \ref{brfss_lev} plots the leverage from M2 for all covariates from ACS data year 2024.


\end{document}